\documentclass[11pt]{article}
\usepackage{amsmath}
\usepackage{amsthm}
\input{epsf.tex}
\newtheorem{theorem}{Theorem}[section]
\newtheorem{proposition}[theorem]{Proposition}
\newtheorem{lemma}[theorem]{Lemma}
\newtheorem{corollary}[theorem]{Corollary}
\newtheorem{define}[theorem]{Definition}

\def\Empty{}

\catcode`\@=11

\def\section{\@startsection {section}{1}{\z@}{-3.5ex plus -1ex minus 
-.2ex}{2.3ex plus .2ex}{\large\bf}}

\def\fnum@figure{{\small Figure \thefigure}}
\def\fakefigure{\def\@captype{figure}}

\long\def\@makecaption#1#2{
    \vskip 10pt 
    \def\FCap{#2} \def\NoCap{\ignorespaces}
    \ifx \FCap\NoCap
       \setbox\@tempboxa\hbox{#1}  
      \else
       \setbox\@tempboxa\hbox{#1: \small \it #2}
    \fi
    \ifdim \wd\@tempboxa >\hsize   
        \unhbox\@tempboxa\par      
      \else                        
        \hbox to\hsize{\hfil\box\@tempboxa\hfil}  
    \fi}

\def\@oddhead{\hbox{}\rightmark \hfil \rm\thepage}
\def\sectionmark#1{\markright {\sc{\ifnum \c@secnumdepth >\z@
      \S\thesection.\hskip 1em\relax \fi #1}}}

\catcode`\@=12

\def\oplabel#1{
  \def\OpArg{#1} \ifx \OpArg\Empty {} \else
  	\label{#1}
  \fi}

\newlength{\saveu}

\input{serg.tex}
\input{leemacro.sty}

\def\centeredepsfbox#1{\centerline{\epsfbox{#1}}}

\begin{document}

\title{Ideal boundaries of pseudo-Anosov flows and
uniform convergence groups, with connections
and applications to large scale geometry}
\author{S\'{e}rgio R. Fenley
\thanks{Reseach partially supported by NSF grant 
DMS-0305313.}
\footnote{Mathematics Subject Classifiction.
Primary: 37D20, 53C23, 57R30, 37C85; Secondary: 
57M50, 37D50, 58D19.}
}
\maketitle

\vskip .2in

{\small{
\noindent
{\bf {Abstract}} $-$ 
Given a general pseudo-Anosov flow in a closed three manifold,
the orbit space of the lifted flow to the universal
cover is homeomorphic to an open disk. We construct a
natural compactification of this orbit
space with 
an ideal circle boundary.
If there are no perfect fits 
between stable and unstable leaves
and the flow is not topologically
conjugate to a suspension Anosov flow, we then show:
The ideal circle of the orbit space 
has a natural quotient space which is a sphere. This sphere 
is a dynamical systems 
ideal
boundary for a compactification of the universal 
cover of the manifold.
The main result is that the fundamental
group acts on the flow ideal boundary as
a uniform convergence group. Using a theorem of Bowditch, 
this yields a proof
that the fundamental group of the manifold is Gromov hyperbolic and it 
shows that the action of the fundamental group
on the flow ideal boundary is conjugate
to the action 
on the Gromov ideal
boundary.
This gives an entirely new proof that the fundamental group
of a closed, atoroidal $3$-manifold which fibers over the circle
is Gromov hyperbolic.
In addition with further geometric analysis, the main result also implies 
that pseudo-Anosov flows
without perfect fits are quasigeodesic flows
and that the stable/unstable foliations of these flows
are quasi-isometric foliations.
Finally we apply these results to (nonsingular) foliations:
if a foliation is $\rrrr$-covered or with one sided
branching in an aspherical, atoroidal three manifold then 
the results above imply that the leaves of the foliation
in the universal cover extend continuously to the
sphere at infinity.
}}

\vskip .2in
\section{Introduction}
The main purpose of this article is to analyse what information
can be obtained
about the asymptotic structure
or large scale geometry
of the universal cover of a manifold using only the 
dynamics of a pseudo-Anosov flow in the manifold.
We introduce a dynamical systems ideal boundary for a large
class of such flows and a corresponding compactification
of the universal cover.
The fundamental group acts
on the flow ideal boundary and compactification with excellent dynamical properties.
These objects are later shown to be strongly related to
the large scale geometry of the manifolds and of the
flows themselves. They also imply results about the geometry
of foliations.

In three manifold theory, the universal cover of the manifold
plays a crucial role. Topologically one is invariably interested
that the universal cover is $\rrrr^3$ \cite{Wa,He}. In terms of geometry,
for example, Thurston showed that a large class
of manifolds are hyperbolic \cite{Th1,Th2,Th3,Mor,
Ot1,Ot2} and the asymptotic or large scale
 structure of the universal cover was very
important for these results.

Our goal is to analyse what can a flow say about the asymptotic structure of the
universal cover of the manifold. Here we consider pseudo-Anosov flows
as they have rich dynamics and have been shown to be strongly connected
to the geometry \cite{Th3,Ot1} and topology of $3$-manifolds \cite{Ga-Oe,Fe6}.
Gabai and Oertel proved for example that the universal
cover of the underlying manifold
is $\rrrr^3$ \cite{Ga-Oe}. 
We will prove that under certain hypothesis the dynamics
of the flow creates a much richer asymptotic structure for 
the universal cover.

In this article all manifolds are connected.

We start by analysing the orbit space of the flow.
Suppose that $\Phi$ is a general pseudo-Anosov flow in
a closed $3$-manifold $M$.
Such flows
are very common 
\cite{Th4,Bl-Ca,Mo2,Mo3,Th5,Fe6,Cal2,Cal3}.
The flow has associated stable and unstable
(possibly singular) $2$-dimensional foliations $\ls, \lu$.  
When there are no singularities the flow is called an Anosov flow.
Let $\wwp$ be the lifted flow to the
universal cover $\widetilde M$
and let $\oo$ be the the orbit space of $\wwp$.
This orbit space is always homeomorphic 
to an open disk \cite{Fe1,Fe-Mo}.
The fundamental group of $M$ acting on $\mi$ by covering
translations, leaves invariant the foliation of $\mi$ by
flowlines of $\wwp$. Hence this induces an action of
the fundamental group on $\oo$.
The stable and unstable foliations of $\Phi$ lifted to
the universal cover also induce one dimensional foliations
in $\oo$.
We first analyse the asymptotic behavior of the orbit space:

\vskip .1in
\noindent
{\bf {Theorem A}} $-$ Let $\Phi$ be a pseudo-Anosov flow
in a closed $3$-manifold $M$.
There is a natural construction of a compactification
$\cd = \oo \cup \partial \oo$, obtained solely from
the stable and unstable foliations in $\oo$. The boundary
$\partial \oo$ is homeomorphic to a circle and the 
compactification $\cd$ is homeomorphic to a disk,
whose boundary circle is $\partial \oo$.
Since the fundamental group of $M$ preserves the stable and
unstable foliations in $\oo$, it follows that $\pi_1(M)$
acts by homeomorphisms on the compactification $\cd$ 
and also along the boundary circle $\partial \oo$.
\vskip .1in

We stress that compactifications of $\oo$ are not unique,
even compactifications to a closed disk.
For example given a point $p$ in  the above mentioned ideal boundary
$\partial \oo$, one can blow each point of the $\pi_1(M)$ orbit
of $p$ to a segment.
By doing this carefully
the ensuing compactification of 
$\oo$ is again a  closed disk where one can define a (non natural)
action of $\pi_1(M)$.

The stable/unstable foliations $\wls, \wlu$ in the universal
cover project to 1-dimensional  foliations
$\oos, \oou$ in $\oo$.
The possible singularities are only of $p$-prong type
with $p \geq 3$ (the condition $p \geq 3$ is necessary for
all the results in this article).
The prototype here is a suspension pseudo-Anosov flow over
a hyperbolic surface. In this case
$\oo$ is identified with a lift of a fiber
and it is possible to prove that 
the ideal circle boundary of $\oo$ constructed
in theorem A is identified with the circle
at infinity of the lift of the fiber.
In this example 
$\oos, \oou$ in $\oo$ correspond
to the stable and unstable foliations of the monodromy of the
fiber lifted to the universal cover of the fiber.
We stress that in
general there is no geometry (even coarse geometry) 
in the space $\oo$.

For general pseudo-Anosov flows, 
an ideal point of $\oo$
will be defined as an equivalence class of nested sequences 
of polygonal paths. A {\em polygonal path}
is a properly embedded, bi-infinite path in $\oo$ made up 
of a finite collection of segments alternatively
in $\oos, \oou$ and 2 rays of $\oos$ or $\oou$ at the ends.
In general one needs to use polygonal paths rather than 
just leaves of $\oos, \oou$ to define
ideal points of $\oo$ because of an obstruction which 
is called a perfect fit, as explained below.
Any ray of a leaf of $\oos, \oou$ is properly embedded in $\oo$ and
defines an ideal point of $\oo$, but 
there are many other points.
There is a natural group invariant topology in $\cd = \oo \cup 
\partial \oo$ $-$ this is a fundamental point here: the ideal points
($\partial \oo$)
and the topology in $\cd$ are constructed using
only the foliations $\oos, \oou$ in $\oo$.
Since these foliations are invariant under the action by the
fundamental group of $M$, then this group acts on $\cd$
by homeomorphisms.
The proof that $\dd$ is homeomorphic to a closed
disk is very involved and extremely long.
We show that $\partial \oo$ has a natural cyclic order and
that
$\partial \oo$ is metrizable, connected and more importantly
it is compact.
The last property is very hard to prove.
Point set topology theorems and additional work show that
$\partial \oo$ is homeomorphic to a circle and
$\cd$ is homeomorphic to a closed disk.
This works for any pseudo-Anosov flow.

We remark that Calegari and Dunfield
\cite{Ca-Du} previously showed that if $\Phi$ is a pseudo-Anosov
flow, then $\pi_1(M)$ acts nontrivially on a circle, with
very important consequences
for the existence question of pseudo-Anosov flows
\cite{Ca-Du}.
Their
construction is very different than ours.
They show that the space of ends of the leaf space
of say $\wls$ is
circularly ordered and maps injectively to a circle.
By collapsing complementary intervals
one gets an action on ${\bf S}^1$.
It is not entirely clear how to use the space of ends
in order to produce an actual compactification of $\oo$,
where the group acts naturally and with good properties.
For example, consider sequences escaping compact
sets in $\oo$ with all points in the
same stable leaf. As seen in the leaf space the points
do not go into any end, but they should have a convergent
subsequence in a compactification of $\oo$.
In this article we produce an actual compactification  
of the orbit space $\oo$ as a closed disk.
In addition very specific properties 
of the compactification
as related to the stable/unstable foliations ($\oos, \oou$)
will be used for the geometric results in the second
part of this article.

\vskip .1in
One main goal in introducing an ideal boundary for $\oo$ is
that it leads to an understanding of the asymptotic behavior
of $\mi$.
Our objective is 
to give a fairly explicit dynamical systems description of 
the asymptotic behavior of the universal
cover.
We do not know how to do this in general $-$ in this article
we can only deal with 
 pseudo-Anosov flows without
perfect fits. 

We first discuss perfect fits and their importance.
An unstable leaf $G$ of $\wlu$ makes a 
{\em perfect fit} with a stable leaf $F$ of $\wls$ if $G$ and $F$ do
not intersect but they ``almost" intersect:
any other unstable leaf sufficiently
near $G$ (and in the $F$ side), will
intersect $F$ and vice versa.
See detailed definition
in section 2 (figure \ref{loz}, a). 
We also use the terminology perfect fits for their
projections to the orbit space.
In the orbit space one can think of a perfect fit
as a proper embedding in $\oo$ of a rectangle 
minus a corner. 
Stable (unstable) leaves correspond to horizontal (vertical)
 segments.
The 2 boundary leaves without an endpoint form a perfect fit 
$-$ one stable leaf (horizontal) and one unstable leaf (vertical).
Perfect fits are very important in the topological
theory of pseudo-Anosov flows, see \cite{Ba1,Ba2,Fe1,Fe2,Fe4,Fe5}.
They occur for instance whenever there are closed orbits
of $\Phi$ which are freely homotopic \cite{Fe4,Fe5} or when the
leaf space of $\wls$ or $\wlu$ is not Hausdorff \cite{Fe4,Fe5}.
Examples of flows without perfect fits are
suspensions (with or without singularities)
and many other interesting examples as described later.

For the results in this article, perfect fits are one
main obstruction to simple definitions and proofs:
For example consider a point $p$ of $\partial \oo$ which is 
associated to the ideal point of
(say) an unstable ray $l$ of $\oou$.
Let $(z_n)_{n \in {\bf N}}$ be a nested sequence of
stable leaves intersecting $l$ and so that the intersection
with $l$ escapes compact sets in $l$.
What one strongly expects and hopes is that
the sequence $(z_n)_{n \in {\bf N}}$
defines the ideal point $p$ associated to $l$.
In particular one expects that the leaves $z_n$ escape compact sets
in $\oo$ as $n$ grows. This occurs in the suspension case
and in many other situations, but in fact it does not always happen.
When it does not
occur, then the sequence $(z_n)$ limits to 
a stable leaf $r'$ in $\oo$ and  one can then show that 
there is a stable leaf
$r$ (possibly $r = r'$), so that 
$r$ and $l$ form a perfect fit in $\oo$.
In this case the sequence $(z_n)$ will 
{\underline {not}} define the ideal point $p$.
Conversely any perfect fit generates a sequence $(z_n)$ as above.
Because of perfect fits
then to define ideal points of $\oo$,
one needs to consider not only leaves of $\oos, \oou$,
but rather
sequences of polygonal paths in $\oos, \oou$.
The definition of ideal points, implies that if
$r$ ray of $\oou$ and $l$ ray of $\oos$ form a perfect fit,
then these rays define the same ideal point of $\oo$.
Suspension Anosov flows (without singular orbits)
are special and have to be treated 
differently, because in that case a sequence of stable
leaves in $\oos$ escaping compact sets approaches 
infinitely many ideal points of $\oo$.

When there are no perfect fits
we construct the flow ideal boundary 
and compactification of $\mi$.
The flow ideal boundary is a quotient of
$\partial \oo$. The assumption of no perfect fits
is fundamental for this result:


\vskip .1in
\noindent
{\bf {Theorem B}} $-$ Let $\Phi$ be a pseudo-Anosov flow
without perfect fits and not topologically conjugate
to a suspension Anosov flow.
Let $\oo$ be its orbit space and $\partial \oo$ be the 
ideal boundary of theorem A.
Consider the equivalence relation in $\partial
\oo$ generated by: two points are in the same class if
they are ideal points of the same stable or unstable
leaf in $\oo$.
Let $\rr$ be the set of equivalence classes with
the quotient topology.
Then $\rr$ is homeomorphic to the 2-sphere. The fundamental group
of $M$ acts on $\rr$ by homeomorphisms.
There is a natural topology in $\mi \cup \rr$
making it into a compactification of $\mi$.
The action of $\pi_1(M)$ on $\mi$ extends to an action on $\mi \cup \rr$.
The quotient map from $\partial \oo \ (\cong {\bf S}^1)$ \  to  \
$\rr \ (\cong {\bf S}^2)$ \ is a group invariant
Peano curve associated to the flow $\Phi$.
All of this uses only the dynamics of the flow $\Phi$.
\vskip .1in

If $x$ in $\partial \oo$ is an ideal
point of (say) a stable leaf in $\oos$, 
then the condition of no perfect fits implies that 
no unstable leaf has ideal point $x$.
Hence if $k$ is the maximum number of prongs 
in singular leaves of $\oos$ (or $\oou$), then any equivalence
class has at most $k$ points.

\vskip .1in
Our goal is to relate the flow ideal compactification
with well known objects in three manifold topology.
We have actions of $\pi_1(M)$ on a circle \ ($\partial \oo$) \
and a sphere  \ ($\rr$).
Motivated by a lot of previous work in $2$ and $3$-dimensional
topology, one
asks whether such actions are convergence group actions.
For example a group that acts as a uniform convergence group
on the circle is topologically conjugate to
a Moebius group \cite{Tu1,Ga2,Ca-Ju} with fundamental consequences
for $3$-manifold theory \cite{Ga2,Ca-Ju}. 
Also a fundamental question of Cannon \cite{Ca-Sw} asks whether
a uniform convergence group acting on a 2-sphere is conjugate to
a cocompact Kleinian group. This is related to the geometrization
of $3$-manifolds.

A compactum is a compact Hausdorff space.
A group $\ga$ acts as a {\em convergence group} on a 
metrisable compactum $Z$ if
for any sequence $(\gamma_n)_{n \in {\bf N}}$ of
distinct elements in $\ga$, there is
a subsequence $(\gamma_{n_i})_{i \in {\bf N}}$
and a source/sink pair $y, x$ so that
$(\gamma_{n_i}(t))_{i \in {\bf N}}$ converges 
uniformly to the constant map with value  $x$ 
in compact sets of $Z - \{ y \}$
 \cite{Ge-Ma}.
Notice that $x, y$ may be the same point.
This is equivalent to 
$\ga$ acting properly discontinuously
on the set of distinct triples $\Theta_3(Z)$ of elements of
$Z$ \cite{Tu2,Bo2}.
In addition 
the action is {\em uniform} if the quotient
of $\Theta_3(Z)$ by the action is compact.
If $Z$ is perfect (no isolated points) then the additional condition is
equivalent to  every point  of $Z$ being a conical
limit point for the action.
A point $x$ in $Z$ is  a {\em conical limit point} if 
there is a sequence $(\gamma_n)_{n \in {\bf N}}$ 
in $\ga$ and $b, c$ distinct in 
$Z$, with $\gamma_n(x)$ converging to $c$
but for every other point $y$ in $Z$ then $(\gamma_n(y))$
converges to $b$.

The action of $\pi_1(M)$ on
$\partial \oo$ is not a convergence action.
Here is the proof: let $g$ non trivial in $\pi_1(M)$ so that 
$g$ fixes a point $x$ in $\oo$. Equivalently $g$ is associated
to a periodic orbit of $\Phi$. Up to taking a power assume
that $g$ leaves invariant all prongs of $\oos(x), \oou(x)$.
Hence it fixes the points in $\partial \oo$ which 
are the ideal points of these prongs. We show in this article
that all these ideal points are distinct points of the
circle. In addition the fixed points alternate between
contracting and expanding fixed points for $g$.
Now consider the sequence $(g^n)$ acting on $\partial \oo$.
The above facts imply that all elements in this
sequence of distinct elements
of $\pi_1(M)$ (or any subsequence) will 
share more than 2 fixed points and hence
the sequence $(g^n)$
does not have a single source/sink pair.
Hence the action of $\pi_1(M)$ on $\partial \oo$
is not a convergence group action.

\vskip .1in
\noindent
{\bf {Main theorem}} $-$ Let $\Phi$ be a pseudo-Anosov flow
without perfect fits and not topologically conjugate to a
suspension Anosov flow. Let $\rr$ be the associated flow ideal
boundary with corresponding compactification $\mi \cup \rr$ of
the universal cover.
Then the action of $\pi_1(M)$ on $\rr$ is a uniform
convergence group.
In addition the action of $\pi_1(M)$ on $\mi \cup \rr$ 
is a convergence group.

\vskip .1in

The main part of the proof is to prove uniform convergence action
on $\rr$.
Here 1-dimensional dynamics (action on the circle $\partial \oo$)
completely encodes the 2-dimensional dynamics (action on $\rr$).
A lot of the proof can be done using only this interplay
and the action on the $2$-dimensional space $\oo$,
but as expected the 3-dimensional setting of the
flow $\wwp$ in the
universal cover of $M$ needs to be used in some crucial
steps.

To prove the convergence group property,
we break into three cases up
to subsequences: \ 1) every $\gamma_n$ is associated to 
a singular orbit of $\Phi$,
\ 2) every $\gamma_n$ is associated to a nonsingular closed
orbit of $\Phi$, \ 3) every $\gamma_n$ acts freely on $\oo$.
For example consider case 2). 
Up to taking squares, the action of $\gamma_n$ in $\partial
\oo$ immediately has 4 fixed points, associated to the two ideal
points of the stable leaf of the periodic orbit and the two unstable ones.
By dynamics of pseudo-Anosov flows,
the stable points are locally attracting for the action
of $\gamma_n$ on $\partial \oo$
and the unstable ones are locally repelling. When there are no
perfect fits, this carries over to the whole of $\partial \oo$.
As the 2 ideal points of a stable leaf are identified in
$\rr$, this produces a source/sink behavior for 
(powers of) one $\gamma_n$.
An extended analysis shows the source/sink behavior for sequences.
The uniform property of the action is achieved by showing that 
every point of $\rr$ is a conical limit point.
The proofs of these results are very involved.

To prove the fact about the action on $\mi \cup \rr$,
consider a sequence of distinct
elements $(\gamma_n)_{n \in {\bf N}}$ of $\pi_1(M)$.
At this point we will already 
know that up to subsequence
it has a source/sink pair $y,x$ for the action
restricted to $\rr$.
We then show that $y,x$ is a source sink pair for the
action on $\mi \cup \rr$.
This depends on a careful analysis of neighborhoods
in $\mi \cup \rr$ of points in $\rr$.
The harder case is when such a point comes from an ideal
point of a leaf of $\oos$ or $\oou$.
The main theorem implies in particular that
the  action of $\pi_1(M)$ on $\rr$ (or 
on $\partial \oo$) is minimal. 

We mention that when 
there are perfect fits it is not at all clear 
what is the resulting structure of the quotient space $\rr$.
For example consider $\Phi$ an $\rrrr$-covered Anosov flow,
see \cite{Fe1}.
There are infinitely many examples where $M$ 
is hyperbolic \cite{Fe1}.
In this case the quotient $\rr$ (of the circle $\partial \oo$) 
as defined in theorem B, is
a union of a circle and two special points: each special point
is non separated from every point in the circle \cite{Fe1,Th5}.
Hence ${\cal R}$ is not even metrizable.
Clearly in this case the quotient $\rr$ does
not provide the expected ideal boundary of $\widetilde M$ (which
is a sphere).

\vskip .1in
This finishes the topological/dynamical systems
 part of the article. In the remainder of
the article 
we use the excellent properties of $\rr$ and $\mi \cup \rr$
to relate them with the large scale geometry of the manifold.
This has geometric consequences for the fundamental group
of the manifold and also for flows and foliations.
In particular we give an entirely new proof that 
the fundamental group of closed, atoroidal
$3$-manifolds that fiber over the circle
is Gromov hyperbolic.


The key tool will be the following:
Bowditch \cite{Bo1}, following ideas of Gromov,
 proved the very interesting theorem
that if $\ga$ acts as a uniform convergence group
on a perfect, metrisable
compactum $Z$, then $\ga$ is Gromov hyperbolic,
$Z$ is homeomorphic to the Gromov ideal boundary $\partial \ga$ and
the action on $Z$ is equivariantly topologically conjugate
to the action of $\ga$ on its Gromov ideal boundary.
This is a true geometrization theorem (in the sense of groups):
the hypothesis are entirely topological on the group action and
there is a strong geometric conclusion.
The main theorem then immediately implies the following:

\vskip .1in
\noindent
{\bf {Theorem D}} $-$ Let $\Phi$ be a pseudo-Anosov flow without
perfect fits and not topologically conjugate to a suspension Anosov flow.
Let $\rr$ be the associated flow ideal boundary of $\mi$ and
$\mi \cup \rr$ the flow ideal compactification. 
Then $\pi_1(M)$ is Gromov hyperbolic and the action of $\pi_1(M)$ on
$\rr$ is topologically conjugate to the action 
on the Gromov ideal boundary $\si$.
In addition the actions on $\mi \cup \rr$ and $\mi \cup \si$ are
also topologically conjugate $-$ by a homemorphism which is the
identity in $\mi$.
\vskip .1in

It was known that the Gromov boundary of 
$\pi_1(M)$ is a sphere because $M$ is irreducible
\cite{Be-Me}.
To prove the last statement of theorem D: Let $\xi$ be the
bijection between $\mi \cup \rr$ and $\mi \cup \si$, which is
the identity in $\mi$ and the conjugacy of the actions in $\rr$.
Clearly this is group equivariant.
We show that the bijection $\xi$ is continuous. This follows
from the convergence group action properties
for the action on $\mi \cup \rr$ plus the conjugacy 
between the actions on $\rr$ and $\si$.
Theorem D means that the constructions of
this article can be seen as  a dynamical systems analogue
to Gromov's geometric constructions in the case of 
this class of pseudo-Anosov flows.

A few remarks are in order here. In theorem D,
the result that $\pi_1(M)$ is Gromov hyperbolic is not new and 
also follows  from a result of Gabai-Kazez \cite{Ga-Ka} and
additional work.
The reason is: if $M$ with a pseudo-Anosov flow is toroidal,
then either there is a free homotopy between closed orbits
of the flow or the flow is topologically conjugate
to a suspension Anosov flow \cite{Fe7}.
The last option is disallowed by hypothesis of theorem D.
If there is a free homotopy between closed orbits
then there are perfect
fits \cite{Fe4,Fe5}.
Hence the hypothesis of theorem D imply that $M$ is atoroidal. 
With further analysis using the topological theory of pseudo-Anosov
flows \cite{Fe4,Fe5} one can then show that $\Phi$ has 
singular orbits.
Therefore the
(singular)
stable foliation blows up to an essential lamination
which is genuine, so \cite{Ga-Ka} 
implies that $\pi_1(M)$ is Gromov hyperbolic.
Gabai and Kazez showed that least area disks in $M$
satisfy a linear isoperimetric inequality. The
proof of this last
fact uses the ubiquity theorem for semi-Euclidean laminations
of Gabai \cite{Ga3}. This is a  deep but very mysterious result.
In particular it provides no 
direct relationship
with the Gromov ideal boundary.

The important new feature of theorem D is that 
it relates  the flow structure with the large scale 
geometric structure. 
Our construction gives
a very explicit description of the Gromov ideal boundary of
$\mi$ $-$ first as a purely dynamical systems object and 
{\underline {a posteriori}} implying that $\pi_1(M)$ is 
Gromov hyperbolic 
and totally  relating the two ideal boundaries.
In particular this is a new approach to  obtain
Gromov hyperbolicity.
There are several
important geometric consequences.
First we obain a new proof of a classical result:

\vskip .1in
\noindent
{\bf {Corollary E}} $-$ Let $\Phi$ be a suspension pseudo-Anosov flow
with at least a singular orbit in a closed $3$-manifold $M$.
Then $\pi_1(M)$ is
Gromov hyperbolic.
\vskip .1in

This theorem has two well known proofs: the original by Thurston
\cite{Th3} and a later proof by Bestvina and Feighn \cite{Be-Fe}.
Thurston's original proof uses quasiconformal maps, Kleinian groups
and the double limit theorem and obviously proves much more $-$
it proves that $M$ admits a hyperbolic metric.
Bestvina and Feighn's proof is a geometric group theory proof
and introduces the extremely useful condition of flaring annuli.
Our proof is entirely new in the sense that it uses dynamical systems
and convergence groups via Bowditch's theorem.

The proof of corollary E is as follows: 
Let $S$ be a cross section of $\Phi$. Since there is a singularity
of $\Phi$, 
$S$ is a hyperbolic surface. We already mentioned that
the orbit space of $\wwp$ is identified with the universal
cover $\widetilde S$ and the foliations $\oos, \oou$ in $\oo$
are identified with lifts $\widetilde f^s, \widetilde f^u$
of the stable and unstable foliations of the monodromy of
the fibration.
According to corollary D all that is needed is to prove that
there are no perfect fits. Notice that this is a topological 
condition. We will check this for $\widetilde f^s, \widetilde f^u$.
Consider $S$ with a hyperbolic metric, hence $\widetilde S$ is
the hyperbolic plane. If there is a perfect fit between
$\widetilde f^s$ and $\widetilde f^u$, then there 
is a ray $l$ of (say) $\widetilde f^s$ so that: if $s_n$ is a sequence 
of unstable leaves (of $\widetilde f^u$) intersecting $l$ and
with $l \cap s_n$ escaping to the appropriate
end of $l$ then $s_n$ does not
escape compact sets in $\widetilde S$ and converges to a
leaf $s$ of $\widetilde f^u$. Now use the fundamental
property that leaves of $\widetilde f^s, \widetilde f^u$ are
uniform quasigeodesics in $\widetilde S$ \cite{Th4,FLP}.
It follows that $s$ is
unique and that $l, s$ have a common ideal point in $\partial \widetilde S$.
This is impossible \cite{Th4,FLP}.
This finishes the proof of corollary E. 
As a remark for future reference, the case of pseudo-Anosov flows
without perfect fits shares many features with the suspension pseudo-Anosov
situation: the property alluded above about ideal points of $l$ and $s$
has an analogue for general pseudo-Anosov flows without perfect
fits. This is the content of the escape lemma (lemma \ref{duo}).
The escape lemma is extremely useful for the analysis of pseudo-Anosov 
flows without perfect fits.

 
\vskip .1in
We should remark that if $M$ is closed, irreducible, aspherical, atoroidal and
with infinite fundamental group  then
Perelman's results \cite{Pe1,Pe2,Pe3}
show that $M$ is hyperbolic.
We do not make use of Perelman's results here.
We stress again that a fundamental goal of this article is to analyse which
geometric information can be obtained solely
from dynamical systems constructions.

%

We now describe other very  important geometric consequences of 
theorem D.
Flow objects (flowlines, stable/unstable leaves, foliations transverse
to the flow) behave very well in the compactification $\mi \cup \rr$.
Since this is homeomorphic to the Gromov compactification, 
%
it is natural to expect that
these objects also have good geometric properties.
First we study 
metric properties of such flows and their
stable/unstable foliations.
In manifolds with Gromov hyperbolic fundamental group the relation
between objects in $\mi$ and their limit sets is extremely
important \cite{Th1,Th2,Th3,Gr,Gh-Ha,CDP} and is related to
the large scale geometry in $\mi$.
A flow in $M$ is {\em quasigeodesic} if flow lines in $\mi$ are
uniformly efficient in measuring ambient distance up
to a bounded multiplicative distortion \cite{Th1,Gr,Gh-Ha,CDP}.
It implies that each flow line is a bounded distance
from the corresponding geodesic which has the same
ideal points. Quasigeodesic flows are very useful
\cite{Ca-Th,Mo1,Mo2,Fe2}.
Usually it is very hard to show that a flow is quasigeodesic
and there is no general construction of quasigeodesic flows
in hyperbolic manifolds $-$ the known class of examples
is relatively small.
Theorem D provides a powerful way to obtain
quasigeodesic flows:

\vskip .1in
\noindent
{\bf {Theorem F}} $-$ Let $\Phi$ be a pseudo-Anosov flow
without perfect fits.
Then $\Phi$ is a quasigeodesic flow in $M$.
In addition $\Lambda^s, \Lambda^u$ are quasi-isometric
singular foliations in $M$.
\vskip .1in

First assume that $\Phi$ is not topologically conjugate
to a suspension Anosov flow.
By theorem D, $\pi_1(M)$ is Gromov hyperbolic.
To prove theorem F we first prove some properties 
in the flow compactification $\mi \cup \rr$:
\ 1) Each flow line $\gamma$ of $\wwp$ has a unique
forward ideal point $\gamma_+$ in $\rr$ and 
a backward ideal point $\gamma_-$;
\ 2) For each $\gamma$ the points $\gamma_-, \gamma_+$
are distinct; \ 3) The forward (backward) ideal point map
is continuous.
Theorem D conjugates the action in $\mi \cup \rr$
to the action in $\mi \cup \si$, hence the
same properties are true in $\mi \cup \si$.
A previous result of the author and Mosher \cite{Fe-Mo} 
then implies that $\Phi$
is quasigeodesic.

{\em Quasi-isometric} behavior for
$\Lambda^s, \Lambda^u$ means that leaves of $\wls$ (or $\wlu$)
are uniformly efficient in measuring distance in $\mi$ 
\cite{Th1,Gr,CDP}. This is
the analogue of quasigeodesic behavior in the two dimensional setting
and again it is extremely useful \cite{Gr,Th1,Th2,Th3}.
For example it implies that leaves of $\wls, \wlu$ are
quasiconvex \cite{Th1,Gr}.
Quasi-isometric foliations are very useful
\cite{Ca-Th,Th5,Fe5,Fe8}.
To prove the second part of theorem F: 
the lack of  perfect fits implies that the leaf spaces of
$\wls, \wlu$
are Hausdorff \cite{Fe4,Fe5}.
Together with the fact that $\Phi$ is quasigeodesic
this now implies that $\Lambda^s, \Lambda^u$ are
quasi-isometric foliations \cite{Fe5}.
This provides a new way to obtain quasi-isometric 
singular foliations in such manifolds.

If now $\Phi$ is topologically conjugate to a suspension
Anosov flow, then quasigeodesic behavior of $\Phi$ and
quasi-isometric behavior of $\ls, \lu$ are easy to prove.

\vskip .1in
Finally we apply these results to (nonsingular) foliations and their asymptotic
properties and we show that pseudo-Anosov flows
without perfect fits are very common.
A foliation $\fol$ in a $3$-manifold is $\rrrr$-covered if the leaf
space $\hp$ of $\fn$ is Hausdorff or equivalently homeomorphic
to the real numbers.
$\rrrr$-covered foliations are very common
\cite{He,Fe1,Th5,Cal1}.
On the other hand 
if $\hp$ is not Hausdorff then it is a simply connected,
non Hausdorff, 1-dimensional manifold with a countable
basis \cite{Ba3}.
Hence it is orientable. The non separated points in $\hp$
correspond to branching in the negative (positive)
direction if they are separated on their
positive (negative) sides.
A foliation $\fol$ has one sided branching if the branching
in $\fn$ is only in one direction (positive or negative).

If $\fol$ is a Reebless foliation in $M^3$ aspherical with 
$\pi_1(M)$ Gromov hyperbolic then 
each leaf $F$ of $\fn$ is uniformly Gromov hyperbolic in its path metric
and has an ideal circle $\pin F$ compactifying 
it to a closed disk $F \cup \pin F$.
The {\em continuous
extension question} asks what is the asymptotic behavior
of the leaves of $\fn$, that is, do they approach
the ideal boundary $\si$ in a continuous way?
This is formulated as follows: 
Does the inclusion $i: F \rightarrow \mi$
extend continuously to $i: F \cup \pin F \rightarrow \mi \cup \si$?
If so then $i$ restricted to $\pin F$ is a continuous
parametrization of the limit set of $F$, which will be locally
connected. When this happens
for all leaves of $\fn$, we say that $\fol$ has
the continuous extension property \cite{Ga1,Ca-Th,Fe5}.
This property is very hard to prove.

We use the geometric tools developed in this article to
prove the following theorem.
For any codimension one foliation $\fol$ if it is not
transversely orientable there is a transversely orientable
lift $\fol_2$ in a double cover $M_2$ of $M$. 
If $\fol$ is transversely orientable we abuse notation
and let $M_2 = M$ and $\fol_2 = \fol$.
If $M$ is aspherical and atoroidal then the author \cite{Fe6}
and Calegari \cite{Cal2} proved that there is a pseudo-Anosov
flow $\Phi$ which is transverse to $\fol_2$ in $M_2$.

\vskip .1in
\noindent
{\bf {Theorem G}} $-$  Let $\fol$ be an
$\rrrr$-covered foliation in an aspherical, atoroidal $3$-manifold $M$.
The pseudo-Anosov flow $\Phi$ transverse to the transversely oriented
foliation $\fol_2$ associated to $\fol$ does not have
perfect fits and is not conjugate to a suspension
Anosov flow. It follows that $\Phi$ is quasigeodesic by Theorem F and
this in turn implies that $\fol_2$ satisfies the continuous extension
property.
This trivially implies that $\fol$ satisfies the continuous extension property.
In addition the stable/unstable foliations of $\Phi$ (in the cover $M_2$)
are quasi-isometric.
\vskip .1in

The aspherical property is used only to get rid of a manifold which is
finitely covered by ${\bf S}^2 \times {\bf S}^1$. 
The problem is that the $\rrrr$-covered
property does not imply that the foliation is Reebless. 
For example consider the foliation $\fol$ of ${\bf S}^2 \times {\bf S}^1$
which is obtained by glueing two Reeb components appropriately. If one
is careful, then $\fol$ is $\rrrr$-covered. On the other hand
the author previously proved that if $\fol$ is $\rrrr$-covered, but
not Reebless, then $M$ is finitely covered by ${\bf S}^2 \times {\bf S}^1$
\cite{Fe8}.
Apart from this special case, the universal cover is homeomorphic
to $\rrrr^3$ and the results of the author and Calegari can be applied.

The continuous extension property was previously proved for:
 1) fibrations
in the seminal work of Cannon and Thurston \cite{Ca-Th}, 2) Finite
depth foliations and some other classes by the author 
\cite{Fe5,Fe8},
3) slitherings or uniform foliations by Thurston \cite{Th5}.
The methods of the proof were very different from
those in this article $-$ in all of the previous cases one always had 
a strong geometric property to start with:
For example in the case of finite depth foliations (not fibrations),
the compact leaf is quasi-isometrically embedded and therefore
quasiconvex.
After 
some work this implies that the almost transverse pseudo-Anosov
flow is 
quasigeodesic. After substantial more 
work this implies the continuous extension
property for the foliation.
The problem in general is that for instance in an arbitrary 
$\rrrr$-covered foliation, the leaves have no good geometric
property to start with $-$ so these methods do not work.
In this article we obtain geometric properties for the flow
directly and solely from the dynamics of the pseudo-Anosov flow
and this can then be applied to the foliations.
Theorem G implies the previous results for fibrations and slitherings.
Theorem G produces new examples of
quasigeodesic flows and quasi-isometric foliations.

In order to prove theorem G assume that $\fol$ is
transversely oriented and  start with a pseudo-Anosov 
flow $\Phi$ transverse to $\fol$ as constructed in
\cite{Fe6} or \cite{Cal2}. We show that 
$\Phi$ is not conjugate to a suspension Anosov flow
and has no perfect fits.
By theorem F, the flow $\Phi$ is quasigeodesic and
its stable/unstable foliations are quasi-isometric. 
By previous results \cite{Fe8}, it follows that $\fol$ has the continuous 
extension property.

We also consider foliations with one sided branching and prove:

\vskip .1in
\noindent
{\bf {Theorem H}} $-$ Let $\fol$ be a foliation with 
one sided branching in $M^3$ aspherical, atoroidal.
Then $\fol$ is transverse to a pseudo-Anosov flow $\Phi$ without perfect
fits and not conjugate to a suspension Anosov flow.
It follows that $\Phi$ is quasigeodesic, its stable/unstable foliations
are quasi-isometric  and
$\fol$ has the continuous extension property.
If $F$ is a leaf of $\fn$, then the limit set of $F$  is
not the whole sphere.
\vskip .1in

Under the conditions of this theorem,
Calegari \cite{Cal3} proved that
$\fol$ is transverse to a pseudo-Anosov flow $\Phi$.
We show that such 
$\Phi$ does not have perfect fits nor is conjugate
to a suspension Anosov flow.
By theorem F, the flow $\Phi$ is quasigeodesic.
This implies that $\fol$ has the
continuous extension property.
The last statement follows from metric properties
of leaves of $\wls, \wlu$.

The geometric applications obtained here (theorems E, F, G and H) were the main
motivation for the construction of the flow ideal
boundary of $\mi$ and the ideal circle of $\oo$.

The open case for the continuous extension question
is contained in the case when $\fol$ branches in both directions.
The case of finite depth foliations 
was resolved very recently
in \cite{Fe8} using work of Mosher, Gabai and the author
\cite{Mo3,Fe-Mo}.
For general foliations with two
sided branching, Calegari \cite{Cal4} constructed a 
very full lamination transverse to $\fol$, like 
the stable/unstable foliation of a flow.
It is possible that in certain situations 
there are 2 laminations, which perhaps are transverse
to each other 
and these can be possibly blowed down to produce
a pseudo-Anosov
flow transverse or almost transverse to $\fol$ \cite{Mo3}.
When the ideal dynamics of the case of a pseudo-Anosov flow
with perfect fits is better understood, then Calegari's results
could be very useful.


The geometric properties of flows and
foliations (theorems F, G and H) are proved at the 
end of the article, in sections 6 and 7. 
The proofs use the main theorem, theorem D and previous results.
Theorems G, H provide a large class of examples of 
pseudo-Anosov flows without perfect fits and also quasigeodesic flows
and quasi-isometric foliations.
The bulk of the article is proving theorem A (section 3), 
theorem B and the main theorem (section 4). 
Gromov hyperbolicity and conjugacy
are proved in section 5.


\vskip .1in
\noindent
{\bf {How to read this article}} $-$ The body of the article has
two main parts: 1) Section 3 $-$ ideal boundary of $\oo$, 2) Section 
4 $-$ flow ideal boundary for flows without perfect fits
and uniform convergence group action.
For those mainly interested in the geometric results (sections 4-7)
we highlight in section 3 where the case without perfect fits
has simplified proofs.

We thank Lee Mosher who told us about Bowditch's theorem.
We also thank the reviewer who did an outstanding job
of very carefully checking the whole article and who had 
inumerable useful comments,  many detailed suggestions 
and corrections which were incorporated 
in this article.

\section{Preliminaries: Pseudo-Anosov flows}
\label{prelim}


Given $M$ let $\widetilde M \rightarrow M$ be a fixed universal cover.

Let $\Phi$ be a flow on a closed
3-manifold $M$. 
We say that $\Phi$ is a {\em pseudo-Anosov flow} if the following are
satisfied:


- For each $x \in M$, the flow line $t \to \Phi(x,t)$ is $C^1$,
it is not a single point,
and the tangent vector bundle $D_t \Phi$ is $C^0$.

- There is a finite number of periodic orbits $\{ \gamma_i \}$,
called {\em singular orbits}, such that the flow is 
``topologically" smooth off of the
singular orbits (see below).

- The flowlines of $\Phi$ are contained in 
two possibly singular $2$-dimensional foliations
$\ls, \lu$ satisfying: Outside of the singular orbits, the
foliations $\ls, \lu$ are not singular, they 
are transverse to each other  and their leaves
intersect exactly along the orbits of $\Phi$.
A leaf containing a singularity is homeomorphic 
to $P \times I/f$
where $P$ is a $p$-prong in the plane and $f$ is a homeomorphism
from $P \times \{ 1 \}$ to $P \times \{ 0 \}$.
We restrict to $p$ at least $2$, that is, we do not allow
$1$-prongs.

- In a stable leaf all orbits are forward asymptotic,
in an unstable leaf all orbits are backwards asymptotic.

Basic references for pseudo-Anosov flows are \cite{Mo1,Mo3}.

\vskip .05in
The singular
foliations lifted to $\mi$ are
denoted by $\wls, \wlu$.
If $x$ is a point in $M$ let $W^s(x)$ denote the leaf of $\ls$ containing
$x$.  Similarly one defines $W^u(x)$
and in the
universal cover $\ws(x), \wu(x)$.
If $\alpha$ is an orbit of $\Phi$, similarly define
$W^s(\alpha)$, 
$W^u(\alpha)$, etc...
Let also $\wwp$ be the lifted flow to $\mi$.

\vskip .1in
We review the results about the topology of
$\wls, \wlu$ that we will need.
We refer to \cite{Fe4,Fe5} for detailed definitions, explanations and 
proofs.
Proposition 4.2 of \cite{Fe-Mo} shows that the
orbit space of $\wwp$ in
$\mi$ is homeomorphic to the plane $\rrrr^2$.
This orbit space is denoted by $\oo \cong \mi/\wwp$. 
Let $\Theta: \mi \rightarrow \oo \cong \rrrr^2$
be the projection map. 
If $L$ is a 
leaf of $\wls$ or $\wlu$,
then $\Theta(L) \subset \oo$ is a tree which is either homeomorphic
to $\rrrr$ if $L$ is regular,
or is a union of $k$ rays all with the same starting point
if $L$ has a singular $k$-prong orbit.
The foliations $\wls, \wlu$ induce singular $1$-dimensional foliations
$\oos, \oou$ in $\oo$. Its leaves are the $\Theta(L)$'s as
above.
If $L$ is a leaf of $\wls$ or $\wlu$, then 
a {\em sector} is a component of $\mi - L$.
Similarly for $\oos, \oou$. 
If $B$ is any subset of $\oo$, we denote by $B \times \rrrr$
the set $\Theta^{-1}(B)$.
The same notation $B \times \rrrr$ will be used for
any subset $B$ of $\mi$: it will just be the union
of all flow lines through points of $B$.
If $x$ is a point of $\oo$, then $\oos(x)$ (resp. $\oou(x)$)
is the leaf of $\oos$ (resp. $\oou$) containing $x$.

\begin{define}
Let $L$ be a leaf of $\wls$ or $\wlu$. A slice leaf of $L$ is 
$l \times \rrrr$ where $l$ is a properly embedded
copy of the real line in $\Theta(L)$. For instance if $L$
is regular then $L$ is its only slice leaf. If a slice leaf
is the boundary of a sector of $L$ then it is called
a line leaf of $L$.
If $a$ is a ray in $\Theta(L)$ then $A = a \times \rrrr$
is called a half leaf of $L$.
If $\zeta$ is an open segment in $\Theta(L)$ 
it defines a flow band $L_1$ of $L$
by $L_1 = \zeta \times \rrrr$.
\end{define}

\noindent
{\bf {Important convention}} $-$ In general a slice leaf is just a slice
leaf of some $L$ in $\wls$ or $\wlu$ and so on.
We also use the terms slice leaves, line leaves,
perfect fits, lozenges and rectangles for the projections of these
objects in $\mi$ to the orbit space $\oo$.

\vskip .1in
If $F \in \wls$ and $G \in \wlu$ 
then $F$ and $G$ intersect in at most one
orbit. 
Also suppose that a leaf $F \in \wls$ intersects two leaves
$G, H \in \wlu$ and so does $L \in \wls$.
Then $F, L, G, H$ form a {\em rectangle} in $\mi$
and there is  no singularity  of $\wwp$ 
in the interior of the rectangle see \cite{Fe4} pages 637-638.
There will be two generalizations of rectangles: 1) perfect fits $=$ 
in the orbit space this is a  properly embedded rectangle
with one corner removed 
and 2) lozenges $=$ rectangle with two opposite corners removed.

\begin{define}{(\cite{Fe2,Fe4})}{}
Perfect fits -
Two leaves $F \in \wls$ and $G \in \wlu$, form
a perfect fit if $F \cap G = \emptyset$ and there
are half leaves $F_1$ of $F$ and $G_1$ of $G$ 
and also flow bands $L_1 \subset L \in \wls$ and
$H_1 \subset H \in \wlu$,
so that 
%
the set 

$$\overline F_1 \cup \overline H_1 \cup 
\overline L_1 \cup \overline G_1$$

\noindent
separates $M$ 
and the joint structure of $\wls, \wlu$ in a complementary
component $R$ is that of
a rectangle as above without one corner orbit.
Specifically, a stable leaf intersects $H_1$ if and only
if it intersects $G_1$ and similarly for unstable
leaves intersecting $F_1, L_1$.

\end{define}

\begin{figure}
\centeredepsfbox{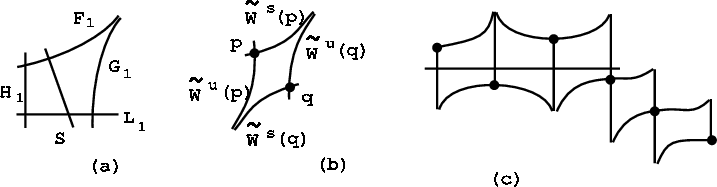}
\caption{a. Perfect fits in $\mi$,
b. A lozenge, c. A chain of lozenges.}
\label{loz}
\end{figure}

We refer to fig. \ref{loz}, a for perfect fits.
We also say that the leaves $F, G$ {\em almost intersect}.

%
%
%
%

\begin{define}{(\cite{Fe2,Fe4})}{}
Lozenges - A lozenge is an open region of $\mi$ whose closure
in $\mi$
is homeomorphic to a rectangle with two corners removed.
More specifically two orbits $\alpha = \wwp_{\rrrr}(p), 
\beta = \wwp_{\rrrr}(q)$ 
form the corners
of a lozenge if there are half leaves $A, B$ of
$\ws(\alpha), \wu(\alpha)$ defined by $\alpha$
and  $C, D$ half leaves of $\ws(\beta), \wu(\beta)$ so
that $A$ and $D$ form a perfect fit and so do
$B$ and $C$. The region in $\mi$ bounded by $A, B, C, D$ is
the lozenge
$R$ and it does not have any singularities.
%
%
%
%
See fig. \ref{loz}, b.
\end{define}


This is definition 4.4 of \cite{Fe5}.
The sets $A, B, C, D$ are the sides of the lozenge.
There may be singular orbits
on the sides of the lozenge and the corner orbits.
Two lozenges are {\em adjacent} if they share a corner and
there is a stable or unstable leaf
intersecting both of the lozenges, see fig. \ref{loz}, c.
Therefore they share a side.
A {\em chain of lozenges} is a collection $\{ \cc _i \}, 
i \in I$, of lozenges
where $I$ is an interval (finite or not) in ${\bf Z}$,
so that if $i, i+1 \in I$, then 
${\cal C}_i$ and ${\cal C}_{i+1}$ share
a corner, see fig. \ref{loz}, c.
Consecutive lozenges may be adjacent or not.
The chain is finite if $I$ is finite.

\begin{define}{}{}
Suppose $A$ is a flow band in a leaf of $\wls$.
Suppose that for each orbit $\gamma$ of $\wwp$ in $A$ there is a
half leaf $B_{\gamma}$ of $\wu(\gamma)$ defined by $\gamma$ so that: 
for any two orbits $\gamma, \beta$ in $A$ then
a stable leaf intersects $B_{\beta}$ if and only if 
it intersects $B_{\gamma}$.
%
%
This defines a stable product region which is the union
of the $B_{\gamma}$.
Similarly define unstable product regions.
\label{defsta}
\end{define}

The main property of product regions is the following, see
\cite{Fe5} page 641:
for any $F \in \wls$, $G \in \wlu$ so that 
$(i) \ F \cap A  \
\not = \ \emptyset \ \ {\rm  and} \ \ 
 (ii) \ G \cap A \ \not = \ \emptyset,
\ \ \ {\rm then} \ \ 
F \cap G \ \not = \ \emptyset$.
There are no singular orbits of 
$\wwp$ in $A$.

%

We abuse convention and say that
a leaf $L$ of $\wls$ or $\wlu$ is {\em periodic}
if there is a non trivial covering translation
$g$ of $\mi$ with $g(L) = L$. This is equivalent
to $\pi(L)$ containing a periodic orbit of $\Phi$, which may or
may not be singular.
In the same way, an orbit 
$\gamma$ of $\wwp$
is {\em periodic} if $\pi(\gamma)$ is a periodic orbit
of $\Phi$.
Finally a leaf $l$ of $\oos$ or $\oou$ is periodic if there
is $g \not =$ id in $\pi_1(M)$ with $g(l) = l$.

We say that 
two orbits $\gamma, \alpha$ of $\wwp$ 
(or the leaves $\ws(\gamma), \ws(\alpha)$)
are connected by a 
chain of lozenges $\{ {\cal C}_i \}, 1 \leq i \leq n$,
if $\gamma$ is a corner of ${\cal C}_1$ and $\alpha$ 
is a corner of ${\cal C}_n$.
If a lozenge ${\cal C}$ has corners $\beta, \gamma$ and
if $g$ in $\pi_1(M) - id$ satisfies $g(\beta) = \beta$,
$g(\gamma) = \gamma$ (and so $g({\cal C}) = {\cal C}$),
then $\pi(\beta), \pi(\gamma)$ are closed orbits
of $\Phi$ which are freely homotopic to the inverse of each
other.

\begin{theorem}{(\cite{Fe5}, theorem 4.8)}{}
Let $\Phi$ be a pseudo-Anosov flow in $M$ closed and let 
$F_0 \not = F_1 \in \wls$.
Suppose that there is a non trivial covering translation $g$
with $g(F_i) = F_i, i = 0,1$.
Let $\alpha_i, i = 0,1$ be the periodic orbits of $\wwp$
in $F_i$ so that $g(\alpha_i) = \alpha_i$.
Then $\alpha_0$ and $\alpha_1$ are connected
by a finite chain of lozenges 
$\{ {\cal C}_i \}, 1 \leq i \leq n$ and $g$
leaves invariant each lozenge 
${\cal C}_i$ as well as their corners.
\label{chain}
\end{theorem}


\begin{figure}
\centeredepsfbox{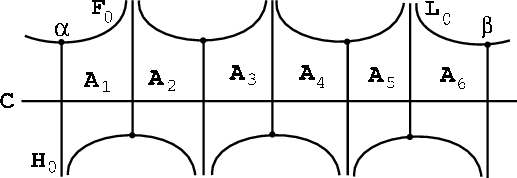}
\caption{
The correct picture between non separated
leaves of $\wls$.}
\label{pict}
\end{figure}

The leaf space of $\wls$ (or $\wlu$) is usually not a Hausdorff
space. Two points of this space are non separated if they
do not have disjoint neighborhoods in the respective leaf
space.
The main result concerning non Hausdorff behavior in the leaf spaces
of $\wls, \wlu$ is the following:

\begin{theorem}{(\cite{Fe5}, theorem 4.9)}{}
Let $\Phi$ be a pseudo-Anosov flow in $M^3$. 
Suppose that $F \not = L$
are not separated in the leaf space of $\wls$.
Then $F$ and $L$  are periodic.
Let $F_0, L_0$ be the line leaves of $F, L$ which
are not separated from each other.
Let $V_0$ be the sector of $F$ bounded by
$F_0$ and containing $L$.
Let $\alpha$ be the periodic orbit in $F_0$ and
$H_0$ be the component of $(\wu(\alpha) - \alpha)$ 
contained in $V_0$.
Let $g$ be a non trivial covering translation
with $g(F_0) = F_0$, $g(H_0) = H_0$ and $g$ leaves
invariant the components of $(F_0 - \alpha)$.
Then $g(L_0) = L_0$.
This produces closed orbits of $\Phi$ which are freely
homotopic in $M$.
Theorem \ref{chain} then implies that $F_0$ and $L_0$ are connected by
a finite chain of lozenges 
$\{ A_i \}, 1 \leq i \leq n$,
consecutive lozenges are adjacent.
They all intersect a common stable leaf $C$.
There is an even number of lozenges 
in the chain, see
fig. \ref{pict}.
In addition 
let ${\cal B}_{F,L}$ be the set of leaves of $\wls$ non separated
from $F$ and $L$.
Put an order in ${\cal B}_{F,L}$ as follows:
The set of orbits of $C$ contained
in the union of the lozenges and their sides is an interval.
Put an order in this interval.
If $R_1, R_2 \in {\cal B}_{F,L}$ let $\alpha_1, \alpha_2$
be the respective periodic orbits in $R_1, R_2$. Then
$\wu(\alpha_i) \cap C \not = \emptyset$ and let 
$a_i = \wu(\alpha_i) \cap C$.
We define $R_1 < R_2$ in ${\cal B}_{F,L}$ 
if $a_1$ precedes
$a_2$ in the order of the set of orbits of 
$C$.
Then ${\cal B}_{F,L}$
is either order isomorphic to $\{ 1, ..., n \}$ for some
$n \in {\bf N}$; or ${\cal B}_{F,L}$ is order
isomorphic to the integers ${\bf Z}$.
In addition if there are $Z, S \in \wls$ so that
${\cal B}_{Z, S}$ is infinite, then there is 
an incompressible torus in $M$ transverse to 
$\Phi$. In particular $M$ cannot be atoroidal.
Also if there are $F, L$ as above, then there are
closed orbits $\alpha, \beta$ of $\Phi$ which
are freely homotopic to the inverse of each other.
Finally up to covering translations,
there are only finitely many non Hausdorff
points in the leaf space of $\wls$.
\label{theb}
\end{theorem}

%

Notice that ${\cal B}_{F,L}$ is a discrete set in this order.
For detailed explanations and proofs, see
\cite{Fe4,Fe5}. 

\vskip .1in
\noindent
{\bf {Scalloped regions}} 

\begin{figure}
\centeredepsfbox{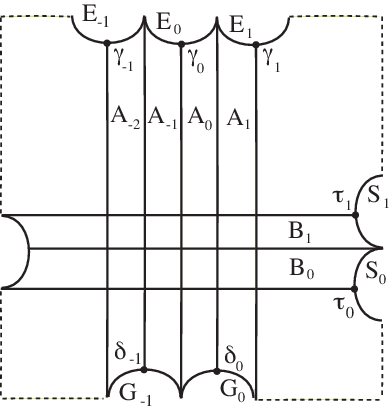}
\caption{A scalloped region ${\cal S}$. 
The collections $\{ E_i \}_{i \in {\bf Z}}, \
\{ G_i \}_{i \in {\bf Z}}$ of stable leaves
are part of the boundary of ${\cal S}$.
In addition $\{ S_i \}_{i \in {\bf Z}}$ are unstable leaves in the boundary 
of ${\cal S}$.
For better viewing we indent a few of the non separated leaves
in (say) $\{ E_i \}_{i \in {\bf Z}}$ into the
square.
Similarly for $\{ G_i \}, \{ S_i \}$. }
\label{scal}
\end{figure}

Suppose that ${\cal E} = \{ E_i \ | \ i \in {\bf Z} \}$ is a
bi-infinite collection of leaves of $\wls$ or $\wlu$ all of
which are non separated from each other and ordered as in
theorem \ref{theb}. 
There is an associated  structure in $\mi$ or $\oo$, which is called
a scalloped region, which we now describe.
Let $\{ A_i \ | \ i \in {\bf Z} \}$
be the bi-infinite collection of lozenges associated to ${\cal E}$ $-$
consecutive $A_i$'s are adjacent.
For simplicity assume that ${\cal E}$ is a collection of 
stable leaves, so that every $A_i$ intersects a fixed stable leaf $\zeta$.
The $A_i$ are chosen so that each $E_i$ has a half leaf in the boundary
of $A_{2i}$ and another half leaf in the boundary of $A_{2i-1}$. 
Each leaf $E_i$ contains a periodic orbit $\gamma_i$.
Let $W_i$
be the half leaf of 
$\wu(\gamma_i)$ which is in the boundary of both $A_{2i}$ and $A_{2i-1}$.
In addition since $A_{2i}$ and $A_{2i+1}$ are also adjacent there is
a stable leaf $G_i$ which has half leaves in the closure of each
of $A_{2i}$ and $A_{2i+1}$. 
Hence $\{ G_i \ | \ i \in {\bf Z} \}$ is another collection of leaves
of $\wls$ non separated from each other. Each $G_i$ contains a periodic
orbit $\delta_i$ and $\wu(\delta_i)$ has a half leaf $Y_i$ which is
in the closure of both $A_{2i}$ and $A_{2i+1}$.
The {\em {scalloped region}} associated to ${\cal E}$ is

$$ {\cal S} \ \ = \ \ \bigcup_{i \in {\bf Z}} \ (A_i \ \cup 
\ W_i \ \cup \ Y_i),$$

\noindent
see fig. \ref{scal}.

Scalloped regions were introduced for Anosov flows in
section 5, 
theorem 5.2 of \cite{Fe3}, but the same analysis works for pseudo-Anosov flows,
mainly because there can be no singularities in the lozenges 
\cite{Fe4}.
It is proved in \cite{Fe3} that such a scalloped region ${\cal S}$ (where the 
$E_i$ are stable leaves) is also the
union of another  bi-infinite collection of lozenges $\{ B_i \ | \ i \in {\bf Z} \}$
and stable half leaves in the boundary of pairs of consecutive lozenges.
All of the lozenges $B_i$ intersect a fixed {\underline {unstable}} leaf.
Therefore the foliations $\wls, \wlu$ restricted to ${\cal S}$
form a product structure in ${\cal S}$,
they both have leaf space which is homeomorphic to  $\rrrr$.
In this way the boundary $\partial {\cal S}$ also has two bi-infinite collections
of leaves of $\oou$. In each collection all leaves are non separated
from each other. Let $\{ S_j \}_{j \in {\bf Z}}$ be the collection
which is in the limit of the sequence $\wu(\gamma_i)$ (or equivalently
$\wu(\delta_i)$) when $i$ converges to plus infinity. The other bi-infinite
collection of unstable leaves is obtained as the limit of $(\wu(\gamma_i))$
as $i$ converges to minus infinity.
We may choose the indexing of the $\{ S_j \}$ so that $S_j$ has one
half leaf in the closure of $B_{2j}$ and another in the closure of $B_{2j-1}$.
Let $\tau_j$ be the periodic orbit in $S_j$. We may also choose the
indexing so that $(\ws(\tau_j))$ converges to the collection
$\{ E_i \}_{i \in {\bf Z}}$
when $i \rightarrow \infty$ and $(\ws(\tau_j))$ converges 
to $\{ G_i \}_{i \in {\bf Z}}$ when $i \rightarrow -\infty$.
We also call a scalloped region 
the projection of ${\cal S}$ to the orbit space $\oo$.

Here is an actual model for a scalloped region in $\oo$.
Let $I, J$ be two properly embedded, order preserving images 
of ${\bf Z}$ into $(-1,1)$ which are intercalated, for example
$J = \{ \pm (1 - \frac{1}{2n}) \ | \ n \geq 1 \}$ and 
$I = \{ \pm (1 - \frac{1}{2n-1}) \ | \ n \geq 1 \}$.
The closure of a scalloped region is a proper embedding of the set

$$V \ =  ([-1,1] \times [-1,1]) \ - \ 
\Big( (J \times \{ 1 \}) \cup (\{ 1 \} \times J)
\cup (I \times \{ -1 \}) \cup (\{ -1 \} \times I)
\cup ( \{ -1, 1 \}  \times \{ -1, 1 \}) \Big)$$

\noindent
into $\oo$ satisfying the following
conditions:
The horizontal and vertical foliations of $\rrrr^2$ restricted to $V$
are mapped to the stable and unstable foliations in $\overline {\cal S}$.
The interior of $V$ maps to the scalloped region.
It is crucial that $I, J$ do not intersect. For example the stable leaf
$(-1/2, 1/2) \times \{ 1 \}$ is one of the $E_i$, we may assume that it is $E_0$.
Then $(0,1)$ is the periodic orbit $\gamma_0$ and $\{ 0 \} \times (-1,1)$ is
the half leaf of $\wu(\gamma_0)$ which is in the boundary of the
lozenges \ $A_{-1} = (-1/2,0) \times (-1,1)$ \ and \ 
$A_0 = (0,1/2) \times (-1,1)$. \ It is crucial
in this particular example that $(0,-1)$ is {\underline {not}} in $V$.
We may assume that $S_0 = \{ 1 \} \times (-1/2, 1/2)$.

In fig. \ref{scal} we indent the region along the boundary stable and
unstable leaves to highlight that they form collections of non separated leaves.

\begin{theorem}{(\cite{Fe5}, theorem 4.10)}{}
Let $\Phi$ be a pseudo-Anosov flow. Suppose that there is
a stable or unstable product region. Then $\Phi$ is 
topologically conjugate to a suspension Anosov flow.
In particular $\Phi$ is nonsingular.
\label{prod}
\end{theorem}


\section{Ideal boundaries of pseudo-Anosov flows}

Let $\Phi$ be a pseudo-Anosov flow in $M$. The orbit
space $\oo$ of $\wwp$ (the lifted flow to $\mi$) 
is homeomorphic to $\rrrr^2$ \cite{Fe-Mo}. 
In this section we construct a 
natural compactification of $\oo$ with an ideal circle
$\partial \oo$ called the ideal boundary of the pseudo-Anosov
flow. We put a topology in $\cd = \oo \cup \partial \oo$ making it
homeomorphic to a closed disk.
The induced action of $\pi_1(M)$ on $\oo$ 
extends to an action on $\oo \cup \partial \oo$.
This works for any pseudo-Anosov flow in a $3$-manifold
$-$ no metric, or topological  assumptions (such as atoroidal) on $M$
or on the flow $\Phi$.
In addition there are no assumptions about perfect
fits for $\Phi$ or concerning topological conjugacy to suspension
Anosov flows.

One key aspect here is that we want to use only the foliations
$\oos, \oou$ to define $\partial \oo$ and its topology.

Before formally defining ideal
points of $\oo$ we analyse some examples.
Given $g$ in $\pi_1(M)$ it acts on $\mi$ and sends flow
lines of $\wwp$ to flow lines and hence acts on $\oo$.
This action preserves the foliations $\wls, \wlu, \oos, \oou$.
Recall that a
$2$-dimensional foliation $\fol$ in a $3$-manifold $N$ is
called $\rrrr$-covered if the leaf space of $\fn$ is
homeomorphic to the real line \cite{Fe1}.
An Anosov flow is $\rrrr$-covered if $\ls$ (or equivalently 
$\lu$ \cite{Ba1}) is $\rrrr$-covered.

\vskip .1in
\noindent
1) Ideal boundary for $\rrrr$-covered Anosov flows. The product
case.

\begin{figure}
\centeredepsfbox{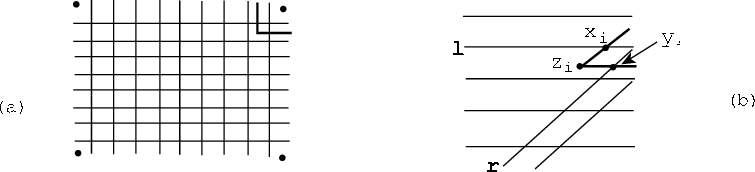}
\caption{
Ideal points for product $\rrrr$-covered
Anosov flow, the dots represent the 4 special points,
b. The picture in skewed case.}
\label{rcov}
\end{figure}

A product Anosov flow is an Anosov flow for which both 
$\ls, \lu$ are $\rrrr$-covered and in addition
every leaf of $\oos$
intersects every leaf of $\oou$ and vice versa \cite{Fe1,Ba1}.
Barbot proved that this implies that
$\Phi$ is topologically conjugate to
a suspension \cite{Ba1}.
Every ray in $\oos$ or $\oou$ generates a point
of $\partial \oo$ and they are all distinct. Furthermore
there are 4 additional ideal points corresponding to escaping
quadrants in $\oo$, see fig. \ref{rcov}, a.
The quadrants are bounded by a ray in $\oou$ and a ray
in $\oos$ which intersect only in their common starting
point (or finite endpoints).
In this case it is straightforward to put a topology in
$\cd = \oo \cup \partial \oo$ so that 
it is a closed disk and covering transformations
act on the extended object. 
If $\Lambda^s, \Lambda^u$ are both transversely orientable,
then
any covering translation $g$ fixes
the 4 distinguished points. It is associated to a periodic
orbit if and only if it fixes 4 additional ideal points:
if $x$ in $\oo$ satisfies $g(x) = x$, then $g$ fixes
the ``ideal points" of rays of $\oos(x), \oou(x)$.
When $\Lambda^s, \Lambda^u$ are not transversely orientable,
there are other restricted possibilities.

We want to define a topology in $\cd$ using only
the structure
of $\oos, \oou$ in $\oo$.
A distinguished ideal point $p$ has a neighborhood
basis determined by (say nested) pairs of rays in $\oos, \oou$ 
intersecting at their common finite endpoint and so that the
corresponding quadrants ``shrink" to $p$.
For an ordinary ideal point $p$, say a stable ideal point
of a ray in $\oos(x)$, 
we can use shrinking strips: the strips are bounded by 2 rays
in $\oos$ and a segment in $\oou$ connecting the endpoints of
the rays. The unstable segment intersects the original stable
ray of $\oos(x)$ 
and the intersections escape in that ray and also shrink
in the transversal direction.
Already in this case this leads to an important concept:

\begin{define}{(polygonal path)}{}
A polygonal path in $\oo$ is a properly embedded, bi-infinite path
$\zeta$
in $\oo$ satisfying:  either $\zeta$ is a leaf of $\oos$ or $\oou$ or
$\zeta$ is the union of a finite collection
$l_1, ... l_n$ 
of segments and  rays in leaves of $\oos$ or $\oou$
so that
$l_1$ and $l_n$ are rays in $\oos$ or $\oou$ and  the other $l_i$
are finite segments. 
We require that 
$l_i$ intersects $l_j$ if and only if $|i-j| \leq 1$.
In addition the $l_i$ are alternatively in $\oos$ and $\oou$.
The number $n$ is the length of the polygonal path.
The points $l_i \cap l_{i+1}$ are the vertices of the path.
The edges of $\zeta$ are the $\{ l_i \}$.
\end{define}

In the product $\rrrr$-covered case, the exceptional ideal
points need neighborhoods basis formed by polygonal paths of length 2 and
all the others need polygonal paths of length 3.

\vskip .1in
\noindent
2) $\rrrr$-covered Anosov flows $-$ skewed case.

This is an Anosov flow so that $\ls, \lu$ are $\rrrr$-covered 
and the following is satisfied:
Topologically the orbit space $\oo$ is homeomorphic to $(0,1) \times \rrrr$,
a subset of the plane, so that stable leaves are 
horizontal segments and unstable leaves are segments making
a constant angle $\not = \pi/2$ with the horizontal,
see fig. \ref{rcov}, b.
A leaf of $\oos$ does not intersect every leaf of $\oou$ and
vice versa \cite{Fe1,Ba2}.
Here again each ray of $\oos$ or $\oou$ defines an ideal point
of $\oo$. However as is intuitive from the picture,
rays of 
$\oos, \oou$ 
which form a perfect fit in $\oo$
should define the same ideal point of $\oo$.
In addition to these ideal points of rays of leaves in $\oos$ or $\oou$,
there should be 2 distinguished ideal points $-$ one
from the ``positive" direction of $\rrrr$ and one from the 
``negative" direction of $\rrrr$. 
Hence $\cd$ is equal to $[0,1] \times \rrrr$ union two points: 
one for the positive end of $\rrrr$ and one for the negative end.
Put a topology in $\cd$ so that $[0,1] \times \rrrr$ is homeomorphic
to a disk minus two boundary points. 
Covering translations act
as homeomorphisms of this disk. A transformation without
fixed points in $\oo$ fixes only the 2 distinguished ideal
points in $\partial \oo$, one attracting and another repelling.
If a transformation $g$ has a fixed point $p$ in $\oo$, then it leaves
invariant the leaf $l = \oos(p)$ of $\oos$.
If $g$ switches the components of $l - \{ p \}$, then $g$ does not
fix any point in $\partial \oo$.
Otherwise there
are infinitely many fixed points, see \cite{Fe1,Ba2}.


A neighborhood basis of the distinguished ideal points can be
obtained from leaves of $\oos$ or $\oou$ which escape in
that direction (positive or negative). 
For non distinguished ideal points, we get sequences of polygonal paths
of length 2 escaping every compact set and ``converging" to this
ideal point, see fig. \ref{rcov}, b.
More precisely if rays $l, r$ of $\oos, \oou$ respectively form
a perfect fit defining the ideal point $p$, then choose 
$x_i$ in $l$ and escaping in the direction of
the perfect fit and similarly chose  $y_i$ in $r$.
Consider the 
polygonal path of length two containing rays in the stable
leaf through $y_i$ and the unstable leaf through $x_i$
(intersecting in $z_i$, see fig. \ref{rcov}, b).

\vskip .1in
\noindent
3) Suspension pseudo-Anosov flows $-$ singular case.

The fiber is a hyperbolic surface.
The orbit space
$\oo$ is identified with the universal cover of the
fiber which is metrically 
the hyperbolic plane $\hh$. 
There is a natural ideal
boundary $\su$ $-$ the circle at infinity of $\hh$.
One expects that $\partial \oo$ and $\su$ should be the equivalent.
But the construction of $\su$
uses the {\underline {metric}} structure on the
surface $-$ in general there is no metric structure
in $\oo$, so again we want to define $\partial \oo$ using only
the structure of $\oos, \oou$. 
From a geometric point of view,
there are some points of $\su$ which are ideal points
of rays of leaves of $\oos$ or $\oou$.
But there are many other points in $\su$. The foliations $\oos, \oou$
can be split into geodesic laminations (of $\hh$) which 
have only complementary regions which are finite sided
ideal polygons.
This implies that given $p$ in $\su$ 
there is always a sequence of leaves
$l_i$ (in $\oos$ or $\oou$) which is nested,
escapes to infinity and ``shrinks" to the ideal point
$p$. In this way one can characterize all
points of $\su$ using only the
foliations $\oos, \oou$ and hence $\partial \oo = \su$ in this case.
Also $\oos, \oou$ define a topology in $\oo \cup \partial \oo$
compatible with the metric topology.

\vskip .1in
Now we analyse a potential difficulty.
Let $l$ be a nonsingular ray (say) in $\oos$
and let $x_i$ in $l$, forming a nested sequence of points in $l$,
escaping compact sets in $l$. 
For simplicity assume that the 
leaves $g_i$ of $\oou$ through $x_i$ are nonsingular.
We would like to say that the sequence $(g_i)$ ``defines" an 
ideal point of $\oo$. 
If the $g_i$ escape compact sets in $\oo$, then this will
be the case. 
However it is not always true that $(g_i)$ escapes in $\oo$.
If they do not escape in $\oo$, then they limit
on a collection of unstable leaves $\{ h_j \ | \
j \in J \}$. 
But there is one of them, call it $h$
which makes a perfect fit with $l$ on that side of $l$.
This non trivial fact is proved in \cite{Fe4}.
The perfect fit $l, h$  is the obstruction to leaves
$g_i$ escaping in $\oo$.

We need a couple of definitions.
A {\em quarter} at $z$ is a component of $\oo - (\oos(z) \cup \oou(z))$.
If $z$ is nonsingular
there are exactly $4$ quarters, if $z$ is a $k$-prong point
there are $2k$ quarters.

\begin{define}{(convex polygonal paths)}{}
A polygonal path $\delta$ in $\oo$ is convex if there is a complementary region $V$
of $\delta$ in $\oo$ so that at any given vertex $z$ of $\delta$ 
the local region of $V$ near $z$ is not a quarter at $z$.
Let $\widetilde \delta  = \oo - (\delta \cup V)$. This region 
$\widetilde \delta$ is the convex
region of $\oo$ associated to the convex polygonal path $\delta$.
\label{convex}
\end{define}

The definition implies that if the region $\widetilde  \delta$ contains
2 endpoints of a segment in a leaf of $\oos$ or $\oou$, then
it contains the entire segment (proved later).
This is why $\delta$ is called convex.
If $\delta$ is a single nonsingular leaf
of $\wls$ or $\wlu$ or if
all the vertices of $\delta$ are singularities, then it is possible
that there are two regions $\widetilde \delta$ which are
convex. In the future the context will make clear which region
we are considering.
If $\delta$ is a polygonal path, $V$ a complementary region and
$p$ a vertex for which $V$ is a quarter at $p$, then $p$ is 
called a non
convex vertex of $\oo - (\delta \cup V)$.

\begin{define}{(equivalent rays)}{}
Two rays $l, r$ of $\oos, \oou$ are {\em {equivalent}} if there
is a finite collection of distinct rays $l_i, 1 \leq i \leq n$, 
alternatively in $\oos, \oou$ so that $l = l_0, r = l_n$ and
$l_i$ forms a perfect fit with $l_{i+1}$ for \ $1 \leq i < n$.
\end{define}

It is important to notice that this is strictly about rays in
$\oos, \oou$ and not leaves of $\oos, \oou$. More specifically
we want consecutive perfect fits to be in the same rays
of the adjoining leaf. This implies for instance that 
if $n \geq 3$ then for all $1 \leq i \leq n - 2$ the
leaves $l_i$ and $l_{i+2}$ are non separated from each other in the
respective leaf space.

\begin{define}{(admissible sequences of paths)}{}
An admissible sequence of polygonal paths in $\oo$ is a sequence of convex 
polygonal paths 
$(v_i)_{i \in {\bf N}}$ so that the associated convex regions 
$\widetilde v_i$ form a nested 
sequence of subsets of $\oo$, which escapes compact sets in $\oo$ and
for any $i$,
the two rays at the
ends of $v_i$ are not equivalent.
\end{define}

The fact that the $\widetilde v_i$ are nested and escape compact sets
in $\oo$ implies that the $\widetilde v_i$ are uniquely defined given
the $v_i$.

\vskip .1in
\noindent
{\bf {Structure of this section}} $-$ 
The construction of the ideal compactification of $\oo$ and the analysis
of its properties is very involved and complex. This will take
all of this very long section, so here is an outline of the section:
Ideal points of $\oo$ will be defined by admissible sequences
of polygonal paths, definition \ref{idpo}.
But many admissible sequences generate the same ideal point,
so we first define a relation in the set of admissible
sequences, definition \ref{rela}. We establish a technical result
called the fundamental lemma (lemma \ref{jumps}) which implies
that the relation above is an equivalence relation, lemma
\ref{equiv}. 
In definition \ref{idpo} we define ideal points of 
$\oo$ producing $\partial \oo$ and with union $\cd = \oo \cup \partial \oo$.
Some special ideal points are defined in
definition \ref{stan} associated to ideal points of rays
of $\oos$ or $\oou$ and in lemma \ref{infin} we deal with 
infinitely many leaves of $\oos$ or $\oou$ all non separated
from each other.
Not every admissible sequence is efficient to study ideal points
of $\oo$ and we define master sequences in definition \ref{master}:
roughly the rays in the polygonal paths of these sequences
approach the ideal point of $\oo$ from 
``both" sides. In lemma \ref{mast} we prove that any ideal point
admits a master sequence and they are used to distinguish points
of $\partial \oo$.
In definition \ref{topology} we define a topology for $\cd = \oo 
\cup \partial \oo$ and in lemma \ref{topol} we prove that this is 
indeed a topology in $\cd$.
We then progressively prove stronger properties of $\cd$:
Lemma \ref{haus} shows that $\cd$ is Hausdorff, 
lemma \ref{fc} shows that $\cd$ is first countable
and lemma \ref{sec} shows that $\cd$ is second countable - this last
one is a bit more complicated than the other ones.
These and the structure of $\cd$ quickly imply that
$\cd$ is regular (lemma \ref{regul}) and hence metrizable.
Then we study compactness properties: first we prove a technical
and very tricky lemma about a special case (lemma \ref{forcing}).
This lemma considerably simplifies the proof of compactness
of $\cd$ (proposition \ref{comp}).
At this point we can quickly prove that the ideal boundary
$\partial \oo$ is homeomorphic to a circle (proposition
\ref{circle}).
We then prove a harder result (theorem \ref{disk}) that
$\cd = \oo \cup \partial \oo$ is homeomorphic to a
closed disk. 
Finally in lemmas \ref{equiva}, \ref{choi}, \ref{infin} and proposition
\ref{options} we prove additional properties of the ideal points
of $\oo$ and which types of admissible sequences are associated
to different types of ideal points.
\vskip .1in

An ideal point of $\oo$ will be determined by an 
admissible sequence of paths.
Clearly this does not work for suspension Anosov flows because
a sequence of escaping leaves of $\oos$ approaches infinitely
many different ideal points.
Hence such flows are special and are treated separately.
We abuse notation and say that $(v_i)_{i \in {\bf N}}$ is nested.
For notational simplicity many times we denote such a sequence
by $(v_i)$.

Two different admissible sequences may define the same ideal
point and 
we first need to decide when two such sequences
are equivalent.
At first it seems that any 2 sequences associated to 
the same ideal point of $\oo$
would have to be eventually nested with each other. 
However it is easy to see that
such is not the case. For example consider a nested sequence of
rays of a fixed leaf $l$.
We will later see how to
extend each ray on one side of $l$ to form an admissible sequence.
Extend them also to the other side to form another admissible
sequence. Intuitively the two sequences should converge to the intrinsic
ideal point of $l$, but clearly they are not eventually nested.

\begin{define}{}{}
Given two admissible sequences of 
chains $C = (c_i)$, \ $D = (d_i)$, we say that 
$C$ is smaller or equal than $D$, denoted by $C \leq D$,
if: for any $i$ there is $k_i > i$ so that
$\widetilde c_{k_i} \subset \widetilde d_i$.
Two admissible sequences of chains $C = (c_i), \ D = (d_i)$
are equivalent and we then write $C \cong D$ if there is
a third admissible sequence $E = (e_i)$ so that
$C \leq E$ and $D \leq E$.
\label{rela}
\end{define}

Ideal points of $\oo$ will be defined as equivalence classes
of admissible sequences of polygonal paths.
Hence we must prove that $\cong$ is
an equivalence class and along the way we derive 
several other properties.
We should stress that the requirement that the chains are
{\underline {convex}} is fundamental for the whole discussion.
It is easy to see in the skewed $\rrrr$-covered Anosov case,
that given any two distinct
ideal points $p, q$ on the ``same side" of the distinguished
ideal points, the following happens:
Let $l, r$ be stable rays defining $p, q$ respectively.
Then there is a sequence of polygonal paths in $\oo$,
that escapes compact sets in $\oo$ and so that each
$\widetilde c_i$ contains
subrays of both $l$ and $r$.
The polygonal paths can be chosen to satisfy
all the properties, except that they are convex. 
On the other hand convexity does imply important properties
as shown in the next lemma.
This key lemma will be used throughout this section.
After this lemma we show that $\cong$ is an equivalence relation.

\vskip .1in
\noindent
{\bf {Singular foliations in surfaces with boundary and index formula}}

Let $\fol$ be a singular foliation
on a compact surface $S$ with boundary, so that interior
singularities are all of $k$-prong type and $k \geq 3$. The foliation
may be tangent to part of the boundary. There is an Euler-Poincare
index formula so that the sum of the indices of the singularities
equals the Euler characteristic of the surface.
An interior singularity with $k$ prongs has index $1 - \frac{k}{2}$.
A boundary singularity has index $\frac{1}{2} - \frac{k}{2} -
\frac{t}{4}$, where $k$ is the number of prongs going into the
surface and $t$ is the number of prongs which are part of the
boundary. 
The possible values of $t$ are $0, 1, 2$.
For example if $k = 0, t = 0$ the singularity is
half of a center, which has index $1/2$.
This will be used for compact subsets of $\oo$, which are 
foliated by $\oos$ or $\oou$.

\begin{lemma}{(fundamental lemma)}{} \
Assume that $\Phi$ is not topologically conjugate to
a suspension Anosov flow.
Let $l, r$ be rays of $\oos$ or $\oou$, which are not equivalent.
Then there is no pair 
of admissible sequence of polygonal paths $E = (e_i)$, \
$F = (f_i)$ so that: \ $\widetilde e_i \cap \widetilde f_i
\not = \emptyset$ (for all $i$) and $\widetilde e_i \cap r \not
= \emptyset$, $\widetilde f_i \cap l \not = \emptyset$,
for all $i$.
\label{jumps}
\end{lemma}

\begin{proof}{}
We assume that both $l$ and $r$ are
rays of $\oos$,
other cases are treated similarly.
By taking subrays if necessary, we may assume that 
$l, r$ are disjoint, have no
singularities  and miss a compact set containing the
base point in $\oo$. Join the initial points of $l, r$ by
an arc $\alpha'$ missing this  big compact set to 
produce a properly embedded bi-infinite
curve  $\alpha = l \cup \alpha' \cup r$, see fig. 
\ref{refu}, b.  
Let $V$ be the component of $\oo - \alpha$ which misses the basepoint.

\begin{figure}
\centeredepsfbox{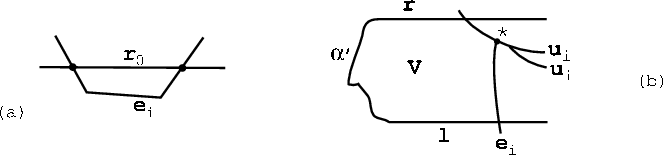}
\caption{
a. Convexity implies connected
intersection of $r$ and $B_i$. 
b. All rays of $u_i$ stay in $V$ forever.
There is a non convex
vertex at $*$.}
\label{refu}
\end{figure}

\vskip .1in
\noindent
{\underline {Case 1}} $-$ $E = F$.

Here we have to show that there is 
no admissible sequence of polygonal paths $E = (e_i)$ 
such that $\widetilde e_i$
always intersects $l$ and $r$. This implies that
the phenomenon described
above (in the skewed Anosov flow case)
for non convex polygonal paths cannot happen for convex polygonal paths.
Suppose this is not true and let $E = (e_i)$ be one such sequence.  
Let $B_i = \widetilde e_i \cup e_i$.

\vskip .1in
\noindent
{\underline {Claim 1}} $-$ If $e_i$ is a convex polygonal path
with region $\widetilde e_i$ and $r$ is a leaf of $\oos$ (or 
of $\oou$),
then $(\widetilde e_i \cup e_i) \cap r$ is connected.

Otherwise there is a compact subarc $r_0$ of
$r$ with $\partial r_0$ in $e_i$ and the rest of
$e_i$ contained in $\oo - B_i$, see fig. \ref{refu}, a.
There is a compact arc $\tau$ in $e_i$ joining the endpoints $x,y$
of $r_0$. Let $D$ be the disc in $\oo$ bounded $r_0 \cup \tau$
and consider the foliation $\oos$ induced in $D$. The singularities
in the interior are $k$ prong type all with negative index.
At $x$ there is a boundary prong of $\oos$ (since $r_0$ is in the
boundary of $D$) so the index is $\leq 1/4$ and similarly for
$y$. 
If there are singularities in the interior of $r_0$ then
they have negative index as $r_0$ is contained in a lef of $\oos$.
Since the Euler characteristic of the disc is $1$ and there
are no half centers in $\tau$, \
  all singularities in $\tau$ have index $\leq 1/4$.
It follows that there must be 
at least {\underline {two}} boundary singularities in $\tau - \{ x, y \}$ with index $1/4$.
Each one of these has to be a point
$z$ so that there is a prong of $\oos$ and
a prong of $\oou$ locally contained in $\tau \subset \partial D$ and no
other prongs of $\oos \cup \oou$ entering $D$. The unstable
prong is transverse to $\oos$.
This shows there is a quarter of $D$ at $z$.
But since $r_0 - \{ x, y \} \subset (\oo - B_i)$ this means
that $\widetilde e_i$ has a non convex vertex at $z$, contradiction
to $e_i$ being convex.
Therefore $B_i \cap r$ is connected
and this proves claim 1.
This is the convexity property of 
$\widetilde e_i$ mentioned after 
definition \ref{convex}.

\vskip .1in
We continue the analysis of case 1.
Notice that $(B_j \cap r)_{j \in {\bf N}}$ is a nested family of 
non empty sets
in $r$. Since $B_j$ escapes compact sets as $j \rightarrow 
\infty$ and $B_j \cap r$ is connected, it follows that
$B_j \cap r$ is a subray of $r$ for any $j$.
If $e_j \cap r$ contains a non trivial  segment, 
then again by convexity and Euler characteristic
it follows that $\widetilde e_j \cap r = 
\emptyset$ contradiction.
Hence $e_i$ intersects $r$ in a single point. 
Let $u'_i = \oou(e_i \cap r)$ be the unstable leaf through
the intersection.
Up to subsequence, we may assume no two $u'_i$ are the same.

Since $r$ has no singularities there are two components of
$u'_i - (u'_i \cap r)$. There is only one of them
denoted by 
$u_i$ which locally enters $V$ at the intersection, see fig. \ref{refu}, b.
There are two subcases:

\vskip .1in
\noindent
{\underline {Case 1.a}} $-$
Some ray
of $u_i$ stays in $V$ for all time.

Let this ray be $s$. 
Then $s$ is properly embedded in $V$ and together with a subray of
$r$ it bounds a subregion $W$ of $V$.
It follows that by taking a bigger $i$ if necessary we may assume that all
rays of $u_i$ stay in $V$ forever $-$ because they 
are in the region $W$ above.
Take the ray $s$ of $u_i$ starting at $u'_i \cap r$ and fartherst from $r$ or
equivalently closest to $l$. Even though $r, l$ are rays and do not separate
$\oo$, 
this makes sense because $V$ is an open disc with boundary $\alpha$
and $l, r$ are disjoint subrays of $\alpha$.
All rays of $u_i$ start in $r$ and the collection of rays of $u_i$ is
(weakly) nested.

In that case, in order for $e_i$ to reach $l$ it leaves $s$ at a point
$*$ where $e_i$ switches to travel along a segment $t$ in $\oos$.
There cannot be any other prong of $\oos(*) \cup \oou(*)$ not in
$\widetilde e_i$: since $s$ is an unstable prong and $t$ is contained in
a stable prong, there would have to be another unstable prong
in $\widetilde e_i$. But this unstable prong is contained in $V$ by construction
and hence not contained in $\widetilde e_i$. Hence this shows that $*$ is a non
convex vertex in $e_i$, see fig. \ref{refu}, b. This is a contradiction to
$e_i$ convex.

\vskip .1in
\noindent
{\underline {Case 1.b}} $-$
For any $i$, all rays of $u_i$ exit $V$.

We first want to show that the sequence $u_i$ does not escape compact
sets in $\oo$. Then we show that a leaf $u$ in the limit of $(u_i)$
has a ray which makes a perfect fit with $r$ and we restart the proof
with $u, l$ in place of the rays $r, l$.

Suppose first that all $u_i$ intersect $l$.
In that case let $z_i$ be the part of $u_i$ between
$l$ and $r$. If the $z_i$ escapes compact sets in $\oo$,
then the region between $l$ and $r$ is an unstable product region
as in Definition
\ref{defsta}. Theorem \ref{prod} then implies that $\Phi$
is topologically conjugate to a suspension 
Anosov flow. This is disallowed by hypothesis
(in fact
the lemma fails for product $\rrrr$-covered Anosov flows).
Hence the $u_i$ does not escape compact sets in $\oo$.
The other option is that the $u_i$ does not intersect
$l$ $-$ hence they intersect $\alpha'$.
Since $\alpha'$ is compact, then in all cases $u_i$ does not
escape compact sets in $\oo$.

The intersection of $\overline u_i$ with $r$ escapes
in $r$, and $(u_i)$ is a  nested collection (as subsets of $V$), so 
$u_i$ converges to a  collection
of (line) leaves of $\oou$. Let $u$ be one of the limit leaves.
Consider the set $B$ of unstable leaves non separated from $u$
and which are either contained in $V$ or intersect $\alpha$.
By theorem \ref{theb}
there is an order in the set $B$ and there are only finitely many unstable leaves
between any given $u$ 
and $r$, so we may assume that $u$ is the leaf in $B$ which 
is the closest one to $r$ in terms of this order.

\vskip .1in
\noindent
{\underline {Claim 2}} $-$ $u$ makes a perfect fit with $r$.

Suppose that $u$ does not make a perfect fit with $r$.
We will produce a product region.
Let $z$ a point in $u$. The stable leaf through $z$ intersects
$u_i$ for a fixed $i$ big. For any other $w$ in $u$ 
then $\oos(w)$ intersects $u_j$ for some $j > i$. We say that
$w$ is closer to $r$ than $z$ if the intersections
$\oos(z) \cap u_j, \oos(w) \cap u_j, \overline u_j \cap r$ are
linearly ordered in $u'_j$.
Hence $\oos(w)$ also intersects the fixed $u_i$.
It follows that 
as $w$ escapes in $u$ in the
direction of $r$, the $\oos(w)$ converge to a stable leaf $r'$ which makes a perfect
fit with $u$. Hence $r, r'$ are distinct.
The region between
$r, r'$ is a product region because
all the $u_j$ ($j \geq i$)  intersect $r, r'$ and there are no limit leaves
of the $(u_j)$ between $r, r'$.
As seen above,
this would imply $\Phi$ is topologically conjugate
to a suspension Anosov flow, contradiction.
This proves claim 2.

\vskip .1in
The rest of case 1 concerns only flows with perfect fits.

We now show that $u$ is not contained in $V$. If $u$ is contained in $V$,
there are two cases: i) $u \subset \widetilde e_i$ for all $i$ $-$ but
this contradicts that $\widetilde e_i$ escapes compact sets of $\oo$;
\ ii) there is $i$ with $u$ not contained in $\widetilde e_i$. But then
$e_i$ has to cross $u$, and since $u$ is contained in $V$, then $e_i$ has
to cross $u$ again in order to intersect $l$.
This produces two intersections of $e_i$ with $u$,
which is disallowed by claim 1.

It follows that 
 there is a
ray of $u$ exiting $V$.
We now restart the argument with $u, l$ instead of $r, l$.
The same arguments as above 
produce a line leaf $v_1$ of $\oos$ making a perfect fit with
$u$ and $v_1$ exiting $V$. In addition $v_1$ is non separated from $r$ 
in the leaf space of $\oos$ $-$ because of the perfect fits
$r \rightarrow u \rightarrow v_1$.
Now iterate to obtain
$v_2, v_3 ...$. This is a nested collection and the sequence $v_j$
cannot accumulate anywhere in $\oo$, since $v_k, v_{k+2}$ are
non separated from each other in the corresponding leaf space.
In addition no two consective unstable leaves in the sequence can
intersect $l$ as they are non separated from each other. It follows that
none of them intersect $l$ and so they all intersect $\alpha'$, which
is compact. This contradicts the fact that they escape in $\oo$.
This proves that 
no escaping sequence of
convex polygonal paths can always intersect both $l$ and $r$.
This finishes the analysis of Case 1.

\vskip .15in
\noindent 
{\underline {Case 2}} $-$ $E \not = F$.

Let $r, l$ as in the statement of the lemma 
and suppose that $E = (e_i), F = (f_i)$ are admissible sequences
with $\widetilde e_i \cap \widetilde f_i \not = \emptyset$,
$r \cap \widetilde e_i \not = \emptyset, \ l \cap \widetilde f_i
\not = \emptyset$, for all $i$. 
As before consider the region $V$ bounded by
$l, r$ and an arc $\alpha'$ connecting them.
By case 1,  $\widetilde e_i$ eventually stops intersecting $l$. Discarding
the initial terms we can assume that $\widetilde e_i \cap l = \emptyset$ and
$\widetilde f_i \cap r = \emptyset$ for all $i$.

\begin{figure}
\centeredepsfbox{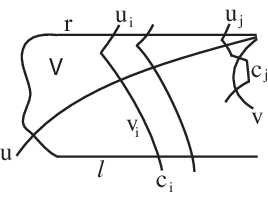}
\caption{
Two polygonal chains and perfect fits.}
\label{relu}
\end{figure}

We construct a polygonal path $c_i$ as follows: first consider the part of
$e_i$ outside of 
$V$. Then add the edges (or pieces of edges) 
of $e_i$ until it first meets $f_i$, then switch to $f_i$ and follow along 
the rest of $f_i$ 
in the direction that intersects $l$. There is only one such direction
as $f_i$ intersects $l$ in a single point and notice that $e_i$ does
not intersect $l$.
This path $c_i$ separates $\oo$ and has a complementary component 
$\widetilde c_i$ which contains subrays of $l, r$. This component contains all of
$V$ except for a subset contained in a compact set of $\oo$.

The vertices of $c_i$ are all convex for $\widetilde c_i$,
except perhaps for the single vertex $p_i$ where $c_i$ changes from $e_i$ to $f_i$.
Once the non convex vertex appears, all
subsequent vertices have to be convex.

As before consider the unstable leaf $u_i$ through
$e_i \cap r$.
If some $u_i$ has a ray which is entirely in $V$, then
as seen in case 1,
for $j > i$ all rays of $u_j$ which enter $V$ must
be entirely in $V$.
This implies that the change from $e_i$ to $f_i$ has to be in $u_i$.
Here is why:
otherwise  the next edge in $c_i$ is $w_i$ an edge still in $e_i$. But since 
$c_i$ eventually has to cross $l$, and $u_i$ is entirely contained in $V$, it follows
that $c_i$ has to intersect $u_i$ at twice. As seen in 
the proof of claim 1,
this implies the
existence of {\underline {two}} non convex vertices in $c_i$. But $c_i$ has only
one non convex vertex, contradiction.


We conclude that all rays of $u_i$ which enter $V$ have to
exit $V$. As seen in case 1 they cannot escape compact sets in $\oo$.
They converge to a collection
of (line) leaves in $\oou$. As in case 1, one of them, call it $u$ makes
a perfect fit with $r$. 
Since $u, r$ make a perfect fit and $\widetilde e_i$ escapes compact
sets, it follows that for $i$ big $e_i$ intersects $u$ and the
second edge of $e_i$ is in leaves $v_i$ of $\oos$ and $v_i$ intersects $u$.

The first possibility here is that $u$ contained in $V$.
Let $W$ be the component of $\oo - u$ contained in $V$.
Since $r, u$ make a perfect fit and $\widetilde f_i \cap \widetilde e_i 
\not = \emptyset$ it follows that $\widetilde f_i$ has to intersect
$W$. Since $u \subset V$, then $c_i$ has to intersect $u$ twice $-$
this is a contradiction as seen before.
The second possibility is that $u$ is not contained in $V$ and intersects
$\alpha$.
Notice that $u$ is a ray equivalent to $r$. We can now restart the
proof of case 2 with $u, l$ instead of $r, l$. The arguments above will produce
a leaf $v$ of $\oos$ making a perfect fit with $u$. Figure \ref{relu}
illustrates the impossible situation that $v \subset V$. In that case some
$c_j$ is forced to have 2 non convex vertices. Hence 
 $v$ intersects
$\alpha$. As in case 1, one can iterate this argument to arrive at
a contradiction.

This finishes the proof of lemma \ref{jumps}.
\end{proof}

\noindent
{\bf {Remarks}} $-$ If $A = (a_i)$ is an admissible sequence and
$B = (a_{i_j})$ is a subsequence, then clearly $B$ is also an
admissible sequence and furthermore $A \leq B$ and $B \leq A$.
It is also immediate from the nesting property that if $A = (a_i)$,
$C = (c_i)$ are admissible sequences, then the condition that
$\widetilde a_i \cap \widetilde c_j \not = \emptyset$ for all 
$i, j$ is equivalent to
$\widetilde a_i \cap \widetilde c_i \not = \emptyset$ for all  $i$.

\begin{lemma}{}{}
Suppose that $\Phi$ is not topologically conjugate to a suspension
Anosov flow.
Then the relation
 $\cong$ is an equivalence relation for admissible sequences of 
polygonal paths.
\label{equiv}
\end{lemma}

\begin{proof}{}
Clearly $\cong$ is reflexive and symmetric.
Suppose now that $A = (a_i), B = (b_i), C = (c_i)$
are admissible sequences of polygonal paths and $A \cong B, \ B \cong C$.
Then there are $D =  (d_i)$ with $A \leq D, \ B \leq D$ and $E = (e_i)$
with $B \leq E, \ C \leq E$.
If for some $i, j$ the $\widetilde d_i$ and $\widetilde e_j$
do not intersect this contradicts $B \leq D$, $B \leq E$.

\vskip .1in
\noindent
{\underline {Claim}} $-$ Let $j$ be given.  Then either there is 
$i > j$ with $\widetilde a_i \subset \widetilde e_j$ or
there is $i > j$ with $\widetilde c_i \subset \widetilde d_j$.

Along the proof we may replace $j$ by a bigger number $-$ by
the nesting property the result follows for the original $j$.
The proof is by contradiction. 
So assume the claim fails.
For each $i$, then 
$\widetilde a_i \not \subset \widetilde e_j$ and 
$\widetilde c_i \not \subset \widetilde d_j$.
Clearly this implies that none of  
$\widetilde d_j, \widetilde e_j$ is contained in the other.
Define

$$Z'  \ := \ \widetilde e_j \cap \widetilde d_j$$

\noindent
This is an open subset of $\oo$, which is non compact as there is $m \geq j$ with
$\widetilde b_m
\ \subset \ \widetilde e_j \cap \widetilde d_j$.
It is conceivable that 
that even though $\widetilde d_j, \widetilde e_j$
are convex, $Z'$ may not be connected.
In any case let $Z$ be the component of
$Z'$ containing $\widetilde b_m$. Obviously $Z$ is non compact. 
Notice that $\partial Z$ is made up of segments or rays in $e_j$ or $d_j$.
In addition $\partial Z$ has at least two infinite rays because $Z$ is
non compact.
It is easy to prove that $\partial Z$ is convex for $Z$ because
of this property for $d_j, e_j$.

We first deal with the following situation. Suppose that 
$\partial Z$ has two bi-infinite components. Then $e_j, d_j$
do not intersect and the
region between $d_j$ and $e_j$ is equal to $Z$.
Let $\alpha$ be an arc intersecting $e_j, d_j$ only in its boundary.
We can assume that $\alpha$ does not intersect $\widetilde b_m$.
Since $(\widetilde d_k)$ escapes compact sets in $\oo$, then
it eventually stops intersecting $\alpha$, so choose $k > j$
with $\widetilde d_k \cup d_k$ not intersecting $\alpha$.
If $e_k$ does not
intersect $d_k$, 
then either $\widetilde e_k \subset \widetilde d_k$
or $\widetilde d_k \subset \widetilde e_k$. 
This is because $\widetilde e_k \subset \widetilde e_j$, \
$\widetilde d_k \subset \widetilde d_j$, \ 
$\widetilde d_k$, $\widetilde e_k$ intersect and
$d_k \cup \widetilde d_k$ does not intersect $\alpha$.
Assume wlog that $\widetilde e_k \subset \widetilde d_k$.
Choose $i > k > j$
with $\widetilde c_i \subset \widetilde e_k$ which is
a subset of $\widetilde d_k$ and hence of $\widetilde d_j$.
This proves the claim in this case.

Therefore by taking a bigger $j$ if necessary we can assume that
$Z$ has only one bi-infinite boundary component.
Let $y_1, y_2$ be the rays of $d_j$ and $z_1, z_2$ be the
rays of $e_j$. 
The bi-infinite component of 
$\partial Z$ has two rays which are contained in
$y_1 \cup y_2 \cup z_1 \cup z_2$.
If there are subrays of both rays in this boundary $\partial Z$
which
are contained in $y_1 \cup y_2$, then it follows that
$\widetilde d_j \cup d_j - (\widetilde e_j \cup e_j)$ 
is contained in a compact set in $\oo$, see
fig. \ref{inters}, a.
Since the decreasing sequence $(\widetilde d_k)_{k \in {\bf N}}$
of open sets in $\oo$
escapes compact sets in $\oo$, then there would be
$k$ with $\widetilde d_k \subset \widetilde e_j$.
But then there is $i$ with 
$\widetilde a_i \ \subset \ \widetilde d_k \ \subset \
\widetilde e_j$
and this would yield the claim in this case.

The remaining possibility to be analysed is that 
one and only one boundary ray of
$\partial Z$ must be
contained in $y_1 \cup y_2$ and one and only
one boundary ray of $\partial Z$ is in $z_1 \cup z_2$.
This last fact also implies that if a boundary
ray is contained in $y_1 \cup y_2$ then it cannot
have a subray in $z_1 \cup z_2$. 
The argument here will be to produce
two fixed rays $r, l$ of $\oos$ or $\oou$
which always intersect $\widetilde d_i, \widetilde e_i$ respectively
and so that $r, l$ are not equivalent.
This will contradict the fundamental lemma.

Let $l_j$ be the boundary ray of $Z$ contained
in $z_1 \cup z_2$. 
Then this ray is in $\widetilde d_j \cup d_j$ and
since it cannot have a subray contained in 
$d_j$ it follows that it has a subray contained in $\widetilde d_j$.
It also follows that the other ray of $e_j$ 
has to be eventually disjoint from 
$\overline Z$.
Similarly there is a ray $r_j$ of $d_j$ contained in
$\widetilde e_j$, see fig. \ref{inters}, b.
Recall that $\widetilde b_m \subset \widetilde d_j
\cap \widetilde e_j$.
Now consider $i \geq j$. 
If $\widetilde d_i \subset \widetilde e_j$ then we are done.
Otherwise

$$\widetilde d_i \cap \widetilde e_j \ \not = \ \emptyset 
\ \ \ {\rm and} \ \ \ 
\widetilde d_i \ \not \subset \ \widetilde e_j$$

\noindent
so the same analysis as above produces a ray of $e_j$ contained
in $\widetilde d_i$. It can only be $l_j \cap \widetilde d_i$ since the other
ray of $e_j$ is disjoint from $d_j \cup \widetilde d_j$,
so certainly disjoint from $d_i \cup \widetilde d_i$.
It now follows that for 
any $i \geq j$ 
there is a subray of
the fixed ray $l_j$ which is contained in
$\widetilde d_i$. Similarly for any $i \geq j$ there
is a subray of the fixed $r_j$ contained in $\widetilde e_i$.

\begin{figure}
\centeredepsfbox{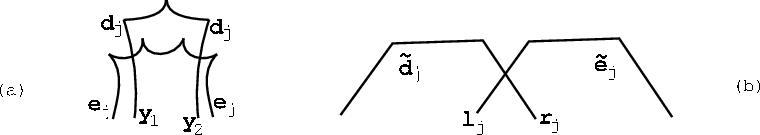}
\caption{
a. The intersection of convex neighborhoods,
b. Intersecting master sequences.}
\label{inters}
\end{figure}

The set $\widetilde d_j \cap \widetilde e_j$ has boundary
which contains subrays of $r_j, l_j$. 
If $r_j, l_j$ are equivalent rays then as there is $i$ with
$\widetilde b_i \subset \widetilde e_j \cap \widetilde d_j$,
the two rays of $b_i$ would be equivalent, contradiction.
Hence $r_j, l_j$ are not equivalent.
But for 
any $i \geq j$, then $\widetilde d_i \cup \widetilde e_i$ is a union
of two convex regions containing 
subrays of $l_j$ and $r_j$ ($j$ is fixed!). 
This is
disallowed by the fundamental lemma \ref{jumps}.
This proves the claim.

\vskip .1in
Suppose then there are infinitely many $j$'s so that for each
one of them,
there is $i(j) > j$ with $\widetilde a_{i(j)} \subset \widetilde e_j$.
Then for any $k$ there is one such $j$ with $j > k$ and so there is
$i(j) > j$ with $\widetilde a_{i(j)} \subset \widetilde e_j \subset
\widetilde e_k$. This means that $A \leq E$ and so $A \cong C$.
The claim shows that if this does not occur, then
there are infinitely many $j$ and
for each such $j$ there is $i(j) \geq j$
and $\widetilde c_{i(j)} \subset \widetilde d_j$. This now implies
that $C \leq D$ and again $C \cong A$. This finishes the proof
that $\cong$ is an equivalence relation.
\end{proof}

We first analyse admissible sequences associated to rays of $\oos$ or
$\oou$ $-$ each ray will define an ideal point of $\oo$. Later we define
general points of $\oo$. We will be interested in the asymptotic
behavior as points escape the ray to infinity. A ray does not separate
$\oo$, but still one can define sides of a ray as follows: let $l$
be a ray of (say) $\oos$. Fix a regular point $p$ in $l$ and consider
the component $W$ of $\oo - \oou(p)$ which contains a subray of
$l$. Then $l \cap V$ separates $V$ and we can talk about the sides
of $l$ in $V$. This depends only on the ray $l$ and not on the point $p$.

\begin{figure}
\centeredepsfbox{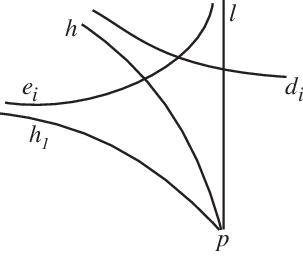}
\caption{
The process of creating standard sequences for rays of 
$\oos, \oou$.
Here the sequence $(d_i)$ of $\oou$ does not escape compact sets
and limits to $h$ leaf of $\oou$ making a perfect fit with
$l$. There is also the sequence $(e_i)$ of leaves of $\oos$ whose intersection
with $h$ escapes in $h$ and $(e_i)$ limits to a leaf
$h_1$ of $\oos$ making a perfect fit with $h$. 
The leaves $l, h_1$ are not separated from each other
in the leaf space of $\oos$.}
\label{hoo}
\end{figure}

\begin{figure}
\centeredepsfbox{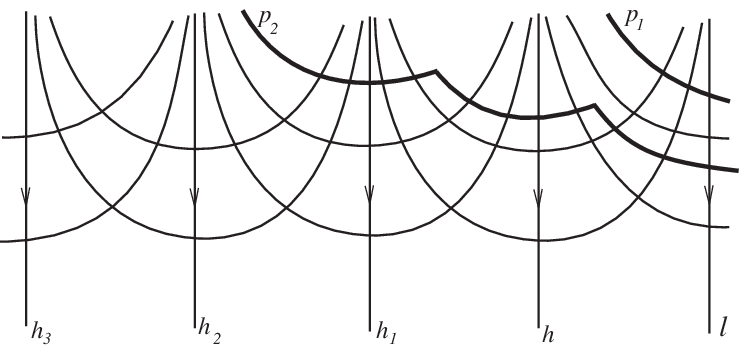}
\caption{
A picture of an infinite perfect fit or a perfect fit horoball.
Here $l, h_1, h_3$ are rays of $\oos$ and $h, h_2$ are rays of $\oou$.
The arrows indicate the direction of the rays.
$l$ and $h_1$ are not separated from each other in the leaf space
of $\oos$ and similarly for $h_1, h_3$ and also for $h, h_2$ (leaf
space of $\oou$ for the last 2). 
The figure is intended to continue indefinitely in both horizontal
directions.
The bold paths $p_1, p_2$ are 2 steps in producing a standard
sequence for the ray $l$.
$p_1$ is a polygonal path of length $1$ and $p_3$ is a polygonal
path of length $3$ (we are only describing what happens
in one side of $l$). }
\label{horo}
\end{figure}

\begin{define}{}{(standard sequences)}
Let $l$ be a ray in $\oos$ or $\oou$. For simplicity assume
that it is in $\oos$. Fix a side of $l$. Let $d_i$ be a nested sequence of
leaves of $\oou$ intersecting $l$ with $d_i \cap l$ escaping $l$.
If $d_i$ escapes compact sets in $\oo$ then $(d_i)$ is
an admissible sequence which is called a standard 
sequence associated to $l$.
If the $d_i$ do not escape in $\oo$, then they limit
on a collection of unstable leaves.
There is one of them, call it $h$
which makes a perfect fit with $l$ on the fixed side of $l$.
Consider now $e_i$ stable (nonsingular) leaves intersecting
$h$ and so that $h \cap e_i$ escapes compact sets in $h$ and
moves in the direction toward the perfect fit with $l$. 
Since $l$ and $h$ form a perfect fit, then for big
enough $i$, the $e_i$ and $d_i$ intersect and form a polygonal path
of length 2, see fig. \ref{hoo}. 

We want to produce an escaping polygonal sequence in that side of
$l$ and we already achieved that with $d_i \cup e_i$ for the region
between $l$ and $h$. Therefore we want to analyse what happens
beyond $h$, that is, the side of $h$ opposite to $l$ or not containing
$l$.
If the rays of  $e_i - h$ in the side of $h$ opposite to $l$
escape in $\oo$ then the polygonal paths made up of
a segment of $d_i$ and a ray of $e_i$ escape compact sets in $\oo$.
Otherwise the rays of $e_i - h$ on that side of $h$
limit to a stable leaf $h_1$ making
a perfect fit with $h$, see fig. \ref{hoo}. Notice that $h_1$ and $l$ are
not separated from each other in the leaf space of $\oos$ $-$ because the
sequence $(e_i)$ converges to both of these leaves.
Now iterate this process.
If this stops after finitely many steps 
then take a sequence of polygonal paths of
fixed length. Otherwise there are infinitely many leaves 
$h_j, j \geq 2$, alternatively in $\oos, \oou$, so that appropriate
rays of $h_j$ make a perfect fit with $h_{j-1}$ and $h_{j+1}$.
In this case 
use polygonal paths of increasing lengths, in order to 
cross over an increasing number
of perfect
fits emanating from $l$, see fig. \ref{horo}.
Do the same for the other side of $l$.
The ensuing sequence $(a_i)$ is an admissible sequence associated
to the ray $l$.
It is called a standard sequence for the ray $l$ of
$\oos$ or $\oou$.
\label{stan}
\end{define}

\vskip .1in
\noindent
{\bf {Remark}} $-$ If there are no perfect fits then $(d_i)$ as
in defintion \ref{stan} is a standard sequence for the ray $l$.

\vskip .1in
There are several other important remarks here and they 
concern only the case with perfect fits.
Along the way we will introduce the concepts of infinite perfect fits
and perfect fit horoballs.
First notice that standard sequences
 for a given ray $l$
are not unique. By construction it is easy to see that the $a_i$ are
convex, the rays of each $a_i$ are not equivalent to each other
 and the sequence $(a_i)$ is nested.
To check whether $\widetilde a_i$ is escaping:
If the $a_i$ have fixed length with $i$ then it is easy to see this.
Otherwise notice that the collection of rays equivalent to a given
ray escapes compact sets in $\oo$, in fact the whole leaves do.
That is, if $h_1, h_2, h_3, h_4... $ are the leaves  produced by the
construction in the definition, then $h_j, h_{j+2}$ are not separated from
each other in the respective leaf space.
Then the sequence $(h_i)$ escapes compact sets
in $\oo$.
So the sequence $(a_i)$ again escapes compact sets.
Hence $(a_i)$ is admissible.
In addition $h_i$ separates $h_k$ from $h_j$ for any
$k < i < j$.

\vskip .1in
\noindent
{\bf {Infinite perfect fits and perfect fit horoballs}}

In the case that the process above does not stop we call the infinite collection
of perfect fits an
{\em {infinite perfect fit}}. 
Associated to this 
one can define a model for a {\em {perfect fit horoball}} in $\oo$ as follows:
take the punctured square $[-1,1] \times [-1,1] - \{ 0, 0 \}$ with
its horizontal and vertical foliations and lift it to its universal
cover $U$. A proper, foliation respecting (horizontal goes to stable,
vertical goes to unstable) embedding of $U$ into $\oo$ gives an intuitive
``neighborhood"  of an ideal point associated to an infinite 
perfect fit as above. Such points clearly seem to have a ``parabolic"
feel as one suspects there is a covering translation which preserves
the perfect fit horoball and acts as a
translation in the collection of the perfect fits.
This is in analogy with Kleinian groups.

Two important questions arise: Is this possible for pseudo-Anosov flows? 
Also is there
a non trivial isotropy group of this infinite perfect fit structure and
why does it not contradict that the action of $\pi_1(M)$ is
cocompact?
First of all this phenomenon does happen, in fact there are several examples,
even for Anosov flows.
The first one is the seminal example of Franks and Williams \cite{Fr-Wi} of 
an intransitive Anosov flow in a closed $3$-manifold.
There is a simple picture of an infinite perfect fit in
the figure in page 164 of \cite{Fr-Wi}.
A second, also famous example, is that of the Bonatti-Langevin \cite{Bo-La}
example of a transitive Anosov flow with a transverse torus and not 
conjugate to a suspension. The structure in the universal cover of this
example is briefly described in \cite{Fe4}.

Once existence of infinite perfect fits is established, one wants to understand
its structure. Notice that infinite perfect fits have in particular
infinitely many pairs of leaves non separated from each other.
The author previously proved \cite{Fe4,Fe5} that up to covering
translations there are only finitely leaves of $\oos$ or
$\oou$ which are not separated from another leaf in the respective leaf space.
Hence given the collection $(h_j)$ produced above so that $h_j$ forms
a perfect fit with $h_{j+1}$, there are $j \not = k$ and $g$ in $\pi_1(M)$
so that $g(h_j) = h_k$.
This implies that the infinite sequence of perfect fits is in fact a
bi-infinite sequence $-$ that is, it extends indefinitely in the
other direction as well. It also justifies the terminology parabolic used above.
In addition if $z$ in $\oos$, $\oou$ is non separated from another leaf, then
the isotropy group of $z$ is non trivial \cite{Fe4,Fe5}. 
In particular this is true of every $h_j$.
With a little more work
this implies that associated to an infinite perfect fit there is a ${\bf Z}^2 \oplus
{\bf Z}^2$ subgroup of $\pi_1(M)$ which leaves the whole structure invariant.
Hence if $M$ is atoroidal, there can be no infinite sequence of perfect fits.

Finally, given the association of parabolic behavior with non compact manifolds,
how does this interact with the fact that $M$ is compact?
In the case of a hyperbolic $3$-manifold and a ${\bf Z} \oplus {\bf Z}$ cusp, then
geodesics escaping to the cusp are asymptotic. In the case of pseudo-Anosov flows,
suppose that leaves $l, h$ of $\oos$ and 
$\oou$ make a perfect fit. We need to analyse the situation in
$\mi$, not $\oo$. Let then (say) $L$ in $\wls$ which projects to $l$ in $\oo$
and similarly $H$ in $\wlu$ projecting to $h$. Then $L, H$ make a perfect fit.
But they are not asymptotic as points escape in $L$ or $H$. If they were, then
in fact $L$ and $H$ would intersect because of the local product structure
of $\ls, \lu$.
In particular $L, H$ would not form a perfect fit.
At this point it is useful to stress once more that the orbit space
$\oo$ is a topological and dynamical object, but it is {\underline {not}} a
metric object. 
Even
though topologically it may seem that rays of $\oos, \oou$ making a perfect
fit are getting close,
this can only be checked in $\mi$, where in fact one sees that their lifts are
not getting close.

\begin{lemma}{}{}
Let $l$ be a ray in $\oos$ or $\oou$ and let $C = (c_i)$ be
a standard sequence associated to $l$.
Let $A = (a_i)$ be an admissible sequence so that for any
$i$, then $\widetilde a_i \cup a_i$ contains a ray
equivalent to $l$.
Then $A \leq C$.
\label{std}
\end{lemma}

\begin{proof}{}
Suppose the lemma is not true and fix an $i$ so that for any
$j$, $\widetilde a_j \not \subset \widetilde c_i$.
Notice first that by the definition of a standard sequence, then
for any $m$ (in particular for $m = i$) and for any ray $s$ 
equivalent to $l$, then $s$ has a subray $s'$ contained
in $\widetilde c_m$.
Since for any $j$, $\widetilde a_j \cap a_j$ contains
such a ray $s$ then $\widetilde a_j \cap \widetilde c_i
\not = \emptyset$. 
If in addition $\widetilde a_j \not \subset \widetilde c_i$, then
as seen in the fundamental lemma, for $j$ big enough, there is
at least one ray of $c_i$ which has a subray contained
in $\widetilde a_j$. By the fundamental lemma,
after discarding finitely many terms in
$(a_j)$ there is a fixed ray $r$ of $(c_i)$
which for every $j$ has a subray contained in $\widetilde a_j$.
Notice that $r$ and $l$ are not equivalent.
We conclude that every $\widetilde c_j$ contains a subray of the
fixed ray $l$
and every $\widetilde a_j$ contains a subray of the
fixed ray $r$.
Since for any $j, m$, \ $\widetilde c_j \cap \widetilde a_m
\not = \emptyset$ this is disallowed by the fundamental lemma.
This finishes the proof of the lemma.
\end{proof}

We now define ideal points of $\oo$.

\begin{define}{}{}
Suppose that $\Phi$ is not topologically conjugate to a suspension
Anosov flow.
A point in $\partial \oo$ or an ideal point of $\oo$
is an equivalence class of
admissible sequences of polygonal paths.
Let $\cd = \oo \cup \partial \oo$.
\label{idpo}
\end{define}

Given $R$, an admissible sequence of polygonal paths, let
$\overline R$ be its equivalence class under $\cong$.
Notice that each ray $l$ in $\oos, \oou$ has admissible sequences 
and the these sequences are all equivalent.
In this way $l$ defines a single point in $\partial \oo$ which 
is denoted by $\partial l$.
This is generalized in the following way: if $l$ is a leaf of $\oos$
or $\oou$, then we denote by $\partial l$ the collection of
ideal points of rays of $l$.
If $l$ is a ray of $\oos, \oou$ associated to infinite perfect fit
then $\partial l$ is called a parabolic ideal point in $\partial \oo$.
We will see later that in this case $\partial l$ is the unique
fixed point of the action of some $g$ in $\pi_1(M)$ which acts
in a ``parabolic" way on $\partial \oo$.

\begin{define}{(master sequences)}{}
Let $R$ be an admissible sequence. An admissible sequence $C$
defining $\overline R$ is a master sequence for $\overline R$
if for any $B \cong R$, then $B \leq C$.
\label{master}
\end{define}

Why master sequences? Ideal points are defined by admissible sequences
of polygonal paths 
and not by sequences of points in $\oo$. Given the admissible sequence $(a_i)$
defining an ideal point $p$, one intuitively expects that a fixed 
$\widetilde a_i$ will at least limit on all points of $\partial \oo$
near $p$ (the topology in $\oo \cup \partial \oo$ will be defined
formally later).
However this is not the case. For example given $l$ a ray in $\oos$ with no
perfect fits associated to it, consider a sequence of regular leaves
$d_i$ in $\oou$
with $d_i \cap l$ escaping in $l$. Then $(d_i)$ defines the 
ideal point $\partial l$. Now fix a
side of $l$ and consider the rays of $d_i - l$ in
this side of $l$. For each $i$, this ray, together with 
an appropriate subray of $l$
forms a convex polygonal path
$b_i$ and
defines an admissible sequence $(b_i)$.
Intuitively $\widetilde b_i$
is $\widetilde d_i$ cut in half by a ray of $l$.
Clearly $(d_i)$ and $(b_i)$ are equivalent, so $(b_i)$ also defines
the same ideal point. But a fixed
$\widetilde b_i$ only accumulates on one side
of $l$. The master sequences are those $(d_i)$ for which an individual
$\widetilde d_i$ ``limits on both sides"
of the ideal point.


\vskip .1in
\noindent
{\bf {Remark}} $-$ 
Recall that a {\em {cyclic order}} on a set $B$ is a partition of the set of 
pairwise distinct triples $(p,q,r)$ into two sets, called the
``positive and negative triples", such that cyclic permutations
in $(p,q,r)$ preserve the sign, non cyclic permutations reverse
the sign and if $(p,q,r)$ and $(r,s,p)$ are positive
triples, then $(q,r,s)$ is also a positive
triple.
\vskip .05in

\begin{define}{}{(order of sets in $\oo$)}
Let ${\cal C} = \{ c_i \}, i \in I \subset {\bf Z}$ be a collection
of properly embedded bi-infinite arcs in $\oo$ 
so that there are components $\widetilde c_i$ of $\oo - c_i$
with $\{ c_i \cup \widetilde c_i \}$
pairwise disjoint.
Suppose that ${\cal C}$ is localy finite: any
compact set in $\oo$ intersects only finitely many of the
$c_i$.
Fix $x \in \oo$ not in any $c_i \cup d_i$ and choose 
paths $\gamma_i$ from $x$ to $c_i$ which are pairwise
disjoint except for $x$. This is all possible since $\oo \cong
\rrrr^2$. Then the germs of the collection $\{ \gamma_i \}$ at
$x$ put a cyclic order in the collection $\{ \gamma_i \}$
and hence on
${\cal C}$. This order is independent of $x$ or the paths
$\gamma_i$. If all $\widetilde c_i$ miss a fixed properly embededded
infinite arc $\gamma$ starting at $x$, then there is a 
linear order in ${\cal C}$. The linear order depends on the path
$\gamma$.
\label{order1}
\end{define}

\begin{lemma}{}{}
Given an admissible sequence $R$,
there is a master sequence
for $\overline R$.
\label{mast}
\end{lemma}

\begin{proof}{}
\noindent
{\underline {Case 1}} $-$ Suppose that for any $A = (a_i), \
B = (b_i)$ in 
$\overline R$ and for any $i, j$ then $\widetilde a_i \cap \widetilde
b_j \not = \emptyset$.

We claim that in this case any $A \cong R$ will serve as a master
sequence.
That is we do not have the situation described above were one slices through
the admissible regions using a fixed ray of $\oos$ or $\oou$.
Choose $A \cong R$ and let $B \cong R$.
We want to show that $B \leq A$.
So by way of contradiction, 

$${\rm assume \  that \ 
there \ is } \ \ i \ \ {\rm  so \ that \
for \ any} \ \  j, \ \ 
\widetilde b_j \ \not \subset \ \widetilde a_i \ \ \ (*)$$

\noindent
This also works for any $k \geq i$, but we will fix $i$ from now on
in case 1.
The contradiction will be obtained by first showing that 
$(*)$ implies that $A$ is associted to an ideal point of a ray
of $\oos$ or $\oou$ and then producing two admissible sequences
in $\overline R$ which fail the hypothesis of case 1.

In case 1, 
$\widetilde b_j \cap \widetilde a_i$ is not empty
for any $j$.
Let $u, v$ be the rays of $a_i$.
Since $\widetilde b_j$ escapes compact sets in $\oo$
as $j \rightarrow \infty$, so does $\widetilde b_j
\cap \widetilde a_i$. 
The arguments of lemma \ref{equiv}, referring to figure \ref{inters}, a;
show that $\widetilde a_i \cup a_i$ 
cannot contain subrays of both rays in $b_j$ and 
in fact 
for $j$ big enough, then 
$\widetilde b_j$ contains at least one  subray $u$ or $v$ and no singular
point.
This implies that $a_i$ cuts $\widetilde b_j$
into at most 3 non compact regions (all of which are convex):
at most one region contained in $\widetilde a_i$ and
at least one and at most 2 disjoint from $\widetilde a_i$.
The regions are convex because one can assume 
$j$ is big enough so that the $b_j \cap a_i$
does not contain any singularity.
Up to discarding finitely many terms we may assume that 
one region contains in its boundary a subray of (say) $u$.
Call this region $\widetilde c_j$ with boundary $c_j$.

There are 2 possibilities: i) For $j$ big enough, the region
$\widetilde c_j$
disappears, that is,
there is no such region with a subray $u$ in
the boundary. In that case there is another  region $\widetilde d_j$
of $\widetilde b_j$ cut along $a_i$
disjoint from $\widetilde a_i$ and containing a subray of $v$ in
the boundary. If $\widetilde d_j$ also eventually disappears,
then some $\widetilde b_k$ is contained in $\widetilde a_i$, 
contrary to assumption in this argument.
So at least one of $(\widetilde c_j)$, $(\widetilde d_j)$ is 
always non empty.
This reduces to the following:
ii) (say) $\widetilde c_j$ is never empty for
any $j$. Then 
$\widetilde b_j$ contains a subray of $u$ for any $j$.
Let $E = (e_k)$ be a standard sequence associated with the ray
$u$. Eliminating finitely many initial terms of $E$ if necessary
we can assume that $u$ cuts every $\widetilde e_k$ into
two components $\widetilde f_k$ and $\widetilde g_k$, which
are convex, with boundaries $f_k$ and $g_k$ respectively and
defining admissible sequences $F = (f_k)$ and $G = (g_k)$.
Assume that $\widetilde f_k \cap \widetilde a_i
= \emptyset$ for all $k$.
Clearly $F \leq E, G \leq E$ and $\widetilde f_k \cap \widetilde g_k
= \emptyset$. 

Suppose that for some 
$m > i$, \ $a_m$ does not have a ray 
equivalent to $u$. Fix this $m$. 
Notice that $\widetilde b_j$ contains a subray of
a fixed ray of $a_m$ and also a fixed subray of $u$
(this is a ray of $a_i$ with $i$ fixed). This is now disallowed
by the fundamental lemma.

The remaining possibility in this case 
is that $a_m$ always has a ray equivalent
to $u$ for any $m$. By 
lemma \ref{std} it follows that $A \leq E$ and so
$R \cong A \cong E \cong F \cong G$. Hence in $\overline R$ there are 
$F = (f_k)$, $G = (g_k)$ with $\widetilde f_k \cap \widetilde g_k 
= \emptyset$ for some $k$. This contradicts the
hypothesis in case 1 and 
implies that $A$ is a master sequence for $\overline R$.


\vskip .1in
\noindent
{\underline {Case 2}} $-$ There are $A, B$ in $\overline R$
and $i$ so that $\widetilde a_i, \ \widetilde b_i$ are disjoint.

Fix this $i$.
In particular $\widetilde a_k, \widetilde b_k$ are disjoint
for $k \geq i$.
Let $C$ be an admissible sequence with $A \leq C, \ B \leq C$. 
We claim that $C$ is a master sequence for the class 
$\overline R$.
Let $D \cong A$. Suppose that $D \not \leq C$.
Hence there is $m$ so that $\widetilde d_j
\not \subset \widetilde c_m$ for any $j$. Fix this $m$.
There are two options: i) There is $k$ with
$\widetilde d_k \cap \widetilde c_k = \emptyset$ ,
\ or \ 
ii) For any $k$, $\widetilde d_k \cap \widetilde c_k \not
= \emptyset$, in which case $\widetilde d_k \cap \widetilde
c_j \not = \emptyset$ for any $k, j$.

In subcase i)
up to deleting a few initial terms we may assume that
$\widetilde d_1 \cap \widetilde c_1 = \emptyset$.
We have $A \cong B \cong D$ with 
$\widetilde a_i, \widetilde b_i,
\widetilde d_i$ disjoint.
Choose $E = (e_j)$ with $C \leq E, D \leq E$.
Assume for simplicity that $i$ is big enough so that 
$\widetilde a_i,
\widetilde b_i,
\widetilde d_i$ are contained in $\widetilde e_1$..
This puts a linear order in $a_i, b_i, d_i$ and we can assume without
loss of generality that $b_i$ is between $a_i$ and $d_i$.
Since 
$b_i$ is between $a_i$ and $d_i$ then:
for any $j$, $\widetilde e_j$ contains subrays of
the
rays of $b_i$ (with $i$ fixed!), which are not equivalent. 
The fundamental lemma \ref{jumps} implies this is impossible.

We now consider option ii).
Since $\widetilde a_i \cup a_i$ and
$\widetilde b_i \cup b_i$ are disjoint and $A \leq C, B \leq C$,
then there is a ray $u$ of $a_i$ and a ray $v$ of $b_i$ so that
for any $j$, $\widetilde c_j$ contains subrays of $u$ and $v$.
A priori $u, v$ can be equivalent.
Since $\widetilde d_j$ is not contained in $\widetilde c_m$
but has to intersect $\widetilde c_m$, we may assume up to eliminating
a few initial terms that $\widetilde d_j$ always contains a subray
of a fixed ray $y$ of $c_m$. 
The rays $y, u$ are not equivalent.
Since $\widetilde d_j \cap \widetilde c_k
\not = \emptyset$ for any $k, j$,
this contradicts the fundamental lemma.
So this cannot happen either.

We conclude that $C$ is a master sequence for $\overline R$.
This finishes the proof of lemma \ref{mast}.
\end{proof}

By definition for any 2 master sequences $A, B$ in the class
class $\overline R$, it follows that
both $A \leq  B$ and $B \leq  A$ hold.

\begin{lemma}{}{}
Let $p, q$ in $\partial \oo$. Then $p, q$ are distinct if and
only if there are master sequences 
$A =  (a_i), B =  (b_i)$ associated
to $p, q$ respectively with 
$(a_i \cup \widetilde a_i) \cap (b_j \cup \widetilde b_j)
= \emptyset$ for some $i, j$. 
\label{idst}
\end{lemma}

\begin{proof}{}
We first show that $p, q$ are distinct if and only if
there are master sequences $A = (a_i), B = (b_i)$, so
that for some $i, j$, $\widetilde a_i \cap \widetilde b_j = \emptyset$.
In the proof we show that the negations are equivalent.
First suppose that $p = q$. Let $A, B$ be any master
sequences associated to $p= q$. Then since $A, B$ are master
sequences associated to the same equivalence class
then $A \leq B$ and $B \leq A$. 
Therefore we can never have $\widetilde a_i \cap \widetilde b_j 
= \emptyset$.
This is the easy implication.

To prove the converse, 
suppose that for any master sequences $A = (a_i)$
and $B =  (b_i)$ associated to $p, q$ respectively and
any $i, j$ then $\widetilde a_i \cap \widetilde b_j \not =
\emptyset$. 
Let $A, B$ be such a pair.
Suppose first that for all $i$, $\widetilde a_i \cap 
\widetilde b_i$ has 2 non compact components.
Then an argument similar to one in the proof of lemma
\ref{mast} shows that there are non equivalent
rays $u, v$ with subrays contained in each $\widetilde a_i
\cap \widetilde b_i$. This is disallowed by the
fundamental lemma.
Similarly if $\widetilde a_i \cap \widetilde b_i$ has
a component with 4 boundary rays for infinitely
many $i$.
On the other hand, $\widetilde b_i \cap \widetilde a_j$ can never
be contained in a compact set or else for some $i' > i$ then
$\widetilde a_j \cap \widetilde b_{i'}  = \emptyset$.
One concludes that $\widetilde a_i \cap \widetilde b_i$ 
eventually has a single non compact component.
Let 
$\widetilde c_i$ be this component of
$\widetilde a_i \cap \widetilde b_i$ and let $c_i =
\partial \widetilde c_i$. Let $C = (c_i)$. 
Clearly $c_i$ is convex and $(c_i)$ is nested.
But a priori , $C$ may not be admissible, that is, the boundary
rays may be equivalent.
Notice that the rays in $c_i$ are subrays of rays of $a_i$ or $b_i$.

The first case is that the rays of $c_i$ are
not equivalent for any $i$. 
Then $c_i$ is a convex
polygonal path, non empty and $C$ is an admissible sequence.
Also $C \leq A, \ C \leq B$, which implies 
that $A \cong C \cong B$ and hence $p = q$.

The second case is
that there is $i$ so that the rays $u, v$ of $c_i$ are
equivalent. 
Notice this can only happen if there are perfect fits.
There is a collection ${\cal Y} = \{ u_0 = u, u_1,..., u_n = v \}$
of rays of $\oou, \oos$
so that $u_k, u_{k+1}$ make a perfect fit for every $k$.
Since the sequence $(\widetilde c_j)$ is nested with $j$, 
the rays of $c_j$ for
$j > j_0$ have to be in the collection ${\cal Y}$.
Up to subsequence we can assume they are all subrays of fixed
rays $r, l$.
Notice that $r \not = l$, or else 
$\widetilde b_j \cap \widetilde a_j
= \emptyset$ for some $j > i$.
Since $r, l$ are equivalent they cannot both be rays of $a_j$
(or both of $b_j$ either).
Hence up to renaming objects, $a_j$ has a subray in $r$ 
and $b_j$ has a subray in $l$, for all
$j > i$, see fig. \ref{renu}.

\begin{figure}
\centeredepsfbox{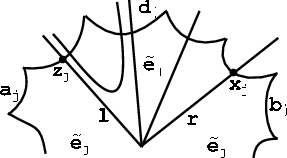}
\caption{
Interpolating chains that intersect
to produce a new convex chain.}
\label{renu}
\end{figure}

Let $z_j = a_j \cap l$, \ $x_j = r \cap b_j$. 
As in the proof of the fundamental lemma notice that $\widetilde b_j$
contains a subray of $r$ and $\widetilde a_j$ contains a subray of $l$.
Then $z_j$
escapes in $l$ and $x_j$ escapes in $r$.
Let $a'_j$ be the component of $a_j - z_j$ not containing a subray of $r$ and
$b'_j$ the component of $b_j - x_j$ not containing a subray of $l$.
The above implies that we can 
connect $z_j, x_j$ 
by a finite convex
polygonal path $d_j$ which extends 
$a'_j \cup b'_j$ to a convex polygonal path $e_j$.
see fig. \ref{renu}.
This is because 
 $l, r$ are
connected by finitely many perfect fits. If $z_j, x_j$ are 
very deep in the rays $l, r$ then we can always connect
$z_j$ and $x_j$ by a convex polygonal path.
Notice that $a_j$ has a subray of $r$ so it goes to $r$,
but $a_j$ may reach $r$ in a point different than $x_j$.
If we just connect this to $x_j$ and then follow
along $b_j$ this will produce a non convex switch in $r$.
That is why we use the interpolating polygonal path $d_j$.
Then the polygonal paths $e_j$ are convex and one can construct
the interpolating polygonal path $d_j$ so that $e_j$ escapes
compact sets as $j \rightarrow \infty$. 
Then $E = (e_j)$ defines an admissible sequence
of chains. Clearly $A \leq E$ and $B \leq  E$
so that $A \cong B$ and again $p = q$.

This finishes the equivalence with the intersection condition
on open sets. 
Finally suppose that $\widetilde a_i \cap \widetilde b_i = \emptyset$
for all sufficiently big $i$, but
$(a_i \cup \widetilde a_i)
\cap (b_i \cup \widetilde b_i) \not = \emptyset$ for any $i$.
This can only happen if there is a ray $l$ of $\oos$ or $\oou$ so that
both $a_i$ and $b_i$ have a subray of $l$.
Let $C = (c_i)$ be a standard sequence for $l$.
By lemma \ref{std} \ $A \leq C$ and $B \leq C$, so $A \cong B$ and
$p = q$.
This proves the lemma.
\end{proof}

We now define the topology in $\cd = \oo \cup \partial \oo$.

\begin{define}{(topology in $\cd = \oo \cup \partial \oo$)}{}
Let ${\cal T}$ be the set of subsets 
$U$ of $\cd = \oo \cup \partial \oo$
satisfying the following two conditions:

(a)  $U \cap \oo$ is open in $\oo$.

(b) If $p$ is in $U \cap \partial \oo$
and $A = (a_i)$ is any master sequence associated
to $p$, then there is $i_0$ satisfying
two conditions: (1) $\widetilde a_{i_0} \subset U
\cap \oo$ and (2) For any $z$ in $\partial \oo$,
if it admits a master sequence $B = (b_i)$ so that
for some $j_0$, one has $\widetilde b_{j_0} \subset 
\widetilde a_{i_0}$ then $z$ is in $U$.
\label{topology}
\end{define}

First notice that if the second requirement works
for a master sequence $A = (a_i)$ with index $i_0$, then
for any other master sequence $C = (c_k)$ defining $p$,
we can choose
$k_0$ with
$\widetilde c_{k_0} \subset \widetilde a_{i_0}$. 
Then $\widetilde c_{k_0} \subset U$.
A point $q$ of $\partial \oo$ which has a master
sequence $B = (b_j)$  and $j_0$ so that

$$\widetilde b_{j_0} \
\subset \ \widetilde c_{k_0}; \ \ \ {\rm then } \ \ \ 
\widetilde b_{j_0} \ \subset \ \widetilde a_{i_0}$$

\noindent
so
$q$ is in $U$. Therefore (b) works for $C$ instead of $A$
with 
$k_0$ instead of $i_0$. So we only need to check
the requirements for a single master sequence.

\begin{lemma}{}{}
${\cal T}$ is a topology
in $\cd = \oo \cup \partial \oo$.
\label{topol}
\end{lemma}

\begin{proof}{}
Clearly $\cd, \emptyset$ are in ${\cal T}$.
Unions:
If \ $\{ V_{\alpha} | \alpha \in {\cal A} \}$  \ is a family 
of sets in 
${\cal T}$, then let $V$ be their union.
If $x$ is in $V$ and $x$ is in $\oo$, there is open
set $O$ in $\oo$ with $x \in O \subset V_{\alpha}
\subset V$ for
some index $\alpha$, hence satisfying condition (a).
Let now $p$ in $V \cap \partial \oo$.
There is $\beta \in {\cal A}$ with $p \in V_{\beta}$.
Let $A = (a_i)$ be a master sequence associated to $p$.
There is $i_0$ with 

$$\widetilde a_{i_0} \ \subset V_{\beta} \cap \oo \ 
\subset V \cap \oo \ \subset \ \oo.$$

\noindent
In addition if $q \in \partial \oo$ and $q$ has
a master sequence $B = \{ b_j \}$ and $j_0$ 
with $\widetilde b_{j_0} \subset \widetilde a_{i_0}$
then $q$ is in $V_{\beta} \subset V$. Hence this $i_0$ works
for $V$ as well. This proves that $V$ is in ${\cal T}$.

Intersections:
Let $V_1, V_2$ be in ${\cal T}$
and $V = V_1 \cap V_2$. Clearly $V_1 \cap V_2 \cap \oo$ is open
in $\oo$. Let $u \in V_1 \cap V_2 \cap \partial \oo$.
Given a master sequence $A = (a_i)$ associated to $u$ there
is $i_1$ with $\widetilde a_{i_1} \subset V_1$ and
if $q$ has master sequence $B = (b_j)$ with 
$\widetilde b_{j_0} \subset \widetilde a_{i_1}$ 
then $q$ is in $V_1$.
Similarly considering $u \in V_2$, there is index
$i_2$ satisfying the conditions for $V_2$. Let $i_0 = max(i_1, i_2)$.
Then $\widetilde a_{i_0}$ is contained in $V_1$ and $V_2$
(since $\widetilde a_i$ are nested).
If now $q$ is in $\partial \oo$ has a master sequence 
$B = (b_j)$ with
$\widetilde b_{j_0} \subset \widetilde a_{i_0}$ for some $j_0$ 
then $q$ is in $V_1$ and is in $V_2$ by choice of $i_1, i_2$.
Therefore $q$ is in $V$. Hence $V$ is in ${\cal T}$.
This shows 
${\cal T}$ is a topology in $\oo \cup \partial \oo$.
\end{proof}

\vskip .1in
\noindent
{\bf {Action of $\pi_1(M)$ on $\cd = \oo \cup \partial \oo$}}

One key remark is that the action of $\pi_1(M)$ on $\oo$ 
preserves the foliations $\oos, \oou$ and sends convex polygonal paths to
convex polygonal paths. If follows that $\pi_1(M)$ acts by homeomorphisms
on $\cd$.

\begin{lemma}{}{}
Suppose $\pi_1(M)$ preserves orientation in $\oo$. 
Then $\partial \oo$ has a natural cyclic order.
\label{order}
\end{lemma}

\begin{proof}{}
Let $p, q ,r $ in $\partial \oo$ pairwise distinct points.
By lemma \ref{idst}, there are master
sequences $A = (a_i), B = (b_i), C = (c_i)$
associated to $p, q, r$ respectively with
$a_1 \cup \widetilde a_1, b_1 \cup \widetilde b_1, 
c_1 \cup \widetilde c_1$ pairwise disjoint.
By definition \ref{order1} there is a cyclic order on $a_1, b_1, c_1$.
This defines a cyclic order on $p, q, r$.
This is independent of the choice of
master sequences (since they are all equivalent).
This order is 
also invariant under
the action of $\pi_1(M)$ on $\oo$, 
since $\pi_1(M)$ preserves orientation in $\oo$.
This defines
a natural cyclic order in $\partial \oo$.
\end{proof}

In general let ${\cal E}$ be the index 2 subgroup of $\pi_1(M)$ 
preserving orientation of $\oo$.
Then ${\cal E}$ preserves a cyclic order in $\partial \oo$
and the elements in $\pi_1(M) - {\cal E}$ reverse this cyclic order.

In any case pick one orientation in $\oo$ that defines a cyclic
order in $\partial \oo$ (invariant only under ${\cal E}$).

\begin{define}{}{}
(the set $U_c$)
For any convex polygonal path $c$ there is an associated open set $U_c$
of $\cd$:
let $\widetilde c$ be the corresponding convex set of $\oo$ (
if $c$ has length 1 there are two possibilities). Let

$$U_c \ = \ \widetilde c \ \cup \ \{ x \in \partial  \oo \ | \
\ \ {\rm there \ is \ a \  master \ sequence} \ \ A = (a_i)
\ \ {\rm with} \ \ \widetilde a_1 \subset \widetilde c \ \}$$

\label{canon}
\end{define}

It is easy to verify that $U_c$ is always an open set in $\cd$.
In particular it is an open neighborhood of any point
in $U_c \cap \partial \oo$.
The rays of $c$ are equivalent if and only if $U_c$ is contained
in $\oo$.
The notation $U_c$ will be used from now on.

Given a cyclic order in $\oo$ and $p, q$ distinct in $\partial \oo$,
let 

$$(p,q) \ := \ \{ \ x \in \partial \oo \ | \ (p,x,q) \ 
{\rm is \ positive \ in \ the \ cyclic \ order \ of} \  \oo \ \}.$$

\noindent
This is the interval from $p$ to $q$ in the cyclic order.
Notice that if one changes the cyclic ordering then $(p,q)$ of the
new cyclic order is $(q,p)$ of the old cyclic order. So the 
collection of intervals is independent of the order.
Let ${\cal Z}$ be the topology in $\partial \oo$ generated by
the intervals. Given $t$ in $(p,q)$ there is a master sequence
$A = (a_i)$ for $t$ \ with \ $U_{a_1} \cap \partial \oo \subset (p,q)$. Hence 
$(p,q)$ is open in the topology of $\partial \oo$.
Conversely if $T$ is open in $\partial \oo$ and $t \in T$,
there is a master sequence $A = (a_i)$ satisfying property
(b) of definition of the topology in $\partial \oo$, so 
that $U_{a_1} \cap \partial \oo \subset T$. The endpoints of the rays
of $a_1$ are $p,q$ and then $t \in (p,q) \subset T$.
So the interval topology is exactly the induced topology
in $\partial \oo$.

\begin{lemma}{}{}
$\cd$ is Hausdorff.
\label{haus}
\end{lemma}

\begin{proof}{}
Any two points in $\oo$ are separated from each other. If $p, q$ are distinct
in $\partial \oo$ choose master sequences $A = (a_i)$
and $B = (b_i)$, where $\widetilde a_1 \cap \widetilde b_1
= \emptyset$. Let $U_{a_1}$ be the open set of $\cd$
associated to
$a_1$ and $U_{b_1}$ associated to $b_1$.
By definition $U_{a_1}$ is an open neighborhood of $p$ and
likewise $U_{b_1}$ for $q$. They are disjoint open sets
of $\cd$.

Finally if $p$ is in $\oo$ and $q$ is in $\partial \oo$,
choose $U$ a neighborhood of $q$ coming from a master sequence 
as above so that
$U \cap \oo$ does not have $p$ in its closure - always
possible because master sequences are escaping sets.
Hence there are disjoint neighborhoods of $p$, $q$.
\end{proof}

Our goal is to show that $\partial \oo$ is homeomorphic
to ${\bf S}^1$ and that $\cd$ is homeomorphic
to a closed disk. We need a few simple results:

\begin{lemma}{}{}
For any ray $l$ of $\oos$ or $\oou$, there is an associated
point in $\partial \oo$.  Two rays generate the same point
of $\oo$ if and only if the rays are equivalent (as rays!).
\label{equiva}
\end{lemma}

\begin{proof}{}
Given a ray $l$ any standard sequence $(c_i)$ associated
to it
defines a point in $\partial \oo$.
Let $r, l$ be rays of $\oos, \oou$.
If they define the same point of $\partial \oo$, then
there is a master sequence $C = (c_i)$ for 
this point. 
Since both standard sequences associated to $r, l$ are
$\leq C$, it follows that 
 every $\widetilde c_i$ contains 
subrays of both
$l, r$.
By the fundamental lemma
(where we use $E = F = C$ in that lemma),
this occurs if and only
if the rays $r, l$ are equivalent.
\end{proof}


\begin{lemma}{}{}
Suppose that $A = (a_i)$ is an admissible sequence 
of polygonal paths and that every 
$a_i$ contains a subray of a fixed ray $l$ of $\oos$ or $\oou$.
Then $A$ is associated to the ideal point $\partial l$ of $l$ and
$A$ is not a master sequence for the point  $\partial l$ of $\partial \oo$.
\end{lemma}

\begin{proof}{}
The point $\partial l$ was defined just before definition \ref{master}.
The first statement was proved in lemma \ref{std}.
For the second statement, notice that each $\widetilde a_i$
is contained in a fixed component of $\oo - l$.
Choose a standard sequence $B$ associated to $l$
and  cut it along $l$. Let $C$ be the admissible sequence
produced so that $\widetilde c_1 \cap \widetilde a_1
= \emptyset$. This shows that $A$ is not a master sequence
for $\partial l$.
%
\end{proof}

\begin{lemma}{}{}
Let $A = (a_i)$ be an admissible sequence defining a point
$p$ in $\partial \oo$. Then one of the following mutually
exclusive possibilities occurs:

(i) There are infinitely many $i$ in 
${\bf N}$ and for each such $i$ there is a ray $l_i$ of $a_i$
which is equivalent to a fixed ray $l$ of $\oos$ 
or $\oou$. Then $p$ is the ideal
point of any of the $l_i$ and $A$ is not a master sequence for $p$.
In fact in this case the hypothesis is true for any $i$ sufficiently
big.

(ii) There are only finitely many rays of paths in the collection
$\{ a_i \}$ which are equivalent
to any given ray of $\oos$ or
$\oou$. In this case $A$ is a master sequence for $p$.
\label{choi}
\end{lemma}

\begin{proof}{}
Most of part (i) was proved in lemma \ref{std}.
The $\widetilde a_i$ are nested and hence  the rays of $a_i$ are split into
two sequences $(r_i), (l_i)$ each of which 
is also ``nested". It is easy to check
that only elements of one of the sequences can be equivalent to $l$.
But if (say) $r_i$ and $r_j$ (with $j > i$) are both equivalent to $l$,
then $r_k$ is equivalent to $l$ for any $i < k < j$. Hence the $r_i$ are
equivalent to $p$ for any sufficiently big $i$. It does not follow 
however that for any big $i, j$, \ $r_i$ and $r_j$ share a subray. This is because 
there may be an infinite perfect fit, so the rays $r_i$ can change with
$i$
escaping in the horoball model of an infinite perfect fit. 
Finally a standard sequence for the ray $l$ and cutting shows that $A$
is not a master sequence for $\partial l$.
This proves (i).
%

To prove part (ii), let $A = (a_i)$ be an admissible
sequence so that 
there are only finitely many rays of $(a_i)$ which are equivalent
to any given ray of $\oos$ or $\oou$.
Suppose by way of contradiction that 
$A$ is not a master sequence  for $p$,
so there is $B \cong A$ 
and $B \not \leq A$. 
Fix some $n$ so that for no $j$, \
$\widetilde b_j \subset \widetilde a_n$.
Hence this is true for any $n' > n$.

The first possibility is 
there are $i, j$, with 
$\widetilde b_j \cap \widetilde a_i = \emptyset$. let
$E = (e_k)$ be an admissible sequence with $A \leq E, B \leq E$.
Choose $k > i, j$, hence $\widetilde b_k \cap \widetilde a_k 
= \emptyset$ and so that 
$a_k$ does not have any rays equivalent to 
any rays of $a_i$. Then any $\widetilde e_m$, $m \geq k$
 contains a fixed 
subray of $b_k$ and a fixed subray of $a_k$ and they are 
not equivalent by choice of $k$.
This is disallowed by the fundamental lemma.

The second possibilithy is that
$\widetilde b_j \cap \widetilde a_i \not
= \emptyset$ for any $i, j$.
Fix $k > n$ so that $a_k$ does not have any 
ray equivalent to a ray of $a_n$.
If the 2 rays of $b_j$ have subrays 
contained in $\widetilde a_n \cup a_n$
then $\widetilde b_j - (\widetilde a_j \cup a_j)$ is contained
in a compact set of $\oo$ and as seen before this implies
that for some $t > j$, then $\widetilde b_t \subset 
\widetilde a_n$, contrary to choice of $n$ in part (ii).
We conclude that for any sufficiently big $m$,
$\widetilde b_m \cup b_m$ is not contained in
$\widetilde a_n \cup a_n$ but has to 
intersect $\widetilde a_k$. 
This implies that for big $m$, \ 
$\widetilde b_m$ has to contain a subray of
a ray of $a_n$ and a subray of a ray of $a_k$.
Again this is disallowed by the fundamental lemma.
This finishes the proof of the lemma.
\end{proof}


\begin{lemma}{}{}
The space $\cd$ is first countable.
\label{fc}
\end{lemma}

\begin{proof}{}
Let $p$ be a point in $\cd$.
The result 
is clear if $p$ is in $\oo$ so suppose that $p$ is in $\partial \oo$.
Let $A = \{ a_i \}$ be a master sequence associated to $p$.
We claim
that $\{ U_{a_i}, i \in {\bf N} \}$ is a neighborhood basis at $p$.
Let $U$ be an open set containing $p$. By definition
\ref{topology} there is 
$i_0$ with $\widetilde a_{i_0} \subset U$
and if $z$ in $\partial \oo$ admits a master
sequence $B = (b_i)$ so that for some $j_0$
then $\widetilde b_{j_0} \subset \widetilde a_{i_0}$ then
$z$ is in $U$.
By the definition of $U_{a_{i_0}}$, it follows that
$U_{a_{i_0}} \subset U$.
Hence the collection $\{ U_{a_i}, i \in {\bf N} \}$ forms
a neighborhood basis at $p$.
\end{proof}

More importantly we have the following:

\begin{lemma}{}{}
The space $\cd = \oo \cup \partial \oo$ is second countable.
\label{sec}
\end{lemma}

\begin{proof}{}
We first construct a candidate for a countable basis.
Since $\oo$ is homeomorphic to ${\bf R}^2$ it has a countable
basis ${\cal B}_1$. 
Let \ ${\cal Z} = \{ l \ | \ l  \ {\rm is \ a \ periodic \ leaf \ of }
\oos \cup \oou \}$.
Let 

$${\cal B}_2 \ = \ \{ \ U_{b_i} \ | \ b_i \in B = (b_i), \
B  \ {\rm  admissible, 
\ where} \  b_i \ {\rm has \ 
all \ sides \ contained \ in \ leaves \ in} \  {\cal Z}\  \}$$

There are countably many leaves in ${\cal Z}$ and so countably
many intersections of these leaves. Since any polygonal path
is a union of a finite number of sides, it now follows that
${\cal B}_2$ is a countable collection of open sets in $\cd$.
We want to show that ${\cal B} = {\cal B}_1 \cup {\cal B}_2$ 
is a basis for the topology in $\cd$.

Let $p$ in $\cd$ and $V$ open set in $\cd$ containing $p$.
If $p$ is in $\oo$ there is $U$ in ${\cal B}_1$ with 
$p \in U \subset V$.
Suppose then that $p$ is in $\partial \oo$.
Choose $A = (a_i)$ a master sequence for $p$.
According to definition \ref{topology} there is $j$ with 
$U_{a_j} \subset V$.

We now modify the sides of the $a_j$ to a convex polygonal path with
sides in ${\cal Z}$. 
The sides of $a_j$ in periodic leaves are left unchanged.
A side in a non periodic leaf is pushed slightly in the
direction of $\widetilde a_j$ to a periodic leaf.
Notice that the union of periodic leaves of $\oos$ (or $\oou$)
is dense in  $\oo$.
The proof is done in 2 steps.
First we do this for the finite sides. The obstruction
to pushing in a side of $a_j$, still intersecting
the same adjacent sides is that there is a singularity
in this side. But then this segment is already in
a periodic leaf and we leave it unchanged.
Do this for all finite sides of $a_i$
to produce a new polygonal path  $b_i$.
Do this  for all $i$.
Given $i$, then 
since $a_j$ escapes in $\oo$ with increasing
$j$, then the finite segments of $a_j$ are eventually  contained
in $\widetilde b_i$.
Hence the finite segments of $b_j$ are contained
in $\widetilde b_i$. One can then take a subsequence of the $(b_i)$
so that $B = (b_i)$ is nested.
The $b_i$ are convex and also 
$(b_i)$ is eventually nested with the $(a_i)$. 
This implies that $B = (b_i)$ is also a master sequence for $p$.

The second step is to modify the rays of $B = (b_i)$ to be in
periodic leaves.
Given $i$, consider one 
ray $l$ of $b_i$ and \ $l_t,  \  t \geq 0$ \ 
leaves of the same foliation as $l$,
with $l_t$ converging to $l$
as $t \rightarrow 0$. 
In addition the 
$l_t$ intersect the side of $b_i$
adjacent to $l$.
Note that this intersection of $l$ and the adjacent side is not
a singular point, otherwise $l$ is periodic and we do not need to change it.
If the $l_t$ converges
to another leaf (in $\widetilde b_i$ or not) besides $l$,
then $l$ is in a non Hausdorff leaf and theorem \ref{theb}
implies
that $l$ is in a periodic leaf and again we leave $l$ as is.
So we may assume that 
as $t \rightarrow 0$
then $l_t$ converges only to the leaf of $\oos$ or $\oou$ containing $l$.
There is $j_i > i$ so that $l$ does not have a subray which is a side of
$b_{j_i}$ $-$ otherwise $B =  (b_j)$
would not be a master sequence, by lemma \ref{choi}.
Then there is $t$ sufficiently small so that $l_t$ separates
$l$ from $b_{j_i}$. This is true because $l_t$ does not 
converge to any other leaf besides $l$.
Choose also one $t$ for which $l_t$ is a periodic leaf
and replace the ray $l$ of $b_i$ by this ray of $l_t$.
After doing this to both rays of $b_i$ this produces a
convex polygonal path $(c_i)$. 
For each $i$ then \ $\widetilde b_{j_i}
\subset \widetilde c_i$, so $\widetilde c_{j_i} \subset
\widetilde c_i$. So after taking a subsequence $C = (c_n)$ is
nested. By the above, $C \cong B$ and $C$ is a master sequence
for $p$. 

Hence we can find $n$ with $U_{c_n} \subset V$.
But all the sides of $c_n$ are periodic. 
This shows that ${\cal B}$ is a basis for the topology of $\cd$ and
finishes the proof of the lemma.
%
\end{proof}

Next we show that $\cd$ is a regular space, which will imply that
$\cd$ is metrizable.

\begin{lemma}{}{}
The space $\cd$ is a regular space.
\label{regul}
\end{lemma}

\begin{proof}{}
Let $p$ be a point in $\cd$ and $V$ be a closed set not containing
$p$. Suppose first that $p$ is in $\oo$. 
Here $V^c$ is an open set with $p$ in $V^c$, so there are
open disks $D_1, D_2$  in $\oo$, so that 
$p \in D_1 \subset \overline D_1 \subset D_2 \subset V^c$,
producing disjoint neighborhoods $D_1$ of $p$ and $(D_2)^c$
of $V$.

Suppose now that $p$ is in $\partial \oo$. 
Since $p$ is not in the closed set $V$, there is an open
set $O$ of $\cd$ containing $p$ and disjoint from $V$. Let $A = (a_i)$ be
a master sequence associated to $p$. Then there is $i_0$
so that $U_{a_{i_0}}$ defined above is contained in $O$.
We claim that the closure of $\widetilde a_{i_0}$ 
in $\cd$ is
$U_{i_0}$ union 
$a_{i_0}$ plus the two ideal points
of the rays in $a_{i_0}$. Clearly the closure of
$\widetilde a_{i_0}$ in $\cd$ intersected with $\oo$ is
obtained by  just
adjoining $a_{i_0}$. An ideal points of a ray $l$ of $a_i$ is 
clearly in the closure as any neighborhood of it contains a subray
of $l$.
Any other point $p$ in $\partial \oo$, if it is in $U_{a_{i_0}}$, then
it is in the closure of $\widetilde a_{i_0}$.
If $p$ is not in $U_{i_0}$ and is not an ideal point of $a_{i_0}$ then
find a master sequence for $p$ disjoint
from master sequences of both ideal points of $a_{i_0}$ and
hence disjoint from $U_{i_0}$. Hence $p$ is not 
in the closure of $\widetilde a_{i_0}$.
This proves the claim.

Choose $j$ big enough so that 
the rays of $a_j$ are not equivalent to
any ray of $a_{i_0}$,
again possible by lemma \ref{choi}. 
By the above it follows that the closure of
$\widetilde a_j$ is contained in $U_{a_{i_0}}$, hence

$$p \in U_{a_j} \ \subset \ \ {\rm closure}(\widetilde a_j)
\ \subset \ U_{a_{i_0}} \ \subset \ O \  \subset \ V^c$$

\noindent
This proves that $\cd$ is regular.
\end{proof}

\begin{corollary}{}{}
The space $\cd$ is metrizable.
\end{corollary}

\begin{proof}{}
Since $\cd$ is second countable and regular, 
the Urysohn metrization theorem 
(see \cite{Mu}
pg. 215) implies that $\cd$ is
metrizable.
\end{proof}

Therefore in order to prove that $\cd$ is compact it suffices
to show that any sequence in $\cd$ has a convergent subsequence.
But it is quite tricky to get a handle on an arbitrary sequence
of points in $\oo$ or in $\partial \oo$ and 
the proof that $\cd$ is compact is hard.
This is the key property of $\cd$.
We first analyse one
case which seems very special, but which in fact implies the
general case without much additional work. Its proof is
very involved because there are many cases to consider.
First a preliminary result involving non separated leaves.
By theorem \ref{theb} this does not occur in the case
without perfect fits.

\begin{figure}
\centeredepsfbox{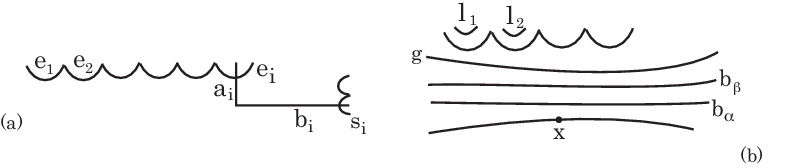}
\caption{
a. Infinitely many non separated leaves
converge to a single ideal point, b. A more interesting situation.}
\label{erg}
\end{figure}

\begin{lemma}{}{}
Let $\{ E_i \}_{i \in {\bf Z}}$ be leaves of (say) $\oos$ which 
are all non separated from each other and ordered as in
theorem 
\ref{theb}. Associated to this collection there are two
ideal points of $\oo$, one for $(E_i)$ with $i$ converging to 
infinity and another one for $(E_i)$ with $i$ converging to
minus infinity.
A master
sequence for any one of them 
is obtained with polygonal paths with length 2.
\label{infin}
\end{lemma}

\begin{proof}{}
As explained in the end of section \ref{prelim},
the collection $\{ E_i \}$ is part of the boundary of a
scalloped region ${\cal S}$.
We will follow the notation from that section.
The region ${\cal S}$ is the union
of infinitely many lozenges $A_i$ and parts
of their boundaries so that a half
leaf of $E_i$ is contained in the boundary of $A_{2i}$ and another 
half leaf of $E_i$ is contained in the boundary of 
$A_{2i-1}$. The lozenges $A_i$ and
$A_{i+1}$ are adjacent for any $i \in {\bf Z}$
and they all intersect
a single stable leaf $C$. 
This is depicted in fig. \ref{scal}.
Let $\gamma_i$ be the periodic orbits in $E_i$.
The collection of lozenges $\{ A_i \}$ also creates another bi-infinite
collection
$\{ G_i \}, i \in {\bf Z}$ of leaves of $\oos$, all of which are
non separated from each other and $G_i$ has a half leaf in the
boundary of $A_{2i}$ and another half leaf in the boundary of
$A_{2i+1}$. Let $\delta_i$ be the periodic orbit in $G_i$.
The boundary of ${\cal S}$ also has two bi-infinite collections 
of non separated leaves from $\oou$: \ $\{ S_j \}_{j \in {\bf Z}}$ and
$\{ T_j \}_{j \in {\bf Z}}$. These are chosen so that $\wu(\gamma_i)$
converges to $\{ S_j \}$ when $i \rightarrow \infty$ and
$\wu(\gamma_i)$ converges to $\{ T_j \}$ when $i \rightarrow -\infty$.
In addition $S_j$ has a periodic orbit $\tau_j$ and we choose the
indexing so that $\ws(\tau_j)$ converges to $\{ E_i \}$ when
$i \rightarrow \infty$ and $\ws(\tau_j)$ converges to $\{ G_i \}$
when $i \rightarrow -\infty$.
The collections $\{ G_i \}, \{ S_j \}$ are ordered with increasing
$i, j$, see also theorem
\ref{theb}.

Now we define the ideal point 
associated to $\{ E_i \}_{i \in {\bf Z}}$ when $i$ converges
to $\infty$.
For each positive $i$ choose rays $a_i$ in $\oou$, $b_i$ in
$\oos$ which intersect only in their starting point $u_i$
which is a point in ${\cal S}$
and $a_i$ intersects $E_i$ and $b_i$ intersects $S_i$,
see fig. \ref{erg}, a.
Let $d_i = a_i \cup b_i$,  let $\widetilde d_i$ be the component
of $\oo - d_i$ which contains $E_k$ for $k > i$ and
$S_k$ for $k > i$. 
The $d_i$ are polygonal paths of length 2.
It follows that $d_i$ is convex for $\widetilde d_i$.
This uses the particular ordering in $\{ E_i \}, \{ S_j \}$ described
above and it also follows that $(d_i)$ is a nested sequence of
polygonal paths.

In the explicit model ($V$) for a scalloped region given in
the end of section \ref{prelim} we can choose 

$$u_i = (1 - \frac{1}{2i-1}, 1 - \frac{1}{2i-1}), \ \ \
a_i \cap \overline {\cal S} \ = \ \{ 1 - \frac{1}{2i-1} \} \times [1 - \frac{1}{2i-1}, 1], \ \ \
b_i \cap \overline {\cal S} \ = \ [1 - \frac{1}{2i-1}, 1] \times \{ 1 - \frac{1}{2i-1} \}.$$

\noindent
Notice that the ray $a_i$ of $\oou$ is clearly not contained  in $\overline {\cal S}$,
only the part contained in $\overline {\cal S}$ has a description in the explicit model.
Similarly the ray $b_i$ of $\oos$ is not contained in $\overline {\cal S}$.
It remains to check that the sequence $(d_i)$ escapes compact sets in $\oo$
as $i \rightarrow \infty$.
In the explicit model
the $a_i$ are subsets of the leaves $\wu(\gamma_i)$. Any point in the limit of 
the sequence $(\wu(\gamma_i))$ is non separated from the
 $\{ S_j \}_{j \in {\bf Z}}$ and hence
has to be in one of the $S_j$. It follows that the part of $a_i$ outside
$\overline {\cal S}$ escapes compact sets in $\oo$. By contruction
the sequence made up of the parts of 
$a_i$ in $\overline {\cal S}$ also  does not limit in $\oo$, hence $(a_i)$ 
escapes compact sets in $\oo$. The same is true for $(b_i)$ so $(d_i)$ 
escapes compact sets in $\oo$ and  so $D = (d_i)$ is admissible and
defines an ideal point $p$ of $\oo$.
This $p$ is associated to the positive 
infinite direction
of the $\{ E_i \}_{i \in {\bf Z}}$. 
By lemma \ref{choi}, $D$ is a master sequence.
Similarly associated to the negative  direction of the $\{ E_i \}$ there is
another ideal point $q$ of $\oo$. 
\end{proof}

An ideal  point $p$ associated to infinitely many non separated leaves
or equivalently to a scalloped region is called
a {\em {corner}} of the scalloped region.

The technical lemma in the special case is the following:

\begin{lemma}{}{} 
Let \ $(l_i), i \in {\bf N}$
be a sequence of line leaves of $\oos$ (or $\oou$)
and let $z_i$ in $l_i$.
Suppose that for each $i$ the set $\oo - l_i$ has a component
$C_i$ so that each $C_i \cup l_i$ contains $\oos(z_i)$ and
also that the collection $\{ C_i \cup l_i \}$ is pairwise disjoint. 
Suppose that the ordering of $l_i$ (see definition 
\ref{order1}) is chosen so that the $l_i$ are linearly
ordered with $i$.
Then in $\cd$, the sequence  $(C_i \cup l_i)$ converges to a point
$p$ in $\partial \oo$.
\label{forcing}
\end{lemma}

\begin{proof}{}
The proof of this lemma is very involved because there are many
possibilities and many places where the leaves $l_i$ can slip through.

Suppose that $l_i$ is always in $\oos$ as other cases
are similar.
If the $l_i$ does not escape compact sets
in $\oo$ when $i \rightarrow \infty$ 
then there are $i_k$ and $z_{i_k}$
in $l_{i_k}$ with $z_{i_k}$ converging to a point $z$.
But then
the $C_{i_k}$ cannot all be disjoint, contradiction.
Hence the $(l_i)$ escapes in $\oo$.

First notice that because the collection $\{ l_i \}$
is linearly ordered with $i$, then if a subsequence
$(l_{i_k} \cup C_{i_k})$ converges to $p$ in $\cd$, then
the full sequence $(l_i \cup C_i)$ also converges to
$p$ in $\cd$.
Choose $z_i$ in $l_i$.

\vskip .1in
\noindent
{\underline {Case 1}} $-$ 
There is an infinite subsequence of the
$(l_j)$,
which we may assume is the original sequence so that $l_j$
are all non separated from $l_1$
(in particular there are perfect fits).

Then the  $\{ l_j, j \in {\bf N} \}$ forms a subcollection of
a collection $\{ z_i \}_{i \in {\bf Z}}$ of non separated leaves
of $\oos$ as in lemma \ref{infin}.
Hence we can find $a_j, b_j$ as in the
previous lemma and for any $i$, 
$l_i$ intersects $a_{j_i}$ where $j_i$ goes to 
infinity with $i$.
As in the lemma let $d_j = a_j \cup b_j$ and
$D = (d_j)$.  Then $D$ is a master sequence 
defining a point $p$ in $\partial \oo$.
In addition given any $j$ then 
for $i$ big enough $l_i$ is contained
in $\widetilde d_j$.
Hence $l_i \cup C_i$ converges to $p$ in $\cd$.

\vskip .1in
\noindent
{\underline {Case 2}} $-$
Up to subsequence, 
for any distinct $i, j$, 
the $l_i$ is separated from $l_j$. 

Let $V = \oo - \cup_{i \in {\bf N}} (C_i \cup l_i)$, an open set in $\oo$.
The procedure will be to inductively construct leaves
$g_n$ 
so that either the sequence $(g_n), n \in {\bf N}$ is
nested with $n$ and escapes compact sets in $\oo$
or is a sequence of non separated leaves.
There are various possibilities for the limiting behavior
of $(g_j)$ which will eventually lead to a proof
that $(l_i \cup C_i)$ converges in $\cd$.

Given $x$ in $\oo$
consider the line leaves $b$
of $\oos$ which separate $x$ from ALL of
the $l_i$. 
For example given $y$ not in the union of $l_i \cup C_i$, then
$\oos(y)$ is disjoint from this union $-$ this is because
no prong of $\oos(z_i)$ is contained in $V$.
For any $x$ in a complementary
region of $\oos(y)$ not interecting this union will have such
line leaves $b$.
A singular leaf has at most two line leaves with this property.
The collection
of line leaves $b$ as above is clearly ordered by separation properties so
we can index then as $\{ b_{\alpha} \ | \ \alpha \in J \}$
where $J$ is an index set.
Put an order in $J$ so that $\alpha < \beta$ if and only if
$b_{\alpha}$ 
separates some point in $b_{\beta}$ from $x$.
Equivalently $b_{\beta}$ separates some point in $b_{\alpha}$ from $x$.
Two such line leaves in the same stable leaf may share the singular point
or a half leaf.
Since the $b_{\alpha}$ cannot escape $\oo$ as $\alpha$ increases
(they are bounded by all the $l_i$) then the
$\{ b_{\alpha} \}$ limits to a collection
of leaves 
of $\oos$ as $\alpha$ grows without bound.

There are 2 options:
1) There are infinitely many line leaves $s_n$ of $\oos$ in the limit of
the $b_{\alpha}$ so that for
each $n$ there is $i_n$ with $s_n$ either equal to $l_{i_n}$ or
separating $l_{i_n}$ from
every $b_{\alpha}$, \
2) There is one line leaf $s$ of $\oos$ in the limit of the $b_{\alpha}$ so that
this single 
$s$ separates infinitely many of the $l_i$ from all of the
$g_{\alpha}$.
Notice that only option 2) can happen when there are no
perfect fits.

\vskip .05in
Consider first option 1. The collection of leaves 
non separated from the $s_n$ is infinite. Because the
$l_i$ are ordered it now follows that each $s_n$ can separate
only finitely many of the $l_i$ from all of the $b_{\alpha}$.
Let $p$ be the ideal point given by lemma \ref{infin}
associated to the direction of the $s_n$ with $n$ increasing.
The proof of lemma \ref{infin} implies that $(l_i \cup C_i)$
converges to $p$.

\vskip .05in
Now consider option 2. Let $g_0 = s$.
The leaf $t$ of $\oos$ containing 
$s$ may have singularities. By the condition of 
pairwise disjointness
of the $l_i \cup C_i$,  there is a single line
leaf $g_1$ of $t$ with a complementary component
$o_1$ in $\oo$ which contains $l_i$ for all $i \geq i_0$.
We will restart the process with the 
$\{ l_i \}, i \geq i_0$ instead of the original sequence. 
We will remember $g_0$ and the leaf $g_1$ which separates $x$ from
all $l_i \cup C_i$, $i \geq i_0$.

%

Restart the process as follows.
Throw out all the leaves until $l_{i_0}$ and redo the process.
This iterative
process produces $(g_j), j \in {\bf N}$ 
which is a 
weakly nested sequence of line leaves.
We explain the weak behavior. For instance in the
first case, after throwing out $l_1$ (or whatever first leaf
was still present), it may be that only $g_1$ is a slice
which separates $x$ from all other $l_i$, see fig.
\ref{onecon}, b. In that
case $g_2 = g_1$. So the $g_j$ may be equal, but
they are weakly monotone with $j$.

If the $(g_j)$ escapes in $\oo$
with $j$ then it defines a point $p$ in $\partial \oo$.
Since each $g_j$
separates infinitely many $l_i$ from $x$ we quickly obtain
as before
that the $l_i$ converge to the point $p$ in $\partial \oo$.

Suppose then that the $(g_j)$ does
not escape in $\oo$. 
The first option is that 
there are infinitely many distinct $g_j$. 
Up to taking a subsequence assume all $g_j$ are distinct
and let $g_j$ converge to $ H = \cup h_k$, a collection
of line leaves of $\oos$. 
By construction, for each $j_0$, the 
$g_{j_0}$ separates some
$l_i$ from $x$ but for a bigger $j$, the $g_j$ does not
separate $l_i$ from $x$, see fig. \ref{onecon}, a.
Also, for each $i$ there is some $j$ so that $g_j$ 
separates $l_i$ from $x$.
In particular there is a component of $\oo - H$ which contains
all the $l_i$.

\begin{figure}
\centeredepsfbox{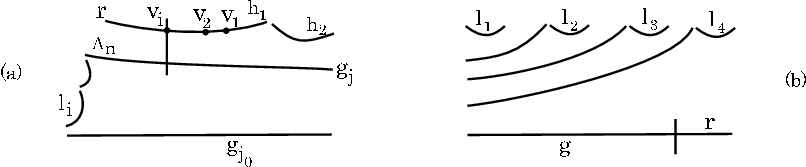}
\caption{
a. Forcing convergence on one side, b. The case that all $g_j$ are
equal.}
\label{onecon}
\end{figure}

We analyse the case there are finitely many line leaves
of $\oos$ in $H$, the other case being
similar.
As seen in theorem \ref{theb} the set of leaves in $H$
is ordered and we choose $h_1$ to be the leaf closest
to the $l_i$.
Also there is  a ray $r$ of $l$ which points in the
direction of the $l_i$, see fig. \ref{onecon}, a.
Let $p$ be the ideal point of $r$ in $\partial \oo$.
We want to show that $l_i \cup C_i$ converges to $p$.

Choose points $v_n$ in $r$ converging to $p$.
For each $n$ then $\oou(v_n)$ intersects $g_j$ for $j$
big enough $-$ since the sequence $g_j$ converges to $H$.
Choose one such $g_{j(n)}$ with $j(n)$ converging to infinity
with $n$.
We consider a convex set $A_n$ of $\oo$ bounded by
a subray of $r$ starting at $v_n$, a segment in 
$\oou(v_n)$ between $h_1$ and $g_{j(n)}$ 
and a ray in $g_{j(n)}$ starting in $g_{j(n)} \cap \oou(v_n)$
and going in the direction of the $l_i$, see fig. \ref{onecon}, a.
We can choose $j(n)$ so that the $(A_n), n \in {\bf N}$
forms a nested sequence.
Let $a_n = \partial A_n$.
Since $h_1$ is the first element of $H$ it follows that 
$(a_n)$ escapes compact sets in $\oo$ and clearly
it converges to $p$ in $\oo \cup \partial \oo$. 
For each $n$ and associated $j$, there
is $i_0$ so that for $i > i_0$ then $g_j$ separates
$l_i$ from $x$. If follows that $l_i \cup C_i$ is contained
in $A_n$ and therefore $(l_i \cup C_i)$ converges to $p$ 
in $\cd$. This finishes the proof in this case.

If $H$ is infinite let $H = \{ h_k, k \in {\bf Z} \}$ with
$k$ increasing as $h_k$ moves in the  direction of the 
$l_i$.
Then $h_i$ converges to a point $p \in \oo$.
A similar analysis as in the case that $H$ is finite
shows that $(l_i \cup C_i)$ converges to $p$ in $\cd$.
Use the convex chains $a_j \cup b_j$ as described in
lemma \ref{infin}.

\vskip .1in
The final case to be considered is that up to subsequence
all $g_i$ are equal
and let $g$ be this leaf. In particular no $l_i$ is equal
to $g$. This can certainly occur as shown in fig. \ref{onecon}, b. 
If we remove finitely many of the $l_i$, then $g$ is still
the farthest leaf separating $x$ from all the remaining $l_i$.
Notice also the $g$ is a line leaf on  the side containing
all the $l_i$.

Consider the collection of leaves ${\cal B}$ of $\oos$ 
non separated from $g$ in the side of $g$ containing the $l_i$.
Let $W$ be the component of $\oo - {\cal B}$ 
which accumulates on all of ${\cal B}$ if ${\cal B} \not = \{ g \}$ and
otherwise
let $W$ be the component of $\oo - \{ g \}$ not containing $x$.

One possibility is that there are infinitely many $i$ so that
$l_i$  is separated from $g$ by an element in ${\cal B}$.
Here we have 2 options. The first option is that there
are infinitely many distinct elements $e$ in ${\cal B}$ for which
there is some $l_i$ with $e$ separating $l_i$ from $g$,
see fig. \ref{respt}, a.
Since the $l_i$ are nested then as seen before this implies
that the $l_i \cup C_i$ converge to some $p$ in $\partial \oo$.
The second option here is that there is some fixed $h'$ in
${\cal B}$ which separates infinitely many $l_i$ from
$g$.
As the sequence $(l_i)$ is nested, this is true for all $i \geq i_0$ for
some $i_0$.
But then $h'$ would eventually  take the place of $g$ 
in the iterative process $-$ that is, some $g_k = h'$ instead of
$g_k = g$. Then $g_k$ is not eventually constant and this was dealt
with previously.

\vskip .1in
The remaining case is that after throwing out a few initial
terms we may assume that
all $l_i$ are contained $W$, see fig. \ref{respt}, b.
Fix an embedded arc $\gamma$ from $g$ to $l_1$ intersecting
them only in boundary points and not intersecting
any other $l_i$.
Let $T$ be the component of $\oo - (g \cup \gamma \cup l_1)$
containing all other $l_i$.
Put an order in ${\cal B}$ so that elements of ${\cal B}$ contained
in $T$ are bigger than $g$ in this order.
For simplicity assume that ${\cal B}$ is finite.
The case where there are infinitely many leaves
non separated from $g$ on that side is very similar with
proof left to the reader.
Let $h$ be the biggest element of ${\cal B}$, which could be $g$ itself.
Let $r$ be the ray of $h$
associated to the increasing direction of the the
$l_i$ and let $p$ in $\partial \oo$ be the ideal point of $r$.
We want to show that $(l_i \cup C_i)$ converges to $p$.

Let $A$ be an arbitrary convex neighborhood of $p$ in $\cd$ bounded
by a convex chain $a$, see fig. \ref{respt}, b.
If $A$ is small enough then $a$ has a ray $r_1$ contained
in $T$.
The rays $r, r_1$ are not equivalent.
Let $h'$ be a  leaf of $\oos$ in $W$ sufficiently close to $g$.
Becauce $h$ is the biggest element in the ordered set ${\cal B}$
then $h'$ has to have a ray contained in $A$.
For $h'$ close enough to $g$, since the $l_i$ are in $T$, then
for some $i_0$ the leaf 
$h'$ separates $l_i, 1 \leq i \leq i_0$ from $g$ and hence from
$x$. By the maximality property of $g$, then for some $j$
the leaf $h'$ does not separate $l_j$ from $g$. Since
$l_j$ is in $T$ this forces $l_j$ to be contained in $A$.
As the $\{ l_i, i \in {\bf N} \}$ forms an ordered collection
this forces $l_i$ to be contained in $A$ for all $i \geq j$.
Since $A$ was an arbitrary neighborhood of $p$ this shows
that $(l_i \cup C_i)$ converges to $p$ in $\cd$.


\begin{figure}
\centeredepsfbox{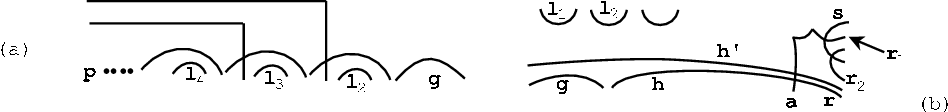}
\caption{
a. The $l_i$ flip to the other side of a
a leaf non separated from $g$, \  b. Convex neighborhood disjoint
from all.}
\label{respt}
\end{figure}

This finishes the proof of lemma \ref{forcing}.
\end{proof}


\begin{proposition}{}{}
The space $\cd$ is compact.
\label{comp}
\end{proposition}

\begin{proof}{}
Since $\cd$ is metrizable,
it suffices to consider the behavior of an arbitrary 
sequence $z_i$ in $\cd$.
We analyse all possibilities and in each case show there is
a convergent subsequence.

Up to taking subsequences there are 2 cases:

\vskip .1in
\noindent
{\underline {Case 1}} $-$ Assume the $z_i$ are all in $\oo$. 

If there is
a subsequence of $z_i$ in a compact set of $\oo$, then
there is a convergent subsequence as $\oo \cong \rrrr^2$.
So assume from now on that $z_i$ escapes compact
sets in $\oo$.
Let $b_i = \oos(z_i)$.
Suppose first there is a subsequence $(b_{i_k})$ converging
to $b$ and assume that all $b_{i_k}$ are in one
sector of $b$ or in $b$ itself.
If a subsequence of $(b_{i_k})$ is constant and hence equal to $b$ then
up to another subsequence the $z_i$ converges in $\cd$ to one of
the ideal points of $b$, done.
Otherwise a small transversal to $b$ in a regular unstable leaf
intersects $b_{i_k}$ for 
$k$ big enough and up to subsequence assume all $z_{i_k}$ are in
one side of that unstable leaf.
Suppose for simplicity there
are only finitely many leaves non separated from $b$ 
in the side containing the $b_i$. 
Let $b'$ be the last one non separated from $b$ in the side the 
$b_{i_k}$ are in and let $p$ be the ideal
point of $b'$ in that direction. 
The argument is similar to one in case 2 of
lemma \ref{forcing}:
let $v_n$ in $b'$ converging to $p$ in $\partial \oo$
with $\oou(v_n)$ regular.
Choose a convex polygonal path $a_n$ made up of the ray in
$b'$ starting in $v_n$ and converging to $p$, then
the segment in $\oou(v_n)$ from $v_n$ to $\oou(v_n) \cap
b_{i_k}$ for apropriately big $k$ and then a ray 
in $b_{i_k}$ starting in this point. As before we can
choose the $\widetilde a_n$ nested with $n$
and so that $(\widetilde a_n \cup a_n)$ escapes compact sets in $\oo$,
so converges
to $p$ in $\cd$.
It follows that 
$z_{i_k}$ converges to $p$ and we are done in this case.
The case of infinitely many leaves non separated from $l$
is treated similarly to what is done in the proof of lemma \ref{forcing}.

Suppose now that the sequence $(b_i), i \in {\bf N}$ escapes compact sets
in $\oo$. 
The goal is to reduce this case to a situation where we can apply lemma
\ref{forcing}.
Fix a base point $x$ in $\oo$ and assume that
$x$ is not in any $b_i$.
Let $l_i$ be the line leaf of $b_i$ (so $l_i$ is a
line leaf of $\oos$) which is the boundary of the component
of $\oo - b_i$ containing $x$. Let $C_i$ be the component of 
of $\oo - l_i$ not containing $x$. If $b_i$ is regular then $C_i$ is a
component of $\oo - b_i$. If $b_i$ is singular then $C_i \cup l_i$ 
contains all the prongs of $b_i$.
In this case it follows that 
$C_i$ escapes compact sets in $\oo$.
If there is a subsequence $(l_{i_k})$ so that 
$(l_{i_k})$ is nested then this defines an admissible sequence
of convex polygonal paths (of length one) converging to an ideal
point $p$.

Otherwise there has to be $i_1$ so that there are only finitely many
$i$ with $C_i \subset C_{i_1}$. Choose
$i_2 > i_1$ with $C_{i_2} \not \subset C_{i_1}$
and hence 
$C_{i_2} \cap C_{i_1} = \emptyset$ and
also so that there are finitely many $i$ with
$C_i \subset C_{i_2}$.
In this way we construct a subsequence $i_k, k \in {\bf N}$
with $C_{i_k}$ disjoint from each other. 
The collection of line leaves

$$\{ \ l_{i_k} \ \ | \ \  k \in {\bf N} \}$$

\noindent
is circularly ordered 
and if we remove one element of the sequence 
(say the first one) then it is linearly ordered. As such it can be
mapped injectively into the set of rational numbers
${\bf Q}$ in an order preserving way. 
Therefore there is another subsequence (call it still
$(l_{i_k})$) for which the set $\{ l_{i_k} \}$ is now linearly
ordered with $k$
$-$ either increasing or decreasing.
We can now apply lemma \ref{forcing} to the sequence $l_{i_k}$
and obtain that
$(l_{i_k})$ converges to a point $p$ in $\partial \oo$ and
hence so does $z_{i_k}$.
It was crucial here that $C_i \cup l_i$ contains all the prongs of
$b_i$ in order to apply lemma \ref{forcing}.

This finishes the analysis of case 1.

\vskip .1in
\noindent
{\underline {Case 2}} $-$ Suppose the $z_i$ are in $\partial \oo$.

We use the analysis of case 1.
We may assume that the points $z_i$ are pairwise distinct.
To start 
we can find a convex polygonal path 
$a_1$ so that 
$\overline U_{a_1}$ contains a neighborhood of $z_1$ in $\cd$ and
also it does not contain any other $z_i$.
Otherwise there is a subsequence of
 $(z_i)$ which converges to $z_1$.
Inductively construct $a_i$ convex polygonal paths with
$\overline U_{a_i}$ neighborhood of $z_i$ in $\cd$ 
and the $\{ \overline U_{a_j} \}, 1 \leq j \leq i$ pairwise disjoint.
By taking smaller convex neighborhoods we can assume that
the $(U_{a_i})$ escapes compact sets in $\oo$ as $i \rightarrow \infty$.
As in case 1 we may assume up to subsequence
that the $\{ a_i \ | \  i \in {\bf N} \}$
forms an ordered set of $\oo$ with the order given by $i$.
Let $w_i$ be a point in $a_i$. Since $a_i$ escapes compact sets in $\oo$,
case 1 implies that there is a subsequence $w_{i_k}$ converging
to a point $p$ in $\partial \oo$.
Consider a master sequence $B = (b_j)$ associated to $p$.
Let $j$ be an integer. If for all $k$ we have that 
$\widetilde a_{i_k} \not \subset 
\widetilde b_j$, then $\widetilde a_{i_k}$ has a point $w_{i_k}$ converging
to $p$ and also has points outside $\widetilde b_j$. This contradicts
the $\widetilde a_{i_k}$ being all disjoint since they are
convex. Therefore $\widetilde a_{i_k} 
\subset \widetilde b_j$ for $k$ big enough $-$ this follows 
because the sequence $(a_{i_k})$ is ordered as a subset of $\oo$.
In fact by increasing the index if necessary then 
\ $U_{a_{i_k}}  \subset {\rm closure} (\widetilde b_j)$
in $\cd$. Since $z_{i_k}$ is in $U_{a_{i_k}}$
this shows that $z_{i_k} \rightarrow p$.
Therefore there is always 
a subsequence of the original sequence 
which converges to a point in $\partial \oo$.

This finishes the proof of proposition \ref{comp}, compactness of $\cd$.
\end{proof}

We now prove a couple of additional properties of $\cd$.

\begin{proposition}{}{}
The space $\partial \oo$ is homeomorphic to a circle.
\label{circle}
\end{proposition}

\begin{proof}{}
The space $\partial \oo$ is metrizable and circularly
ordered. 
Also $\partial \oo$ is
compact, being a closed subset of
a compact space $-$ since $\oo$ is open in $\cd$.
We now show that $\partial \oo$ is connected, no points
disconnect the space and any two points disconnect the
space.

Let $p, q$ be distinct points in $\partial \oo$. Choose disjoint
convex neighborhoods \ $\overline {U_a}, \ \overline {U_b}$ 
\ of $p, q$ defined
by convex polygonal paths $a, b$. There are ideal points of $\oo$ in 
$\overline {U_a}$ distinct from $p$, hence there is a point in $\partial
\oo$ between $p, q$. Hence any ``interval" in $\oo$ is a linear
continuun, being compact and satisfying the property that between
any two points there is another point.
This shows that $\partial \oo$ is connected and also
that no point in $\partial \oo$ disconnects it.
In addition as $\partial \oo$ is circularly ordered, then any two 
points disconnect $\partial \oo$. By theorem I.11.21, page 32 of Wilder
\cite{Wi}, 
the space $\partial \oo$ is homeomorphic to a circle.
\end{proof}

We are now ready to prove that $\cd$ is homeomorphic to a disk.

\begin{theorem}{}{}
The space $\cd = \oo \cup \partial \oo$ is homeomorphic to
the closed disk $D^2$.
\label{disk}
\end{theorem}

\begin{proof}{}
The proof will use classical results of general topology,
namely a theorem of Zippin characterizing the closed
disk ${\bf D}^2$, see theorem III.5.1, page 92 of Wilder \cite{Wi}.

First we need to show that $\cd$ is a Peano continuun, see
page 76 of Wilder \cite{Wi}.
A Hausdorff topological space $C$ is a {\em Peano space}
if it is not a single point, it is second countable,
normal, locally compact, connected and locally connected.
Notice that Wilder uses the term perfectly
separable (definition in page 70 of \cite{Wi}) 
instead of second countable.
If in addition $C$ is compact then $C$ is
a {\em Peano continuun}.

By proposition
\ref{comp} our space $\cd$ is  compact, hence
locally compact.
It is also Hausdorff  $-$ lemma \ref{haus} $-$ hence normal.
By lemma \ref{sec} it is second countable and it is clearly
not a single point.
What is left to show is that $\cd$ is connected and
locally connected.

We first show that $\cd$ is connected. Suppose not and let
$A, B$ be a separation of $\cd$.
Since $\partial \oo$ is connected (this is done in
the proof of proposition \ref{circle}),
then $\partial \oo$ is contained in either $A$ or $B$,
say it is contained in $A$.
Then $B$ is contained in $\oo$.
If $B \not = \oo$, then $B, A \cap \oo$ disconnect $\oo$,
contrary to $\oo \sim \rrrr^2$. If $B = \oo$, then
$A = \partial \oo$ and so $\oo$ is closed in $\oo \cup \partial \oo$,
which is not true. It follows that $\cd$ is connected.


Next we show that $\cd$ is locally connected.
Since $\oo \cong \rrrr^2$, then $\cd$ is locally connected
at every point of $\oo$. Let $p$ in $\partial \oo$ and
let $W$ be a neighborhood of $p$ in $\cd$.
If $A = (a_i)$ is a master sequence associated
to $p$, there is $i$ with 
$\overline {U_{a_i}}$ contained in $W$ and 
$U_{a_i}$ is a neighborhood of $p$ in $\cd$.
Now $U_{a_i} \cap \oo = \widetilde a_i$ is homeomorphic
to $\rrrr^2$ also and hence connected.
The closure of $\widetilde a_i$ in $\cd$ is
$\overline {U_{a_i}}$. Since 

$$\widetilde a_i \ \ \subset \ \ U_{a_i} \ \ \subset \ \ 
\overline {U_{a_i}}$$

\noindent
then $U_{a_i}$ is connected.
This shows that $\cd$ is locally connected and hence that
$\cd$ is a Peano continuun.

To use theorem III.5.1 of \cite{Wi} we need the
idea of spanning arcs.
An arc in a topological space $X$ is a subspace homeomorphic
to a closed interval in $\rrrr$.
Let  $ab$ denote an arc with endpoints $a, b$.
If $K$ is a point set, we say that $ab$ {\em spans} $K$
if $K \cap ab = \{ a, b \}$.
We now state theorem III.5.1 of \cite{Wi}.

\begin{theorem}{(Zippin)}{}
A Peano continuun $C$ containing a $1$-sphere $J$ and satisfying
the following conditions below is a closed $2$-disk with
boundary $J$:

(i) $C$ contains an arc that spans $J$,

(ii) Every arc that spans $J$ separates $C$,

(iii) No closed  proper subset of an arc spanning $J$ 
separates $C$.
\end{theorem}

Here $E$ separates  $C$ mean that $C - E$ is not connected.

In our case $J$ is $\partial \oo$.
For condition (i) let $l$ be a nonsingular leaf in
$\oos$ or $\oou$.
Then $l$ has 2 ideal points in $\partial \oo$ which are distinct.
The closure $\overline l$ is an arc that spans $\partial \oo$.
This proves (i).

We prove (ii). Let $\zeta$ be an arc in $\cd$ spanning $\partial \oo$.
Then $\zeta \cap \oo$ is a properly embedded copy of $\rrrr$ in
$\oo$. Hence $\oo - (\zeta \cap \oo)$ has exactly two components
$A_1, B_1$. In addition $\partial \oo - (\zeta \cap \partial \oo)$
has exactly two components $A_2, B_2$ and they are 
connected, since 
$\partial \oo$ is homeomorphic to a circle by proposition
\ref{circle}.
If $p$ is in $A_2$  and $A = (a_i)$ is a master sequence for
$p$, then  by definition of the topology in $\cd$ there
is $i$ so that 
that $U = U_{a_i}$ is disjoint from $\zeta$ as $\zeta$ is closed
in $\cd$ and $p \not \in \zeta$.
Then $U \cap \oo  = U_{a_i} \cap \oo = \widetilde a_i$ is connected.
Hence $U \cap \oo$ is contained in either $A_1$ or $B_1$.
This also shows that a small neighborhood of $p$ in $\partial \oo$
will be contained in either $A_2$ or $B_2$.
By connectedness of $A_2, B_2$, then after
switching $A_1$ with $B_1$ if necessary  it follows
that: for any $p \in A_2$ there is a neighborhood
$U$ of $p$ 
in $\cd$ with $U \cap \zeta = \emptyset$ and
$U \cap \oo \subset A_1$.
Similarly $B_2$ is paired with $B_1$.
Let $A = A_1 \cup A_2$ and $B = B_1 \cup B_2$.
The arguments above show that $A, B$ are open
in $\cd$ and therefore they form a separation 
of $\cd - \zeta$.
This proves (ii).

Since $\oo - (\zeta \cap \oo)$  has exactly
two components $A_1, B_1$ then $\zeta \cap \oo$ is
contained in $\overline A_1 \cap \overline B_1$ 
and so $\zeta \subset \overline A \cap \overline B$.
It follows that no proper subset of $\zeta$ separates
$\cd$. This proves property (iii).

Now Zippin's theorem implies that $\cd$ is
homeomorphic to a closed disk. This finishes the
proof of theorem \ref{disk}.
\end{proof}

Notice that $\pi_1(M)$ acts on $\oo$ by homeomorphisms.
The action preserves the foliations $\oos, \oou$ and also 
preserves convex polygonal paths, admissible sequences, master sequences
and so on. Hence $\pi_1(M)$ also acts by homeomorphisms
on $\cd$. 
The next result will be very useful in the following section.

%
%

\begin{proposition}{}{}
Let $\Phi$ be a pseudo-Anosov flow in
$M^3$ closed. Let $p$ be an ideal point of $\oo$.
Then one of the 3 mutually exclusive options occurs:

1) There is a master sequence $L = (l_i)$ for $p$ where
$l_i$ are slices in leaves of $\oos$ or $\oou$.

2) $p$ is an ideal point of a ray $l$ of $\oos$ or
$\oou$ so that $l$ makes a perfect fit with another
ray of $\oos$ or $\oou$. There are master sequences
which are standard sequences associated to the ray $l$ 
in $\oos$ or $\oou$ as described in definition \ref{stan}.

3) $p$ is a corner of a scalloped region
as described in section \ref{prelim}.
Then a master sequence for $p$ is obtained as
described in lemma \ref{infin}.

In addition the only conclusion that applies if there are
no perfect fits is conclusion 1).
\label{options}
\end{proposition}

\begin{proof}{}
The point $p \in \oo$ is fixed in this proof.
We first show that cases 1) - 3) are mutually exclusive.
Case 2) it is disjoint from case 1).
This is because any master sequence $E = (e_i)$ in case 2) has
to have $\widetilde e_i$ containing part of a fixed perfect fit
for $i$ big enough.
In particular the polygonal paths $e_i$ have to have
at least 2 sides for $i$ big enough, so this cannot be case 1).
Suppose now that $p$ is a point of type 3). Consider a master
sequence $D = (d_i)$ where $d_i = a_i \cup b_i$, $a_i$ a
ray in $\oou$ and $b_i$ a ray in $\oos$ as described
in lemma \ref{infin}. Notice that all $a_i$ intersect a common
unstable leaf.
If there is a master sequence $L = (l_j)$ as in 1) then
the $l_j$ have to weakly intercalate with the $d_i$. But 
then they have to separate leaves of $\oou$ intersecting
a common leaf of $\oos$ and vice versa. This is impossible.
The same argument can be used to rule out case 2): consider
a master sequence $E = (e_j)$  as in case 2). The weak intercalation 
property of $d_i$ with this sequence implies that the
polygonal paths $e_j$ have to be eventually of length 2 and
both leaves have to be leaves intersecting the scalloped 
region. Hence $E$ is an admissible sequence as in Case 3) and
does not converge to an ideal point of a ray $l$ associated to
a perfect fit.

\vskip .1in
Now we prove that one of options 1) - 3) has to occur.
Fix a basepoint $x$ in $\oo$. Let $A = (a_i)$ be a master sequence
defining $p$. Since $(a_i \cup \widetilde a_i)$ escapes
compact sets in $\oo$, we may throw out a few initial terms
if necessary and assume that $x$ is not in the closure of any
$\widetilde a_i$.
Each $a_i$ is a convex polygonal path,
$a_i = b_1 \cup ... b_n$ where $b_j$ is
either a segment or a ray in $\oos$ or $\oou$.
For simplicity we omit the dependence of the $b_j$'s on the index $i$.

\vskip .1in
\noindent
{\underline {Claim}} $-$ For each $i$ there is some $b_j$ as above, with 
$b_j$ contained in a slice $z$ of a leaf of $\oos$ or $\oou$, so that
$z$ separates $x$ from $\widetilde a_i$.

In this claim $i$ is fixed.
Given $j$ let $y$ be an endpoint of $b_j$. Wlog assume
$b_j$ is in a leaf of $\oos$ and $y$ is in $b_{j+1}$ also.
Since $a_i$ is a convex polygonal path, we can extend
$b_j$ along
$\oos(y)$ beyond $y$ and entirely outside $\widetilde a_i$.
The hypothesis that $\widetilde a_i$ is convex is necessary,
for otherwise at a non convex switch any continuation 
of $b_j$ along $\oos(y)$ would have to enter
$\widetilde a_i$.
If one encounters a singular point in $\oos(y)$ (which could be
$y$ itself), then continue along the prong closest to $b_{j+1}$.
This produces a slice $c_j$ of $\oos(y)$ with
$b_j \subset c_j$. 
There is a component $V_j$ of
$\oo - c_j$ containing $\widetilde a_i$. Since we choose the prong
closest to $b_{j+1}$ then

$$\bigcap_{j = 1}^n V_j \ \  = \ \ \widetilde a_i$$

\noindent
Since $x$ is not in $\widetilde a_i$, then there is at least one $j$ with
$x$ not in $V_j$ and so $c_j$ separates $x$ from $\widetilde a_i$.
Let $z$ be this slice $c_j$.
This proves the claim.

\vskip .1in
Using the claim then
for each $i$ produce such a slice and denote it by $l_i$.
Let $\widetilde l_i$ be the component of $\oo - l_i$ containing
$\widetilde a_i$. Up to subsequence assume all the $l_i$ are
in (say) $\oos$. Since $A$ is a master sequence for $p$, we may
also assume, by lemma
\ref{choi},  that all the $l_i$ are disjoint from each other.

We now analyse what happens to the $l_i$. The first
possibility is that the sequence $(l_i)$ escapes compact
sets in $\oo$. 
Then  this sequence defines an ideal point of $\oo$.
As $\widetilde a_i \subset \widetilde l_i$, it
follows that $L = (l_i)$ is an admissible sequence for $p$
and $A \leq L$.
Since $A = (a_i)$ is a master sequence for $p$, then
given $\widetilde a_i$, there is $j > i$ with 
$l_j \cup \widetilde l_j \subset \widetilde a_i$
and so $L \leq A$.
It follows that $L = (l_i)$ is also a master sequence for
$p$.
This is case 1).

Suppose from  now on that for {\underline {any}} master sequence
$A = (a_i)$ for $p$ and {\underline {any}} $(l_i)$ as constructed above,
then $(l_i)$ does not escape compact sets of $\oo$.
Then $(l_i)$ converges to a family of non separated
line leaves in $\oos$:  \ ${\cal C} = \{ c_k,  \ k \in I
\subset {\bf Z} \}$.
If there are no perfect fits then ${\cal C}$ is a singleton
by theorem \ref{theb}.
Assume ${\cal C}$ is ordered as described in theorem \ref{theb}.
Here $I$ is either $\{ 1, ... k_0 \}$ or is ${\bf Z}$.

Choose $x_i \in b_i = a_i \cap l_i$.
These points will be used for the remainder of the proof.
Since $x_i$ is in $a_i$ and $(a_i)$ is a master sequence
for $p$, the definition of the topology in $\cd$ (definition 
\ref{topology}) implies that $x$ converges to the fixed point $p$
in $\cd$.
Here we need to differentiate between the set ${\cal C}$ of leaves and
the set $\cup {\cal C}$ of points in the leaves in ${\cal C}$.
For any $y$ in $\cup {\cal C}$, then $y \in c_k$ for some $k$ and
$\oou(y)$ intersects $l_i$ for $i$ big enough in a point denoted by
$y_i$. Similarly for $z$ in $\cup {\cal C}$ define $z_i = \oou(z) \cap l_i$.
This notation will be used for the remainder of the proof.

\vskip .1in
\noindent
{\underline {Situation 1}} $-$ Suppose there are $y, z \in \cup {\cal C}$ so
that for big enough $i$, $x_i$ is between $y_i$ and $z_i$ in $l_i$.

We refer to fig. \ref{bett}, a.
Suppose that $z$ is in $c_{j_0}$, $y$ in $c_{j_1}$, with $j_0 \leq j_1$
in the given order of ${\cal C}$.
If $j_0 = j_1$ then the segment $u_i$ of $l_i$ between
$z_i, y_i$ converges to the segment
in $\oos(z)$ between $z$ and $y$.
Then $x_i$ does not escape compact sets, contradiction
to $A = (a_i)$ being an admissible sequence.

\begin{figure}
\centeredepsfbox{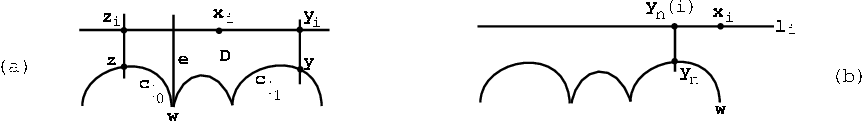}
\caption{
a. The case that $x_i$ is between some unstable leaves,
b. The case $x_i$ escapes to one side.}
\label{bett}
\end{figure}

For any $k$ the leaves $c_k, c_{k+1}$ are non separated from
each other and there is a leaf $e$ of $\oou$
making perfect fits with both $c_k$ and $c_{k+1}$.
This defines an ideal point $w$ of $\partial \oo$
which is an ideal point of equivalent rays  of 
$c_k, c_{k+1}$ and $e$, see
fig. \ref{bett}, a ($k = j_0$ in the figure).
Consider the open region $D$ of $\oo$ bounded by the ray of
$c_{j_0}$ defined by $z$ and going in the $y$ direction,
the segment in $\oou(z)$ from $z$ to $z_i$, the segment 
$u_i$ in $l_i$ from $z_i$ to $y_i$, the segment in
$\oou(y)$ from $y_i$ to $y$, the ray in $\oos(y)$ defined
by $y$ and going towards the $z$ direction and the leaves
$c_k$ with $j_0 < k < j_1$ (this last set is empty if
and only if $j_1 = j_0 + 1$).
By the remark above, the only ideal points of $D$ in $\partial \oo$,
that is the set $\overline D \cap \partial \oo$ (closure in $\cd$), \ 
are those associated to the appropriate
rays of $c_k$ with $j_0 \leq k \leq j_1$.
Since $(x_i)$ converges to $p$ which is in $\partial \oo$, then 
$p$ is one of these points.
So $p$ is an ideal point of  a ray of $\oos$ or $\oou$ which makes
a perfect fit with another leaf. There is a master sequence
which  is a standard sequence associated to $p$.
This is case 2 of the proposition.

\vskip .1in
\noindent
{\underline {Situation 2}} $-$ For any $y, z$ in $\cup {\cal C}$ the $x_i$ is eventually
not between the corresponding $y_i, z_i$.

Let $y \in \cup {\cal C}$. Then up to subsequence the $x_i$ are in one side of
$y_i$ in $l_i$, say in the side corresponding to increasing $k$ in
the order of ${\cal C}$ (this is in fact true for any big $i$
as $x_i$ converges in ${\cal D}$).

Suppose first that ${\cal C}$ is an infinite collection of non separated
leaves.
Let $w$ be the ideal point associated to the infinite collection 
${\cal C}$ and in the increasing direction of ${\cal C}$ 
as in lemma \ref{infin}.
We follow the notation of lemma \ref{infin}:
let $d_m = a_m \cup b_m$ and let
$D  = (d_m)$ be a master sequence associated to 
$w$ as described in lemma \ref{infin}.
Fix $m$. Then $x_i$ is eventually in 
$\widetilde d_m$. Therefore $x_i$ converges to $w$
and it follows that $w = p$.
Here we are in case 3).

Finally suppose that ${\cal C}$ is finite. 
Let $w$ be the ideal point of the ray of $c_{k_0}$ corresponding
to the increasing direction in ${\cal C}$.
Let $y_n$ in $c_{k_0}$ converging to $w$, see fig. \ref{bett}, b.
Let 

$$y_n(i) \ \ = \ \ \oou(y_n) \cap l_i$$

\noindent
Fix $n$. Then eventually in $i$, the
$x_i$ is in the component of 
$l_i - y_n(i)$ in the $w$ side, see fig. \ref{bett}, b.
Consider a standard sequence defining $w$ so that:
it is arbitrary in the side of $\oo - c_{k_0}$ 
not containing $x_i$ and in the other
side we have an arc 
in $\oou(y_n)$ from $y_n$ to $y_n(i)$ and 
then a ray in $l_i$ $-$ which contains $x_i$ for $i$ big.
Since $c_{k_0}$ is the biggest element in ${\cal C}$, there is no
leaf of $\oos$ non separated from $c_{k_0}$ in that side of $C_{k_0}$.
Hence the $l_i$ cannot converge (in $\oo$) to anything on that side and
those parts of $l_i$ escape in $\oo$.
As the $x_i$ are in these subarcs of $l_i$ then $x_i \rightarrow w$ in
$\cd$ and so $p = w$.

Let $r_n = \oou(y_n)$. If $r_n$ escapes compact sets in $\oo$
as $n \rightarrow \infty$, then
it defines a master sequence for $p$ and we are in
case 1). Otherwise $r_n$ converges to some $r$
making a perfect fit with $c_{k_0}$ and we are in case 2).
This finishes the proof of proposition \ref{options}.
\end{proof}

\section{Flow ideal boundary and
compactification of the universal cover}

For the remainder of the article, unless otherwise stated,
 we will only consider 
pseudo-Anosov flows without perfect fits and not conjugate
to suspension Anosov flows.
In this section we  compactify 
the universal cover $\mi$ with a sphere at infinity 
using only dynamical systems tools.


\begin{lemma}{(model pre compactification)}{}
Let $M$ be a closed $3$-manifold with a pseudo-Anosov
flow without perfect fits and not conjugate
to suspension Anosov. There is a compactication $\cd \times [-1,1]$ of $\mi$
which is 
a topological product.
\end{lemma}

\begin{proof}{}
Recall that $\cd$ is a compactification of the
orbit space  $\oo$ of $\wwp$ and $\cd$ is
homeomorphic to a closed disk. Consider 
$\cd \times [-1,1]$ with the product topology.
This is compact and homeomorphic to a closed $3$-ball.
The set $\mi$ is homeomorphic to the interior
of $\cd \times [-1,1]$ which is 
$\oo \times (-1,1)$. In fact choose a cross section
$f_1: \oo \rightarrow \mi$ and a homeomorphism
$f_2: (-1,1) \rightarrow \rrrr$. This produces
a homeomorphism

$$f:  \oo \times (-1,1) \rightarrow \mi, \ \ \ \
f(x,t) \ = \ \wwp_{f_2(t)} (f_1(t))$$

\noindent
Clearly the topology in $\mi$ is the same as
the induced topology from $\oo \times (-1,1)$.
In this way $\mi$ can be seen as an open dense subset of
$\cd \times [-1,1]$ and $\cd \times [-1,1]$ is a 
compactification of $\mi$.
\end{proof}

This  construction
is reminiscent of the one done by Cannon-Thurston
\cite{Ca-Th} for fibrations.
Notice that this construction works for any pseudo-Anosov flow,
even with perfect fits.

\vskip .1in
\noindent
{\bf {Important remark}} $-$ 
We should stress that this precompactification $\cd \times [-1,1]$ is
{\underline {far from natural}}, because in general
it is very hard to put a topology in 
$\partial \oo \times (-1,1)$ which is group
equivariant. 
In other words the section $f_1: \oo \rightarrow \mi$ is
not natural at all.
The interior of $\cd \times [-1,1]$ is
homeomorphic to $\mi$ and clearly $\pi_1(M)$ acts on this open
set.
The topology in $\cd \times \{-1, 1\}$ is what you would
expect, since it is homeomorphic to the topology
of $\cd$, which is group equivariant.
But the topology in $\partial \oo \times \mmp$
is really not well defined. 
Using the section $f_1$ we can define a trivialization
of $\partial \oo \times \mmp$, connecting it to
$\mi \cong \oo \times (-1,1)$.
The problem here 
is that given a covering translation $h$ of $\mi$,
there is no guarantee that it will extend continuously
to $\partial \oo
\times (-1,1)$ \ (but it does extend naturally and continuously
to $\cd \times
\{ -1, 1 \}$).
This problem is easily seen even in the case of suspension
pseudo-Anosov flows.
Instead of using the lift of a fiber as a section
$\oo \rightarrow \mi$, use a section which goes
one step lower (with respect to the fiber) in
certain directions.
From the point of view of the new trivialization of
$\cd \times \mmp$ certain covering translations
will not extend to $\cd \times \mmp$.

But this will not be a problem for us, because we will
collapse $\partial \cd \times [-1,1]$, identifying each vertical
interval $\{ z \} \times \mmp$ ($z$ in $\partial \oo$)
to a point. In fact one could have adjoined to $\mi$ 
just the top and bottom $\cd \times \{ -1, 1\}$.
However 
it is much easier to describe sets and neighborhoods
in the $\cd \times \mmp$ model as above, making many
arguments simpler. 
The topology of the quotient space will be
completely independent of the chosen section/trivialization 
and will depend only on the pseudo-Anosov flow.

\vskip .1in
The compactification of $\mi$ we desire will be a quotient 
of $\cd \times [-1,1]$,
where the identifications occur only in the boundary
sphere.
First we work only in the boundary of $\cd \times [-1,1]$ and later
incorporate $\mi$.

We will use a theorem of Moore concerning cellular
decompositions. A decomposition $G$ of a space
$X$ is a collection of disjoint nonempty closed sets
whose union is $X$. There is a quotient space
$X/G$ and a map $\nu: X \rightarrow X/G$. The points
of $X/G$ are just the elements of $G$. The
point $\nu(x)$ is the unique element of $G$
containing $x$. The topology in $X/G$ is the
quotient topology: a subset $U$ of $X/G$ is open
if and only if $\nu^{-1}(U)$ is open in $X$.

A decomposition $G$ of $X$ satisfies the 
{\em upper semicontinuity property} provided
that, given $g$ in $G$ and $V$ open in $X$
containing $g$, the union of those $g'$ of
$G$ contained in $V$ is an open set in $X$.
Equivalently $\nu$ is a closed map.

A decomposition $G$ of a closed $2$-manifold
$B$ is {\em cellular}, provided that 
$G$ is upper semicontinuous and provided 
each $g$ in $G$ is compact and has a non separating
embedding in the Euclidean plane $E^2$.
The following result was proved by
R. L. Moore for the case of a sphere:

\begin{theorem}{(approximating cellular maps, Moore's theorem)}{}
Let $G$ denote a cellular decomposition 
of a $2$-manifold $B$ homeomorphic to a sphere.
Then the identification map $\nu: X \rightarrow X/G$ 
can be approximated by homeomorphisms. In particular
$X$ and $X/G$ are homeomorphic.
\end{theorem}


\begin{theorem}{(flow ideal boundary)}{}
Let $\Phi$ be a pseudo-Anosov flow in $M^3$ closed
which is not topologically conjugate 
to a suspension Anosov flow and there are
no perfect fits between leaves of $\wls, \wlu$.
Let $\cd \times [-1,1]$
be the model pre
compactification of $\mi$.
Then $\partial  (\cd \times [-1,1])$ has 
a quotient $\rr$ which is a $2$-sphere where
the group $\pi_1(M)$ acts by homeomorphisms.
The space $\rr$ and its topology are completely
independent of the model precompactification
$\cd \times [-1,1]$ and depend only on the flow $\Phi$.
\label{flobo}
\end{theorem}

\begin{proof}{}
The topology
in $\cd \times \{ -1, 1 \}$ is well defined by the obvious 
bijections $\cd \rightarrow \cd \times \{ 1 \}$, \
$\cd \rightarrow \cd \times \{ -1 \}$. 
The structure of $\oos \times \{ 1 \}$ in $\cd \times \{ 1 \}$
is then equivalent to
that of $\oos$ in $\cd$, etc..
We will stress where needed that
arguments are independent of parametrization/trivialization of 
$\partial \oo \times (-1, 1)$.

We construct a cellular decomposition ${\cal R}$ of $\partial \cd \times [-1,1]$ as
follows. The cells are one of the following types:

(1) Let $l$ be a leaf of $\oos$ with ideal points
$a_1, ..., a_n$ in $\partial \oo$. Consider the cell
element 

$$g_l \ \  = \ \ l \times \{ 1 \} \ \cup \ \bigcup_{1 \leq i \leq n}
\ a_i \times [-1,1]$$

(2) Let $l$ be a leaf of $\oou$ with ideal points
$b_1, ..., b_m$ in $\partial \oo$. Consider the cell
element 

$$g_l \ \  = \ \ l \times \{ -1 \} \ \cup \ \bigcup_{1 \leq i \leq n}
\ b_i \times [-1,1]$$

(3) Let $z$ be a point of $\partial \oo$ which is not an ideal point
of a ray of $\oos$ or $\oou$. Consider the cell element
$g_z = z \times [-1,1]$.

Later on we will think of ${\cal R}$ as a set of points with
the quotient topology induced by the map from $\partial (\cd \times [-1,1])$
to ${\cal R}$.

Since every point in $\oo$ is in a leaf of $\oos$, then elements
of type (1) cover $\oo \times \{ 1 \}$. 
Similarly elements of type (2) cover $\oo \times \{ -1 \}$.
Finally elements of type (3)  cover the rest of
$\partial \oo \times \mmp$.
Cover here means the union contains the set in question.
Under the hypothesis of no perfect fits, no two rays of
$\oos$ or $\oou$ have the same ideal point. This implies that
distinct elements of type (1), (2) or (3)
 are disjoint from each
other.
This defines the decomposition ${\cal R}$ of $\partial (\cd \times [-1,1])$.

We now show that ${\cal R}$ is a cellular
 decomposition of $\partial
(\cd \times [-1,1])$.
Any element of type (3) is homeomorphic to a closed
interval, hence compact.
An element $g$ of type (1) is the union of finitely
many closed intervals in $\partial \oo \times [-1,1]$
and a set $(l \cup \pin l) \times \{ 1 \}$
in $\cd \times \{ 1 \}$.
The set $l \cup \pin l$ (contained in $\cd$) is
homemorphic to a compact $k$-prong in the plane.
Therefore $g$ is compact and homeomorphic
to $l \cup \pin l$.
In addition, any $g$ in ${\cal R}$ has a non separating embedding
in the Euclidean plane.

Next we prove that ${\cal R}$ is upper semicontinuous.
Let $g$ in ${\cal R}$ and $V$ an open set in 
$\partial (\cd \times [-1,1])$ containing
$g$.
Let $V'$ be the union of the $g'$ in ${\cal R}$ with
$g' \subset V$. We need to show that $V'$ is
open in $\partial (\cd \times [-1,1])$.
Since $g$ is arbitrary 
it suffices to show that $V'$ contains an open
neighborhood of $g$ in $\partial (\cd \times [-1,1])$.
We do the proof for elements of type (1) see fig. \ref{decomp}, 
the other
cases being very similar.

\begin{figure}
\centeredepsfbox{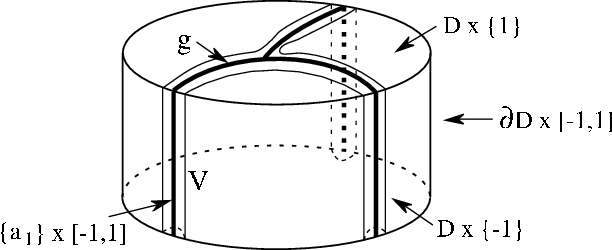}
\caption{An element of $g$ type (1) in $\partial \cd \times [-1,1]$ and
a neighborhood of it.}
\label{decomp}
\end{figure}


Let $g$ be generated by the leaf $l$ of $\oos$, 
let $a_1, ..., a_n$ be the ideal points
of $l$ in $\partial \oo$. 
For each $i$
there is a neighborhood $J_i$ of $a_i$ in
$\partial \oo$ with $J_i \times [-1,1]$ contained
in $V$.
This is because $\partial \oo \times [-1,1]$ is
homeomorphic to a closed annulus.
This conclusion is independent of the parametrization we choose
for 
$\partial \oo \times\mmp$.

Let $(p_k)_{k \in {\bf N}}$ be a sequence of points in 
$\partial (\cd \times \mmp)$
converging to some point
$p$ in $g$.
Let $g_k$ be the element of ${\cal R}$ containing $p_k$.
We show that for $k$ big enough, then $g_k$ is
contained in $V$ and therefore $p_k$ has to be contained
in $V'$. Hence $V'$ contains
an open neighborhood of $g$ in $\partial \cd \times [-1,1]$.
This will prove the upper semicontinuity property of the cellular
decomposition.

Up to a subsequence we may assume that 
all $p_k$ are either  in \ \
I) $\oo \times \{ 1 \}$,
\  \ II)  $\oo \times \{ -1 \}$ or \ \ III) $\partial \oo \times \mmp$.
We analyse each case separately:

\vskip .1in
\noindent 
{\underline {Case I}} $-$ Suppose first that $p_k$ is in $\oo \times \{ 1 \}$.

Hence $p_k \in \cd \times\{ 1 \}$.
Up to subsequence and reordering $\{ a_i \}$,
assume that $p_k$ are in 
a sector of $l \times \{ 1 \}$ defined by $b \times \{ 1 \}$
where $b$ is a line leaf of $l$ with ideal points $a_1, a_2$.
Then $g_k$ is an element of type (1) and is the union 

$$g_k \ = \ l_k \times \{ 1 \} \  \cup \ 
\cup _j \ (\{ w_{kj}\} \times \mmp)$$

\noindent
where $w_{kj}, 1 \leq j \leq j_0$ are the 
ideal points of $l_k$, a leaf of $\oos$.
Notice that $g_k$ is contained in the set
$(\cd \times \{ 1 \}) 
\cup (\partial \oo \times \mmp)$.
%

We need the following result which is also useful later.
It shows the strength of the no perfect fits hypothesis.

\begin{lemma}{(the escape lemma)}{}
Let $\Phi$ be a pseudo-Anosov without perfect fits and not
conjugate to a suspension Anosov flow.

i) Let $(l_n)_{n \in {\bf N}}$ 
be a sequence of 
leaves 
or slices of leaves of (say) $\oos$. 
Suppose that $(l_n)$ converges to a line leaf $l$
of (say) $\oos$.
It follows that the ideal points
of $l_n$ converge to the ideal points of $l$,

ii) Under the hypothesis of i),  
if $x_{n_k} \in l_{n_k}$ converges to $x$ in $\cd$, then
$x$ is in $l \cup \partial l$.

iii) Let $l_n$ in $\oos$ or $\oou$.
Suppose there are $x_n, y_n$ in $l_n \cup \partial l_n$ so that
$x_n, y_n$ converge to distinct points of $\partial \oo$. Then $l_n$ converges to
a leaf $l$. In particular $l_n$ does not escape compact sets
in $\oo$.

\label{duo}
\end{lemma}

\begin{proof}{}
Suppose i) is not true. 
Let $p$ in $l$.
Hence there is an ideal point $a_1$ of $l$ in $\partial \oo$ and 
there are $r_n$ rays of $l_n$
starting at $p_n$ and in the direction of the ray in $l$ with
ideal point $a_1$ so that: $b_n = \partial r_n$ does not
converge to $a_1$. 
This also works up to subsequences.
We may assume that $(l_n)$ is nested.
By separation properties $b_n$ is weakly monotone in $\partial \oo$
and converges to a point $c \not = a_1$.
Consider the interval $(c,a_1)$ of $\partial \oo$ not containing $a_2$.
Suppose first that this interval has an ideal point of a
leaf $e$ of $\oos$ or $\oou$.
The leaf $e$ is a barrier for the leaves $l_n$,
so this implies that $l_n$ also converges to another leaf besides $l$.
Since there are no perfect fits, this is impossible by
theorem \ref{theb}.
We are left with the possibility that $(c,a_1)$ does
not have an ideal point of a leaf of $\oos$ or $\oou$.
But this is also impossible: let $z$ in $(c,a_1)$.
If $z$ is ideal point of leaf of $\oos$ or $\oou$ we are
done. 
Since there are no perfect fits, option 1) of
proposition \ref{options} has to occur and
there is a neighborhood
system of $z$ defined by a sequence of stable leaves.
This shows that any neighborhood of $z$ in $\partial \oo$ has points
which are ideal points of leaves of $\oos$.
These arguments show that these ideal points of $l_n$ converge
to $a_1$. This proves i).

Proof of ii). Up to taking a subsequence we assume the
statement is for 
$x_n$ in $l_n$.
Since the leaf space of $\oos$ is Hausdorff, if $x$ is in $\oo$ then
$x$ is in $l$.
Suppose that $x$ is in $\partial \oo$.
Using the notation from part i) suppose that $x_n$ are in
the rays $r_n$ as in part i). Suppose that $x_n$ does not
converge to $a_1$ and instead converges to $c \not = a_1$.
Let $U, V$ small disjoint neighborhoods of $a_1, c$. 
By conclusion i) already proved, for $n$ big
$r_n$ has ideal point in $U$. Fix one such $n$ 
 and so $r_n$ is
entirely in $U$ except for an initial compact segment $t$.
For any $m > n$ the $r_m$ is constricted to be in the union
of two sets $S_1$ and $S_2$: 
1) $S_1$ the compact region of $\oo$ which is
bounded by a polygon made of 4 arcs: A) $t$, B) a compact
arc $l'$ in $l$ from $p$ to a point in $U$, C) a compact arc in 
$U$ from the the end of $l'$ to the end of $t$ and
D) a very small arc from the beginning of $t$ to the 
beginning of $l$;  
\ 2) the second set $S_2 = U$.
Since $(x_m)$ escapes compact sets in $\oo$, then
for big $m$, $x_m$ cannot
be in $S_1$ so it has to be in $S_2 = U$, contradiction
to $x_m \in V$. This shows
that $x_i$ converges to $a_1$.

Proof of iii). 
Wlog assume that $l_n$ are leaves of $\oos$.
Let $x_n, y_n$ converging to distinct points $x, y$ of 
$\partial \oo$. 
If $x_n$ is in $\partial \oo$ one can choose a point in $l_n$ arbitrarily
near $x_n$, so we may assume that all $x_n, y_n$ are in $\oo$.
Let $r_n$ be the arc in $l_n$ from $x_n$ to $y_n$. If the
sequence $(r_n)$ escapes compact sets in $\oo$, 
then it limits on at least one of the intervals
$(x,y)$ or $(y,x)$ both of which are non degenerate.
But that would imply that this interval does not contain
an ideal point of a ray of $\oos$ or $\oou$ $-$ this was proved
to be impossible in the proof of part i).
Since $l_n$ does not escape compact sets in $\oo$, there is
a subsequence $(l_{n_k})$ and $p_{n_k}$ in $l_{n_k}$ with
$p_{n_k}$ converging to a point $p$ in $\oo$. Let $l$ be the
leaf of $\oos$, with $p \in l$. Hence the sequence
$(l_{n_k})$  converges to $l$ (and to no other leaf
when there are no perfect fits).
Since $x_{n_k}$ is in $l_{n_k}$ and converges to $x$, part ii)
shows that $x$ is an ideal point of $l$ and so is $y$.
Notice in addition
that in the case of no perfect fits there is only one
leaf of $\oos$ with ideal point $x$.
But these arguments can be applied to 
{\underline {any}} subsequence  of $(l_n)$ to show that such
a subsequence has another subsequence converging to a leaf $l'$ which
has an ideal point $x$. But as remarked before this implies
that $l = l'$. It follows that the original sequence $(l_n)$
has to converge to $l$. This finishes the proof of iii).
\end{proof} 


Notice that conclusion iii) is false for suspension Anosov flows.

\vskip .1in
\noindent
{\underline {Continuation of the proof of theorem \ref{flobo} }}

Recall the setup in case I: $p_k \in \oo \times \{ 1 \}$ converge to
$p$ in $g = g_l$. The $p_k$ are in $l_k \times \{ 1 \}$ with
$l_k$ all in a sector of $b$ line leaf of $l$ with ideal points
$a_1, a_2$; 
\ $V$ is a neighborhood of $g$ in $\partial (\cd \times [-1,1])$.
Let $p_k = y_k \times \{ 1 \}$.

\vskip .1in
\noindent
{\underline {Case I.1}} $-$ Suppose $p \in \oo \times \{ 1 \}$.

Then $p_k$ converges to $p = y \times \{ 1 \}$.
By lemma \ref{duo} part ii) any limit point of $x_{n_k}$ with
$x_{n_k}$ in $l_{n_k}$ is in $b \cup \{ a_1, a_2 \}$. Hence 
$l_n \times \{ 1 \} \subset V$ for $n$ big.
Lemma \ref{duo} part i), the ideal points of rays in $l_k$ also
converge to $a_1$ or $a_2$ and so $w_{kj} \times [-1,1] \subset V$ for
$k$ big. It follows that $g_k \subset V$ for $k$ big in this case.

\vskip .1in
\noindent
{\underline {Case I.2}} $-$ 
Suppose $p \in \partial \oo \times \{ 1 \}$.

Wlog assume that $p$ is $a_1 \times \{ 1 \}$.
In this case suppose first that $(l_k)$ does not escape 
compact sets in $\oo$. Assume up to subsequence
that $(l_k)$ converges to a line leaf $s$ of $\oos$.
Then we may assume that $l_k$ is nested. Since there are no perfect fits,
there is only one such leaf $s$ in the limit. 
As $p_k \in l_k \times \{ 1 \}$, lemma \ref{duo} part ii) shows
that the limit of $y_k$ is an ideal point of a ray of $s$.
This limit is $a_1$ so $a_1$ is an ideal point of $s$.
This shows that $s, l$ have rays with same ideal points.
By definition the rays are equivalent.
But since there are no perfect fits, then $s = l$.
This reduces the proof to case I.1.

Finally we suppose that $(l_k)$ escapes compact sets in $\oo$.
Since there are $p_k$ in $l_k \times \{ 1 \}$ converging to $a_1 \times \{ 1 \}$
we claim that $g_k \cap (\cd \times \{ 1 \})$ converges to $a_1 \times \{ 1 \}$.
Otherwise up to subsequence there are $z_k$ in $l_k$  with
$z_k$ converging to $v \not = a_1$. Hence $l_k$ has arcs
with endpoints in $y_k \rightarrow a_1$ and $z_k \rightarrow v$.
The escape lemma (lemma \ref{duo} part i)  implies that
$l_k$ does not escape compact sets in $\oo$, contradiction.
This finishes the analysis of case I).

%

The next case in the proof of theorem \ref{flobo}
is:

\vskip .1in
\noindent
{\underline { Case II}} $-$  Suppose that $p_k$ is in $\oo \times \{ -1 \}$.

There is an asymmetry here because in case I, $g$ and $g_k$ are cells
to type (1), whereas in case II, $g$ is of type (1) and $g_k$ is
of type (2). So we cannot just revert the direction of the flow and
use the proof of case I to prove case II.

In this case 
it follows that $g_k$ is contained in \  $(\cd \times \{ -1 \})
 \cup  (\partial \oo \times \mmp)$. \
Since $p_k$ converges to $p$ in $g$ and $g$ is
contained in \  $(\cd \times \{ 1 \}) \cup (\partial \oo 
\times \mmp)$, it follows that $p$ is
in $\partial \oo \times \{ -1 \}$ and
$p$ is say $(a_1,-1)$, where $a_1$ is one of the
ideal points of $l$.

Here $a_1$ is an ideal point of a ray in $\oos$ and there are no perfect
fits and no non separated leaves of $\oos$ or $\oou$.
Therefore proposition \ref{options} shows that
there is a neighborhood system of $a_1$ in
$\cd$ defined by a sequence $(r_n)_{n \in {\bf N}}$
of {\underline {unstable}} leaves (this is option (1) in 
proposition \ref{options}).
Since $V$ is open in $\partial (\cd \times [-1,1])$, then
for $n$ big enough $r_n \times \{ -1 \}$ is contained
in $V$.
The element $g_k$ is of the form 

$$(s_k \times \{ -1 \}) \ \ \cup  \ \ \cup_j \ (\{ b_{kj} \} \times \mmp)$$

\noindent
where the $s_k$ are leaves of $\oou$ with points
converging to $a_1$. 
Since both $r_n$ and $s_k$ are unstable leaves, they cannot
intersect transverely. 
It now follows that
for $k$ big enough $s_k \times \{ -1 \}$ is contained
in $V$. 
There is an interval $J$ in $\partial \oo$ with $a_1$ in the interior
of $J$ and with $J \times [-1,1] \subset V$ $-$ this is because $V$ is
open and $a_1 \times [-1,1]$ is contained in $V$.
Hence the endpoints $b_j$ have to be in $J$ for $k$ big enough.
It follows that $g_k$ is entirely contained in $V$.
This finishes the analysis of case II.

\vskip .1in
\noindent
{\underline {Case III}} $-$   Suppose that $p_k$ is in $\partial \oo \times \mmp$.

Then $p_k$ converges to 
$p = (c,t)$ where $c$ is in $\partial \oo$.
Hence $V$ contains $J \times [-1,1]$ for some interval $J$ in $\partial \oo$,
so that $J$ contains $c$ in its interior.
Here $g_k$ can be type (1), (2) or (3).
If $g_k$ is of type (3) then for $k$ big enough the
$g_k$ is contained in $J \times \mmp$ and
hence in $V$.

If $g_k$ is of type (2), then 
it has vertical stalks $b_{kj} \times [-1,1]$ which
are eventually contained in $J \times [-1,1]$.
Hence $b_{kj}$ is an ideal point of a leaf $s_k$ in
$\oou$. As $k$ varies, one of the ideal points of 
$s_k$ (namely $b_{kj}$) converges to $a_1$, which is an ideal point of $l$.
The proof then proceeds as in case II to show that eventually
$g_k$ is entirely contained in $V$.

Finally if $g_k$ is of type (1), then 
as seen in part I), $g_k$ is
contained in $V$ for $k$ big enough.

This proves that $V'$ is open.
We conclude that the decomposition satisfies the
upper semicontinuity property.
By Moore's theorem it follows that $\rr$ is a
sphere.

%

\vskip .1in
So far we have not really used the topology in $\partial \oo \times
\mmp$.
We still need to 
show that the topology of $\rr$ is
independent of the choice of the trivialization
$\partial \oo \times \mmp$ and that the fundamental
group acts naturally by homeomorphisms on $\rr$.

To see the first statement, notice that
the quotient  map
$\partial \cd \times [-1,1]   =  \partial (\cd \times \mmp) \ \rightarrow \rr$
can be done in two steps: first collapse
each vertical stalk $\{ z \} \times \mmp$ to a
point where $z$ is in $\partial \oo$ and then
do the remaining collapsing of leaves of $\oou$ in 
$\cd \times \{ -1 \}$
and leaves of $\oos$ in  $\cd \times \{ 1 \}$.
After the first collapsing we 
have $\cd \times \{ 1 \}$ union $\cd \times \{ -1 \}$
glued along the points $\{ w \}  \times \{ -1, 1\}$.
The topology now is completely determined since the
topology on the top $\cd \times \{ 1 \}$ and the bottom
$\cd \times \{ -1 \}$ is completely determined
by the topology in $\cd$.
The fundamental group acts by homeomorphisms
in this object and preserves the foliations
stable on the top and unstable on the bottom.
Therefore the second collapse produces a sphere 
$\rr = \partial \cd \times [-1,1]/{\cal R}$.
The topology in $\rr$ is independent of any choices.
The fundamental group acts by homeomorphisms
on the quotient space $\rr$, since after the first
collapse it acts by homeomorphims and preserves
the elements of the decomposition.
This finishes the proof of the theorem \ref{flobo}.
\end{proof}


We now show that the action of $\pi_1(M)$
in $\rr$ has excellent properties, that is,
it is a 
uniform convergence group action.
A topological space $X$ is a {\em compactum}
if it is a compact Hausdorff topological space.
%
Let $X$ be a compactum and $\Gamma$ a group
acting by homeomorphisms  on $X$.
Let $\Theta_3(X)$ be the space of distinct triples of
$X$ with the subspace topology the product space
$X \times X \times X$.
Then $\Theta_3(X)$ is locally compact
and there is an induced action of $\Gamma$ on $\Theta_3(X)$.
Here {\em local uniform convergence} means uniform convergence in
compact sets.
For simplicity we state results for 
$X$ metrisable (in the general case one uses nets instead of
sequences \cite{Bo2}).
Notice we identify the group with the action.

\begin{define}{(\cite{Ge-Ma})}{}
$\Gamma$ is a convergence group if the
following holds:
If $(\gamma_i)_{i \in {\bf N}}$ is an infinite sequence of distinct
elements of $\Gamma$, then one can find points $a, b$ in $X$
and a subsequence $(\gamma_{i_k})_{k \in {\bf N}}$ of
$(\gamma_i)$, such that the maps $\gamma_{i_k} |_{X - \{ a \}}$
converge locally uniformly to the constant map with value $b$.
\end{define}

Notice that it is not necessary that $a, b$ are distinct,
which in fact does not happen always. 
It is simple to see that 
this is equivalent to the following property:
the action of $\Gamma$ on $\Theta_3(X)$ is 
{\em properly discontinuous} \cite{Tu2,Bo2}.
This means that 
for any compact subset $K$ of $\Theta_3(X)$, 
the set $\{ \gamma \in \Gamma \ | \ \gamma K \cap K \not 
= \emptyset \}$ is finite \cite{Tu2,Bo2}.
The action of $\Gamma$ is {\em cocompact} if $\Theta_3(X)/\Gamma$ is 
a compact space.
If the action is a convergence group and cocompact it is
called an {\em uniform convergence group} action.

\begin{define}{(conical limit points)}{}
Let $\Gamma$ be a group action on a metrisable compactum $X$.
A point $z$ in $X$ is a conical limit point for the action of
$\Gamma$
if there are distinct points $a, b$ of $X$ and a sequence
$(\gamma_i)_{i \in {\bf N}}$ in $\Gamma$ such that
$\gamma_i z \rightarrow a$ and $\gamma_i y \rightarrow b$
for all $y$ in $X - \{ z \}$.
\end{define}

Here it is crucial that $a, b$ are {\underline {distinct}}
for otherwise the convergence group property would
yield the result for many points.
Basic references for conical limit points are 
\cite{Tu3,Bo2}.
It is a simple result that if $\Gamma$ is a uniform convergence
group action then every 
point of $X$ is a conical limit point \cite{Tu3,Bo2}. 
The opposite implication is highly non trivial and was proved independently 
by Tukia \cite{Tu3} and Bowditch \cite{Bo1}.
Recall that $X$ is perfect if it has no isolated points.

\begin{theorem}{(\cite{Tu3,Bo1})}{}
Suppose that $X$ is a perfect, metrisable compactum
and that $\Gamma$ is a convergence group
action on $X$.
If every point of $X$ is a conical limit point for
the action, 
then $\Gamma$ is a cocompact action.
Consequently $\Gamma$ is a uniform convergence group action.
\label{coni}
\end{theorem}

Hence both properties of uniform convergence group action
can be checked by analysing sequences
of elements of $\Gamma$.
Our main technical result is the following:

\begin{theorem}{}{}
Let $\Phi$ be a pseudo-Anosov flow in $M^3$ closed so
that $\Phi$ does not have perfect fits and 
is not topologically conjugate to
a suspension Anosov flow.
Consider the induced quotient $\rr$ of 
$\partial (\cd \times \mmp)$
and the induced action of $\Gamma = \pi_1(M)$ on $\rr$.
Then $\Gamma$ is a uniform convergence group. 
\label{conver}
\end{theorem}

We first prove that $\pi_1(M)$ acts as a convergence group on
$\rr$ using the sequences formulation
and then we show that every point of $\rr$ is a
conical limit point for the action of $\Gamma$ on $\rr$.
The space $\rr$ is homeomorphic to a sphere, hence
it is a perfect, metrisable compact space and theorem \ref{coni}
can be used.

First we define an important  map which will be used throughout
the proofs in this section.
Recall there is a continuous quotient map
$\nu: \partial (\cd \times \mmp) \rightarrow \rr$.
Identify $\partial \oo$ with $\partial \oo \times \{ 1 \}$
by $z \rightarrow (z,1)$ 
in $\partial (\cd \times [-1,1])$. Then there is an induced map:

$$\varphi: \ \partial \oo \rightarrow \rr, \ \ \
\varphi(z) \ = \ \nu((z,1)) \ \ \ \ (*)$$

The map $\varphi$ is continuous.
Every $g$ of ${\cal R}$ contains intervals
of the form $\{ y \} \times \mmp$ where $ y \in \partial \oo$, 
so $\varphi$ is surjective.
Hence $\varphi$ encodes all of the information of the map
$\nu$. 
In addition $\pi_1(M)$ acts on $\partial \oo$.
The proof will use deep knowledge 
about the action of $\Gamma = \pi_1(M)$ on the
circle ${\bf S}^1 = \partial \oo$ in
order to obtain information about the action
of $\Gamma$ on $\rr$.

Notice that the map $\varphi$ is  group equivariant producing examples 
of group invariant sphere filling curves.

\vskip .1in
\noindent
{\bf {Remarks}} $-$ 1) A very important fact is the following.
Suppose that $x, y$ distinct in $\partial \oo$ are identified
under $\varphi$, that is $\varphi(x) = \varphi(y)$.
Because of the no perfect fits condition, there are no
distinct leaves of $\oos, \oou$ sharing an ideal point
in $\partial \oo$. This implies there is a leaf $l$ of
$\oos$ or $\oou$ so that $x, y$ are ideal
points of $l$.
In particular there are at most $k$ preimages 
under $\varphi$ of any point, where $k$ is the
maximum number of prongs at a singular point
of $\oos$ or $\oou$.

2) (important convention)  
Recall that $\hhs, \hhu$ are the leaf spaces of
$\wls, \wlu$ respectively.
\ If $\gamma$ is an element of $\pi_1(M)$ then $\gamma$
acts as a homeomorphism in all of the spaces
$\mi$, $\oo$, $\partial \oo, \rr, \hhs$ and $\hhu$.
For simplicity,
the same notation $\gamma$ will be used for all
of these homeomorphisms. The context will make
it clear which case is in question.
With this understanding, the fact that $\varphi$
is group equivariant means that for any
$\gamma$ in $\pi_1(M)$ then 

$$\gamma \circ \varphi  \ = \ \varphi \circ \gamma$$

\noindent
where the first $\gamma$ acts on $\rr$ and the
second acts on $\partial \oo$. 
The reader should be aware that 
this convention will be used throughout this section.

\vskip .1in
Recall that if $l$ is a ray or leaf
of $\oos$ or $\oou$, then $\partial l$ denotes the
ideal point(s) of $l$ in $\partial \oo$.
Before proving theorem \ref{conver}, we first show
in the next 2 lemmas that for any
$\gamma$ in $\pi_1(M)$, the
action of $\gamma$ on $\partial \oo$ and
$\rr$ is as expected.
In the first lemma we do not assume that there are
no perfect fits.

\begin{lemma}{}{}
Suppose that $\Phi$ is pseudo-Anosov flow 
not conjugate to suspension Anosov.
Suppose there is no infinite collection of leaves of
$\wls$ or $\wlu$ which are all non separated from
each other.
Let $\gamma$ in $\pi_1(M)$ with no fixed points
in $\oo$. Then the action of $\gamma$ on
$\partial \oo$ either 1) has only 2 fixed points
one attracting and one repelling 
and is of hyperbolic type or 2) it 
has
a single fixed point in $\partial \oo$, which is
of parabolic type.
In the second case, the fixed point of $\gamma$ is a parabolic
point in $\partial \oo$ associated to a perfect fit
horoball.
Finally if there are no perfect fits  only option 1) can occur.
\label{ideal}
\end{lemma}

\begin{proof}{}
If $\gamma$ leaves invariant a leaf $F$ in $\hhs$, then 
there is an orbit $\widetilde \alpha$ in  $F$ with
$\gamma(\widetilde \alpha) = \widetilde \alpha$.
Then $\gamma$ does not act freely on $\oo$, contradiction.

The space $\hhs$ is what is called
a non Hausdorff tree \cite{Fe6,Ro-St}.
Very roughly a non Hausdorff tree is a ``one-dimensional" space
with a tree like behavior, except that one allows non Hausdorff 
behavior. It is simply connected and is the union of countably
many ``segments".
Since $\gamma$ acts freely on $\hhs$ then theorem A of
\cite{Fe6} implies that $\gamma$ has a translation
axis for its action in $\hhs$.
The transformation $g$ leaves invariant this axis and acts
as a translation on it.
The points in the
axis are exactly those leaves $L$ of $\wls$ so
that $\gamma(L)$ separates $L$ from $\gamma^2(L)$.
This implies that the $\{ \gamma^n(L), \ n \in {\bf Z}
\}$
form a nested collection of leaves.

As explained in \cite{Fe6},
the axis does not have to be properly embedded
in $\hhs$ $-$ that is there may be
$a_n$ in the axis, escaping in the axis, but not
escaping compact sets in the leaf space $\hhs$.
Let $L$ be in the axis.
If \ $(\gamma^n(L))_{n \in {\bf N}}$ \ 
does not escape compact sets in
$\oo$, then by the nested property, the
$\gamma^n(L)$ converges to some $F$ in $\wls$
as $n$ converges to infinity.
If $\gamma(F) = F$ we have an invariant leaf in $\wls$,
contradiction. If $\gamma(F), F$ are distinct 
let ${\cal B}$ be the set of leaves of $\wls$
non separated from $F$ in the side the 
$\gamma^n(L)$ are limiting to.
By theorem 
\ref{theb}, the set ${\cal B}$ is order isomorphic
to either ${\bf Z}$ or $\{ 1, ..., j \}$ for some
$j$. The first option is disallowed by hypothesis.
Consider the second option. The transformation
$\gamma$ leaves ${\cal B}$ invariant. If $\gamma$ preserves
the order in ${\cal B}$ then as ${\cal B}$ is finite,
$\gamma$ will have
invariant leaves in $\wls$, contradiction.
If $\gamma$ reverses order in ${\cal B}$, then 
there are consecutive elements $F_0, F_1$ in
${\cal B}$ which are swaped by $\gamma$.
There is a unique unstable leaf $E$ which
separates $F_0$ from $F_1$. This $E$ makes a perfect
fit with both $F_0$ and $F_1$, see theorem \ref{theb}.
By the above this leaf $E$ is invariant under
$\gamma$ again leading to a contradiction.
This argument shows that the axis of $\gamma$ is properly embedded
in $\hhs$.

Let $L_0$ be in the axis of $\gamma$ acting on $\hhs$.
As $\gamma(L_0)$ separates $L_0$ from $\gamma^2(L_0)$, 
there is a unique line leaf $L$ of $L_0$ so that the sector
defined by $L$ contains
$\gamma(L_0)$
(if $L_0$ is nonsingular then $L = L_0$).
Recall that $\Theta: \mi \rightarrow \oo$ is the projection
map: it sends a point $x$ in $\mi$ to the orbit
of $\wwp$ containing $x$.
Then \ $(\gamma^n(\Theta(L)))_{n \in {\bf N}}$ \ 
is a nested sequence
of convex polygonal paths, which escapes in $\oo$.
Hence this sequence defines a unique ideal point $b$ in
$\partial \oo$.
Similarly 
\ $(\gamma^{-n}(\Theta(L)))_{n \in {\bf N}}$ \ defines an ideal
point $a$ in $\partial \oo$.
Notice that 

$$\gamma^n(\partial \Theta(L)) \ \rightarrow \ b \  \ {\rm as} \ \
n \rightarrow \infty \ \ \ {\rm and} \ \ \ 
\gamma^n(\partial \Theta(L)) \ \rightarrow \ a \ \ {\rm as} \ \  
n \rightarrow -\infty \ \ \ (**)$$

\noindent
Clearly $\gamma(a) = a, \ \gamma(b) = b$.
For any other $z$ in $\partial \oo$, then
either $z$ is an ideal point of some $\gamma^n(\theta(L))$ or
$z$ is in an interval of $\partial \oo$ defined by
ideal points of $\gamma^n(\Theta(L))$ and 
$\gamma^{n+1}(\Theta(L))$ for
some $n$ in ${\bf Z}$.
It follows that 
property (**) above also holds for $z$.
%
\noindent

If $a, b$ are distinct then the above shows that $a, b$ form
a source/sink pair for $\gamma$ and $\gamma$ has hyperbolic dynamics
in the circle $\partial \oo$.

If $a = b$ then $\gamma$ has parabolic dynamics in $\partial \oo$
with $a$ its unique fixed point.
In addition $\Theta(L)$ has a ray $l$ with ideal point $a$.
The collection $\{ \gamma^n(l) \}_{n \in {\bf Z}}$
of pairwise distinct rays
all have ideal point $a$. 
By lemma \ref{equiva}
any two elements in this collection
are connected by a chain of perfect fits.
Then $\{ \gamma^n(l) \}_{n \in {\bf Z}}$ is an infinite
perfect fit and is associated to a perfect fit horoball.
The perfect fit horoball is invariant under $\gamma$.
This finishes the proof of the lemma.
\end{proof}

Notice that if there are infinitely many leaves of $\wls$
or $\wlu$ not separated from each other, then
there are covering translations acting freely on $\oo$ and
leaving invariant a scalloped region ${\cal S}$, see \cite{Fe3}.
If $\gamma$ is of this type then $\gamma$ will fix the 4 ideal
points in $\partial \oo$ associated to the scalloped
region ${\cal S}$.
Hence the hypothesis
in lemma \ref{ideal} is needed and this is the only additional
possibility that can occur in general:
if the $\gamma^n(L)$ does not escape compact sets for either
$n \rightarrow \infty$ or $n \rightarrow -\infty$, then the 
proof in the lemma shows that $\gamma^n(L)$ converges
to a bi-infinite collection of leaves non separated from
each other. Let ${\cal S}$ be the associated scalloped
region. Here $\gamma$ acts a translation in each collection
of non separated leaves in $\partial {\cal S}$. It follows
that $\gamma$ has exactly $4$ fixed points in $\partial \oo$.
Finally if $\gamma$ has a fixed point in $\oo$, then there
are many more possibilities for the set of fixed points
in $\partial \oo$, in particular it can be infinite.

\begin{lemma}{}{}
Suppose that $\Phi$ does not have perfect fits and 
is not conjugate to suspension Anosov.
For each $\gamma \not =$ id in $\pi_1(M)$, there are distinct $y, x$ in 
$\rr$ which are the only fixed points of $\gamma$ in $\rr$
and $x, y$ 
form a source/sink pair ($y$ is repelling, $x$ is attracting).
\label{soursi}
\end{lemma}

\begin{proof}{}
As with almost all the proofs in this section,
the proof will be a strong interplay between the
pseudo-Anosov dynamics
action on $\partial \oo$ and the induced action
on $\rr$.
By remark 1) the only identifications of the map 
$\varphi$ come from the ideal points of leaves of
$\oos$ or $\oou$.

Any $\gamma$ in $\pi_1(M)$ has at most one
fixed point in $\oo$: 
if $\gamma$ fixes 2 points in $\oo$,
then it produces 2 closed orbits of $\Phi$ which
are freely homotopic to each other (or maybe
freely homotopic to the inverse of each other
or certain powers).
By theorem \ref{chain}, the lifts of the closed orbits
are connected by a chain of lozenges and
 this 
produces perfect fits in the universal
cover $-$ disallowed by hypothesis.

Suppose first that $\gamma$ is associated to a periodic orbit
of $\Phi$ $-$ singular or not. Also $\gamma$ need not 
correspond to an indivisible
closed orbit.
Let $\beta$ be the orbit of $\wwp$ with $\gamma(\beta) = \beta$ and 
$b = \Theta(\beta)$ be the single fixed point of $\gamma$ in $\oo$.
Suppose without loss of generality
that $\gamma$ is associated to an orbit of $\Phi$ being
traversed in the forward direction.
We will show that the set of fixed points of $\gamma$
(or a power of $\gamma$)
in $\partial \oo$ is the union 
$\partial \oos(b) \cup \partial \oou(b)$ and also that
$\partial \oos(b)$ is the set of attracting fixed points for $\gamma$ and
$\oou(b)$ is the set of repelling fixed points
of $\gamma$.

\begin{figure}
\centeredepsfbox{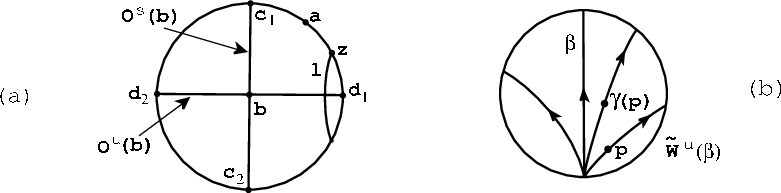}
\caption{a. The action of $\gamma$ in $\cd$ and $\partial \oo$,
b. Action of $\gamma$ in $\wu(\beta)$.}
\label{disp}
\end{figure}

Assume first that $\gamma$ leaves invariant the prongs
of $\oou(b)$ and $\oos(b)$ and that $\gamma$ is
nonsingular. 
Let $c_1, c_2$ in $\partial \oo$ be the  ideal points
of $\oos(b)$ and $d_1, d_2$ the ideal points of $\oou(b)$,
see fig. \ref{disp}, a.

Notice that $\varphi(c_1) = \varphi(c_2)$ and similarly
$\varphi(d_1) = \varphi(d_2)$. Let $x = \varphi(c_1), y
= \varphi(d_1)$. 

Clearly $\gamma$ fixes $c_1, c_2, d_1, d_2$.
Let $I$ be the interval of $\partial \oo$ with
endpoints $c_1, d_1$ and not containing $d_2$.
Since there are no perfect fits, then option (1) of proposition
\ref{options} has to occur. As $d_1$ is an ideal point
of a unstable leaf, then $d_1$ has a neighborhood
system in $\cd$ formed by stable leaves, all of which
have to intersect $\oou(b)$.
Let $l$ be one such leaf with ideal point $z$ in $I$.

The action of $\gamma$ in the set of orbits of 
$\wu(\beta)$ is contracting, 
see fig. \ref{disp}, b.
This is because $\gamma$ is associated with the 
forward flow direction.
Therefore $\gamma^n(l)$ converges
to $\oos(b)$ as $n$ converges to infinity.
It follows that $\gamma^n(z)$ converges to $c_1$ and
so $\gamma^n(a)$ converges to $c_1$. This shows that
$\gamma$ has only 2 fixed points in $I$ and
$d_1$ is repelling, $c_1$ is attracting.
The other intervals of $\partial \oo$ defined by
$\partial \oos(b) \cup \partial \oou(b)$ are
treated in the same fashion.


We claim that $y, x$ form the source/sink pair for the
action of $\gamma$ in $\rr$.
Here 

$$\gamma(x) \ = \ \gamma(\varphi(c_1)) \ = \ \varphi(\gamma(c_1))
\ = \ \varphi(c_1) \ = \ x$$

\noindent
and similarly $\gamma$ fixes $y$.
For any other $w$ in $\rr$
there is $z$ in $\oo - \{ c_1, c_2, d_1, d_2 \}$
with $w = \varphi(z)$. 
Without loss of generality assume that $z$ is in $I$.
Then $\gamma^n(z)$ converges to $c_1$ and

$$\gamma^n(w) \ = \ \gamma^n(\varphi(z)) \ = \ 
\varphi(\gamma^n(z)) \ \rightarrow \  \gamma(c_1) \ = \ x$$

\noindent
Similarly $\gamma^n(w) \rightarrow y$ when $n \rightarrow -\infty$.
So if $\gamma$ leaves invariant the components
of $\partial \oo - \{ c_1, c_2, d_1, d_2 \}$ then
$y, x$ for a source/sink pair for the action of $\gamma$
in $\rr$.

In the general case take a power of $\gamma$ so that in $\partial \oo$
it fixes all points in $\partial \oos(b), \partial \oou(b)$
and preserves orientation in $\partial \oo$.
Then
apply the above arguments.
The arguments show 
that, as a set, $\partial \oos(b)$ is invariant
and attracting for the action of $\gamma$ in $\partial \oo$
and $\partial \oou(b)$ is invariant and repelling
for the action.
All the points in $\partial \oos(b)$ are mapped to $x$
by $\varphi$ and all points in $\partial \oou(b)$ are
mapped to $y$. Hence $y, x$ is the source/sink pair
for the action of $\gamma$ in $\rr$.
This finishes the analysis of the case when $\gamma$
does not act freely in $\oo$.

We now analyse the case that $\gamma$ acts freely in $\oo$.
The previous lemma produces $a, b$ which are a 
source/sink pair for the action of $\gamma$ 
on $\partial \oo$.
Since there are no perfect fits, 
the previous lemma shows that $a \not = b$.
In fact the arguments of the previous lemma show that none of $a, b$
can be the ideal point of a ray of a leaf
of $\oos$ or $\oou$.
Therefore $\varphi(a), \varphi(b)$ are also distinct.

Given $L$ in the axis of $\gamma$ in
$\hhs$, let $l = \Theta(L)$. 
The ideal points of $\oos(l)$ 
separate $a$ from $b$ in $\partial \oo$.
Then the source/sink property for the action
of $\gamma$ on $\partial \oo$ immediately
translates into a source/sink property for the action
of $\gamma$ on $\rr$ with source $\varphi(a)$
and sink $\varphi(b)$.
This finishes the proof of lemma \ref{soursi}.
\end{proof}

We now prove the first part of theorem \ref{conver}.

\begin{theorem}{}{}
Suppose that $\Phi$ does not have perfect fits and 
is not conjugate to a suspension Anosov flow.
Then $\pi_1(M)$ acts on $\rr$ as a convergence group.
\label{cogr}
\end{theorem}

\begin{proof}{}
Let $\gamma_i$ be a sequence of distinct elements of $\ga$.
Up to subsequence we can assume that either 

1) each $\gamma_i$
is associated to a singular closed orbit of $\Phi$;

2) each $\gamma_i$ is associated to a nonsingular closed
orbit of $\Phi$; 

3) each $\gamma_i$ is not associated to
a closed orbit of $\Phi$.

\noindent
Notice that 3) is equivalent to $\gamma_i$ having no fixed 
points in the orbit space $\oo$.
There is some similarity between cases 1) and 2) which will
be explored as we go along the proof.

\vskip .2in
\noindent
{\underline {Case 1}} $-$ Suppose the $\gamma_i$ are all
associated to singular orbits of the flow $\Phi$.

Let $\alpha_i$ be orbits of $\wwp$ with
$\gamma_i(\alpha_i) = \alpha_i$.
There are only finitely many singular orbits of $\Phi$,
so we may assume up to subsequence that all 
$\pi(\alpha_i)$ are the same.
We may also assume that
$\gamma_i$ are associated to (say) the
positive flow direction of $\alpha_i$, that is,
if $p_i$ in $\alpha_i$ then $\gamma_i(p_i) = \wwp_{t_i}(p_i)$
with $t_i$ bigger than zero.
Let $x_i = \Theta(\alpha_i)$ and
$l_i = \oos(x_i)$.

\vskip .1in
\noindent
{\underline {Case 1.a}} $-$
$(l_i)$ does not escape compact 
sets in $\oo$.

It could be that, up to subsequence, $l_i$ is constant.
This means that there is $\gamma$ in $\pi_1(M)$ so that 
$\gamma_i = \gamma^{n_i}$ and $|n_i|$ converging to
$\infinity$.
By the previous lemma
there is a source/sink pair for the sequence
$(\gamma_i)$.

Hence we may assume that 
up to subsequence all $l_i$ are
distinct and converge to a line leaf $l$ of $\oos$.
Up to subsequence assume the $l_i$ are nested and all
in a fixed sector of $l$.
Let $u, v$ be the ideal points of $l$.

\begin{figure}
\centeredepsfbox{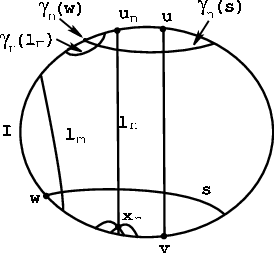}
\caption{The case of line leaves converging to a limit.}
\label{ca1}
\end{figure}

\vskip .1in
\noindent
{\underline {Claim 1}} $-$ There is an ideal point (say)
$v$ of $l$
so that all ideal points
of $l_i$ except for one converge to $v$.
The remaining ideal point of $l_i$ converges to $u$.

Otherwise up to subsequence there are at least 2 ideal
points $u^1_i, u^2_i$ of $l_i$ converging to $u$ and
likewise to $v$. Let $x_i$ be the singular point
of $l_i$. There is at least one 
{\underline {unstable}} prong of $\oou(x_i)$
with an ideal point in $\partial \oo$ between $u^1_i, u^2_i$
very near $u$ and similarly an unstable  prong of $\oou(x_i)$
with ideal
point very near $v$. Their union is a slice
$s_i$ of $\oou(x_i)$ with one ideal point near $u$ and
one ideal point near $v$. 
This slice is {\underline {not}} a line
leaf of $\oou(x_i)$ since there are 2 prongs
of $\oos(x_i)$ on both sides of this slice.
The sequence $(s_i)_{i \in {\bf N}}$
is nested and is bounded by $l$.
Hence it converges to a leaf $s$ of $\oou$.
By lemma \ref{duo} the ideal points
of $s_i$ 
converge to the ideal points of $s$ and
hence $s$ has ideal
points $u, v$. 
But $u$ is also an ideal point of the line leaf $l$ of $\oos$.
Since there are no perfect
fits, no 
leaves of $\oos, \oou$ share an ideal point.
This proves claim 1.

\vskip .1in
Since at least 2 ideal points of $\oos(x_i)$ converge
to $v$ (as $i \rightarrow \infty$) and ideal points
of $\oos(x_i), \oou(x_i)$ alternate in $\partial \oo$, then
at least one ideal point of $\oou(x_i)$ converges to $v$
as $i \rightarrow \infty$.
Suppose for a moment that not all endpoints of $\oou(x_i)$
converge to $v$. Then up to subsequence assume one
of the endpoints converges to $w$ distinct from $v$.
By the escape lemma (lemma \ref{duo}) up to subsequence
$(\oou(x_i))$ converges to a leaf $\delta$ of $\oou$
which has an ideal point $v$.
But $v$ is also an ideal point of line leaf $l$ of $\oos$,
contradiction to no perfect fits by lemma \ref{equiva}.
We conclude that all ideal points of $\oou(x_i)$
converge to $v$.

In order to finish the analysis of case 1.a it is enough to analyse
the following situation, which we state as a separate case as it will
be useful later on:

\vskip .1in
\noindent
{\underline {Case 1.b}} $-$ Suppose that $\oou(x_i)$ escapes compact 
sets in $\oo$, but $\oos(x_i)$ does not escape compact sets
in $\oo$.

Up to subsequence suppose that $\oos(x_i)$ converges to a line
leaf $l$ of $\oos$. Since $\oou(x_i)$ escapes compact sets 
it converges to an ideal point of $l$, which we denote by $v$
(again this follows from lemma \ref{duo}).
Let $u$ be the other ideal point of $l$.

Let $Z_i$ be the component of $\partial \oo - \partial \oou(x_i)$
which contains $u$.
In this case 
$(Z_i)$ converges 
to the set
$\partial \oo - \{ v \}$.
Let $u_i$ be the ideal point of $\oos(x_i)$ very close to $u$.
Suppose first up to subsequence that $\gamma_i(Z_i)$ is not equal
to $Z_i$ for all $i$.
Then $\gamma_i(Z_i)$ is an arbitrary small interval
very close to $v$. This shows that \ $\gamma_i | (\partial \oo - v)$ \
converges locally uniformly to $v$ and
so in $\rr$ it follows that
\ $\gamma_i | (\rr - \varphi(v))$ \ converges locally
uniformly to $\varphi(v)$.
So we assume from now on that $\gamma_i(Z_i) = Z_i$ for
all $i$ and hence $\gamma_i(u_i) = u_i$.
As the $\gamma_i$ are associated to 
positive direction of the flow
then the ideal points of $l_i =  \oos(x_i)$ are attracting for
the action of $\gamma_i$ in $\partial \oo$ (lemma \ref{soursi}).

\vskip .1in
\noindent
{\underline {Claim 2}} $-$
\ $\gamma_i | (\partial \oo - v)$ \ converges
locally uniformly to  $u$.

We already know that $\gamma_i(Z_i) = Z_i$ for all $i$.
As $v$ is an ideal point of a leaf of $\oos$ and $\Phi$ has
no perfect fits then $v$ has a neighborhood basis defined
by unstable leaves.
So it suffices to show that for a fixed unstable leaf
$s$ intersecting $l$, the
endpoints of $\gamma_i(s)$ converge to $u$.
Assume for simplicity that $s$ is nonsingular.

Notice first that it may be that the sectors of $l_i$
are not invariant under $\gamma_i$. 
A priori it may seem that this cannot happen because $\gamma_i(Z_i) = Z_i$.
But in fact 
this occurs when $\gamma_i$ acts in an orientation reversing
way on $\oo$ or equivalently on $\partial \oo$.
Then the other components of $\partial \oo - \partial \oou(x_i)$ 
are not $\gamma_i$ invariant (there are $\geq 2$ such 
other componets as $x_i$ is singular),
and the components of $\partial \oo - \oos(x_i)$ are also not invariant.

To analyse claim 2,
notice that $\gamma_i(s)$ intersects $l_i$. If one endpoint
of $\gamma_i(s)$ converges to $u$ (as $i \rightarrow \infty$),
then as seen above (using the escape lemma) the
other endpoint of $\gamma_i(s)$ also 
converges to $u$ and so $\gamma_i | (\partial \oo - \{ v \})$
converges locally uniformly to $u$ as desired.

\begin{figure}
\centeredepsfbox{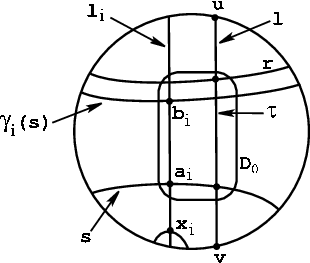}
\caption{
The case when $\gamma_i(s)$  converges to a leaf $r$.} 
\label{ego}
\end{figure}

The remaining case is up to subsequence $\gamma_i(s)$ converges 
to a leaf $r$
of $\oou$. 
Here $u$ cannot be in $\partial r$ and so $r$
intersects $l$.
Let $\tau$ be the segment of $l$ between
$s$ and $r$ and $D_0$ a neighborhood of it in
$\oo$. Let $D$ be the image of a smooth section 
$c_1: D_0 \rightarrow \mi$ of $\Theta$ restricted to $D_0$.
Recall the orbits $\alpha_i$ of $\wwp$
with $\gamma_i(\alpha_i) = \alpha_i$.
Let $\beta_i = \pi(\alpha_i)$, closed orbits of $\Phi$.
Then $\ws(\alpha_i) \cap D$ are segments
of bounded length. Let 

$$p_i \ = \ \ws(\alpha_i) \cap D \cap (s \times \rrrr),
\ \ \ a_i \ = \ \Theta(p_i), \ \ \ b_i \ = \ \Theta(\gamma_i(p_i)).$$

\noindent
In $D$ we have a segment $r_i$ of bounded length from
$\wwr(p_i)$ to a point in
$\gamma_i(\wwr(p_i))$. 
This is a segment in a {\underline {stable}}
leaf which contracts in positive flow direction.
Flow forward $p_i$ by time $t_i$
until it is distance $1$ from $\alpha_i$
along $\ws(\alpha_i)$. 
Notice that $p_i$ is far from $\alpha_i$ for $i$ big
since $x_i$ escapes compact sets in $\oo$ $-$ hence
$t_i >> 1$.
The segment $r_i$ flows to
a segment of arbitrary small length
under $\wwp_{t_i}$ since $r_i$ has bounded
length and $t_i$ is very big. This is 
a contradiction: the endpoints of $\wwp_{t_i}(r_i)$
both project in $M$ to the
same orbit in  $W^s(\beta_i)$ and the same local
sheet of the foliation $\ls$, but not the same local
flowline of $\Phi$.
Hence these endpoints cannot be
too close since the endpoint $\pi(\wwp_{t_i}(p_i))$ is
distance $1$ from $\beta_i$ in $W^s(\beta_i)$.
We conclude that this cannot happen.

It follows that $\gamma_i(s)$ cannot converge to a leaf
intersecting $l$ and so as seen before, $\gamma_i(s)$ converges
to $u$ in $\cd$ and the endpoints of $\gamma_i(s)$ also do.
This proves claim 2.

\vskip .1in
This completes the analysis of case 1.b, and hence also of case 1.a
that is,
 when the $l_i = \oos(x_i)$ do not 
escape compact sets in $\oo$. 
The same proof applies when $\oou(x_i)$ do not escape compact sets.

\vskip .1in
\noindent
{\underline {Case 1.c}} $-$ 
The sequences $\oos(x_i), \oou(x_i)$ escape compact
sets in $\oo$ and up to subsequence all
ideal points of $\oos(x_i), \oou(x_i)$ converge to the
 point $v$ of
$\partial \oo$.

We can assume that $v$ has a neighborhood basis defined
by (say) stable leaves.
Given a compact set $C$ in \ $\partial \oo - v$ \ let $s$ be a
nonsingular stable leaf with ideal points very close to $v$ and
separating $v$ from $C$ in $\cd$.
For $i$ big enough all the ideal points of
$\oos(x_i), \oou(x_i)$ are separated from $C$ by
$\partial s$.
Then $s$ is contained in a single interval of 
\ $\partial \oo - (\partial \oos(x_i) \cup \partial \oou(x_i))$ \ 
where $\gamma_i$ does not have fixed points.
If $\gamma_i$ leaves this interval invariant then since $\gamma_i(s)$
does not intersect $s$ transversely, then $\gamma_i(s)$
has both ideal points closer to $v$ than those of $s$
and so $\gamma_i(C)$ is very close to $v$ in $\cd$.
If $\gamma_i$ does not leave that interval invariant then as
seen above $\gamma_i(C)$ is also very close to $v$.
As $s$ is arbitrary this shows 
$\gamma_i(C) \rightarrow v$ uniformly. 
Therefore in $\rr$ it follows that
$\gamma_i | ({\cal R} - \varphi(v))$ converges
locally uniformly to $\varphi(v)$.

This finishes the
analysis of case 1: the $\gamma_i$ are associated to singular orbits.

\vskip .2in
\noindent
{\underline {Case 2}} $-$ $\gamma_i$ is associated to nonsingular
periodic orbits.

This is very similar to case 1 and we can use a lot of the
previous analysis.
We also use the following fact, which is a uniform statement
that orbits in leaves of $\Lambda^u$ are backwards asymptotic:

\vskip .1in
\noindent
{\underline {Fact}} $-$ Let $\Phi$ be a pseudo-Anosov
flow in $M^3$. For each $a_0 > 0$ and $\epsilon > 0$
 there is time $t_0 > 0$
so that if $p, z$ are in the same leaf of $\wlu$ and
there is a path $\delta$ in $\wu(p)$  from $p$ to $z$
with length bounded above by
$a_0$, then there is
a path from $\wwp_{t}(p)$ to $\wwr(z)$ in $\wu(p)$ of length
less than $\epsilon$ for all $t \leq -t_0$.

\vskip .08in
Equivalently the orbits $\wwr(p), \wwr(z)$ are $\epsilon$
close to each other backwards  of $\wwp_{-t_0}(p)$.
This is proved in pages 486-487 of \cite{Fe6}.
Notice it is not at all implied that
$\wwp_{-t_0}(p)$ and
$\wwp_{-t_0}(z)$ are $\epsilon$-close,
which may not be true since $p, z$ may
be out of phase.

\vskip .1in
\noindent
{\underline {Case 2.a}} $-$ Suppose that both $\oos(x_i)$ 
and $\oou(x_i)$ escape compact sets in $\oo$.

This is very similar to the singular situation. A proof exactly
as in case 1.c yields the result.

\vskip .1in
\noindent
{\underline {Case 2.b}} $-$ Suppose that exactly one
of $\oos(x_i)$ or $\oou(x_i)$ escapes compact sets.

Wlog assume that $\oou(x_i)$ escapes compact sets and $\oos(x_i)$
converges to a line leaf $l$ of $\oos$. 
Then a proof exactly as in Case 1.b yields the result.

\vskip .1in
\noindent
{\underline {Case 2.c}} $-$ Assume up to subsequence that $x_i$ converges to
$x$ in $\oo$.

If $x_i = x$ for infinitely many $i$ then lemma \ref{soursi}
finishes the proof. So we may assume up to
subsequence that $x_i$ are all nonsingular,
distinct from each other and all in the
same sectors of $\oos(x)$ and $\oou(x)$.
This did not occur in the previous case because the set
of singular points in $\oo$ is a discrete subset of $\oo$.
Let $e$ be the boundary of this sector of $\oou(x)$ $-$ a line leaf of 
$\oou(x)$.
Assume wlog that up to subsequence $\gamma_i$ is associated 
to positive flow direction in $\alpha_i$. 
Hence $\partial \oos(x_i)$ is the attracting fixed point
set for $\gamma_i$ and
$\partial \oou(x_i)$ is the repelling fixed point set
for the action of $\gamma_i$ on $\partial \oo$.


We will show that $\partial \oou(x),
\partial \oos(x)$ forms a source/sink set for the sequence
$\gamma_i$ acting on $\partial \oo$. Then
\ $a = \varphi(\partial \oou(x)),
\ b = \varphi(\partial \oos(x))$ \ forms a source/sink pair
for the sequence $\gamma_i$ acting on $\rr$.
For simplicity assume that $\gamma_i$ preserves the components
of $\oos(x_i) - x_i, \ \oou(x_i) - x_i$. A similar proof works in
the general case.

Let $\alpha_i = \{ x_i \} \times \rrrr$ and
$\pi(\alpha_i)$ closed orbits of $\Phi$.
Assume all $\pi(\alpha_i)$ are distinct.
Let $v$  be a point in $\partial e$ (which is 
a subset of $\partial \oou(x)$).
For any small neighborhood $A$ of $v$ in $\cd$
let $l$ nonsingular stable leaf intersecting $\oou(x)$ 
and contained in $A$.
As $\oou(x_i)$ converges to $e$
(a line leaf in) $\oou(x)$ then for $i$ big
enough $\oou(x_i)$ intersects $l$ and has an ideal
point $v_i$ near $v$.
Since $v_i$ is
a repelling fixed
point for $\gamma_i$ then $\gamma_i(l)$ is closer
to $\oos(x_i)$ than $l$ is. 
Here $\oos(x_i)$ is close to $\oos(x)$ as
well.
Let $L_i$ = $\gamma_i(l) \times \rrrr$, a leaf of $\wls$.

The fact that is going to be used here is that the
lengths of the periodic orbits $\pi(\alpha_i)$ converge
to infinity, which occurs because they are all
distinct orbits.
Draw a disk $D$ transverse to $\wwp$
containing segments $r_i$ in $\wu(\alpha_i)$ from $p_i$ in
$\alpha_i$ to

$$z_i \ = \ (l \times \rrrr) \cap \wu(\alpha_i) \cap D$$

\noindent
and $r_i$ transverse to $\wwp$ in $\wu(\alpha_i)$.
We can assume the $r_i$ converges to $r$, which is
a segment in
$\wwu(p)$  (here $p = \{ x \} \times \rrrr$)
and so  the $r_i$
have diameter uniformly bounded above.
Consider $\gamma_i(r_i)$ which are segments of
diameter bounded above, connecting $\gamma_i(p_i)$ to
$\gamma_i(z_i)$.
Notice that $\gamma_i(z_i)$ is in $L_i$.
Choose 

$$t_i \in \rrrr \ \ {\rm with} \ \ 
\gamma_i(p_i) \ = \ \wwp_{t_i}(p_i).
\ \ \ {\rm Then} \ \ t_i \rightarrow \infty
\ {\rm and} \ \ p_i \ = \ \wwp_{-t_i}(\gamma_i(p_i))$$

\noindent
By the fact above there are segments
from $p_i$ to $\wwr(\gamma_i(z_i))$ in $\wu(\alpha_i)$ with
diameter converging to $0$ as $i \rightarrow \infty$.
As the $p_i$ are converging to the point $p$
in $\{ x \} \times \rrrr$,
this shows that  $\gamma_i(l)$ is converging to 
(a line leaf of) $\oos(x)$.


This shows
that $\partial \oou(x)$ is the repelling fixed point
set for $\gamma_i$ and
$\partial \oos(x)$ is the attracting set. This
finishes the analysis of case 2.

\vskip .2in
\noindent
{\underline {Case 3}} $-$ All the $\gamma_i$ act freely on $\oo$.

This case is extremely long.

By lemmas \ref{ideal} and \ref{soursi} each $\gamma_i$ acts
on $\partial \oo$ with only two distinct fixed
points $v_i, u_i$ forming a source/sink pair, that is,
hyperbolic dynamics in $\partial \oo$.
Assume up to subsequence that $u_i$ converges to $u$ 
and $v_i$ converges to $v$ in $\partial \oo$.
It may be that $u$ is equal to $v$.
Ideally we would like to show that 
\ $\gamma_i | (\partial \oo - v)$ \ converges locally 
uniformly to $u$. Very surprisingly this is not
true in general, see the counterexample after the end of the
proof.

We first consider the situation that $u = v$.
This is dealt with exactly as in case 1.c.

Hence from now on suppose that $u \not = v$.
Assume wlog that $v$ is not an ideal point of a leaf of $\oos$ and hence
by lemma \ref{options}, \ $v$ has a neighborhood system
defined by stable leaves. Let $l$ be
a non singular stable leaf with ideal points near $v$,
separating it from $u$.
This uses the fact that $u \not = v$.
If some subsequence of $(\gamma_i(l))$ escapes compact
sets in $\oo$, then by the escape lemma 
(lemma \ref{duo} part (iii)), the ideal points
of $\gamma_i(l)$ have to be very near each other. Then these ideal points
have to be very near $u_i$ and hence very near $u$.
If this happens for $l$ arbitrarily near $u$, then 
this implies the convergence property: compact sets of $\partial \oo - \{ v \}$
converge to $u$ under $\gamma_i$.
Hence by way of contradiction assume for the remainder of case 3:

\vskip .1in
\noindent
{\underline {Running hypothesis for the remainder of case 3}} $-$
Up to subsequence
suppose that there is $l^c$ with ideal points
very near $v$ and separating it from $u$, so 
that $\gamma_i(l^c)$ converges to a line leaf $l^d$ of
some leaf of $\oos$.
\vskip .1in

There are 2 possibilities.

\begin{figure}
\centeredepsfbox{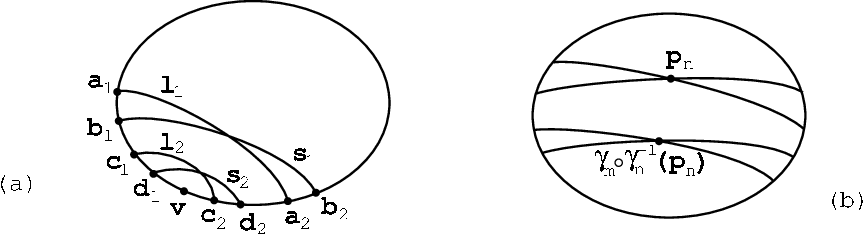}
\caption{a. Set up in $\oo$, b. Producing fixed points.}
\label{case3}
\end{figure}

\vskip .18in
\noindent
{\underline {Case 3.1}} $-$
The point $v$ is not an ideal point of a leaf of $\oou$.

Then there is a neighborhood system of $v$ defined by
unstable leaves as well. 
For a stable leaf $l$ as above let
$\partial l = \{ a _1, a_2 \}$, where we suppress the dependence
on $l$ for notational simplicity..
Consider the collection of unstable
leaves $\{ s \in \oou \ | \ s \cap l \not = \emptyset \}$. 

We claim that if
$l$ is close to $v$ then so are all the possible $s$.
Otherwise vary $l$ and take 
limits of $l$ approaching $v$ and also take limits of such
$s$ with one ideal point not close to $v$, then using the escape lemma
one produces an unstable leaf with
ideal point $v$, contrary to assumption.

In the same way
if $s \cap l$ is near $a_1$ in $\cd$ then $s$ is near
$a_1$ in $\cd$ and has all ideal points  near $a_1$.
Otherwise consider a sequence $s_n \in \oou$ with $s_n \cap l$
converging to $a_1$. If $s_n$ does not escape in $\oo$, then
the escape lemma produces an unstable leaf with ideal point
$a_1$, contrary to hypothesis in this case. Since
$s_n$ escapes compact sets and has $s_n \cap l$ converging
to $a_1$, lemma \ref{duo} again implies that the ideal
points of $s_n$ also converge to $a_1$.
Similarly if $s \cap l$ is near $a_2$ in $\cd$ then $s$ is near
$a_2$ in $\cd$.
It follows that there is a unique unstable leaf $s'$ intersecting
$l$
so that $s'$ has a singularity in $W$ and
has at least 2 prongs contained in $W$ and
enclosing $v$. Enclosing $v$ means that if $b'_1, b'_2$ are
the ideal points of these 2 prongs then
$a_1, b'_1, v, b'_2, a_2$ are all distinct and circularly
ordered in $\partial \oo$ 
(under some circular order in $\partial \oo$).
There is then one prong of $s'$ exiting $W$ so that
together with a prong inside $W$ it
defines a small neighborhood of $v$.  The union of these
two prongs is a slice $s_1$ in $s'$. 
Let $\partial s_1 = \{ b_1, b_2 \}$ with $b_1$ an
ideal point of $W$.
Let $l_1 = l$.
This was the first step of the process, which is going
to be done twice.
We know that $\gamma_i(l_1)$ does not limit to $u$ and
we can assume up to subsequence that $\gamma_i(l_1)$ converges
to $l_0$ a stable leaf with no limit point in $u$.

Now redo the process above to obtain a leaf $l_2$ of
$\oos$ and a slice $s_2$ of $\oou$ which are closer
to $v$.
Let $\partial l_2 = \{ c_1, c_2 \}$ and
$\partial s_2 = \{ d_1, d_2 \}$.
By doing this procedure 3 or 4 times,
we can arrange the construction so that for instance
$a_1, b_1, c_1, d_1, v, c_2, d_2, a_2, b_2$
are all distinct and circularly ordered in $\partial \oo$,
see fig. \ref{case3}, a.

The $\gamma_i(l_2), \gamma_i(s_2)$ do not escape
to $u$, because they are bounded by $l_d$.
Let $j = 1, 2$.
We may assume that the sequence $(\gamma_i(l_j))$ 
is nested and converges to
$l'_j$ and likewise $(\gamma_i(s_j))$ is nested
and converges to $s'_j$,
as $i \rightarrow \infty$ for $j = 1, 2$. 
Because of the set up of the ideal points as above then
$l'_1$ has no common ideal point with $l'_2$.
If for instance $\lim \gamma_i(a_1) = \lim \gamma_i(c_1)$ then
it is also equal to $\lim \gamma_i(b_1)$ and
one produces one unstable leaf $s'_1$ sharing an ideal
point with a stable leaf $l'_1$ $-$ disallowed by no perfect fits.
It follows that all four limits of ideal points are distinct.
Fix $n$ very big and let $m > > n$.
Since $\gamma_m(l_1), \gamma_n(l_1)$ are both very near $l'_1$
and $\gamma_m(s_1), \gamma_n(s_1)$ are very near $s'_1$
then 

$$\gamma_m(l_1 \cap s_1) \ \ \ {\rm is \ very \ near} \ \ \
\gamma_n(l_1 \cap s_1) = p_n$$

\noindent
Or $\gamma_m \circ \gamma^{-1}_n (p_n)$ is
very near $p_n$, see fig. \ref{case3}, b.
If $l'_1 \cap s'_1$ is singular assume up to subsequence
that all $\gamma_i(l_1 \cap s_1)$ are in the intersection
of closures of sectors of 
$\oos(l'_1 \cap s'_1)$  and $\oou(l'_1 \cap s'_1)$.
With these conditions and the fact that
$\gamma_m, \gamma_n$ are distinct, then
the shadow lemma for pseudo-Anosov flows \cite{Han,Man}
implies that $\gamma_m \circ \gamma^{-1}_n$
has a fixed point very near $p_n$.
Similarly there is a fixed point of $\gamma_m \circ \gamma^{-1}_n$
near $\gamma_n(l_2 \cap s_2)$. 
Since $l'_1 \cap s'_1, l'_2 \cap s'_2$ are different, then
for $n, m$ sufficiently big these two fixed points
are different.
But then $\gamma_m \circ \gamma^{-1}_n$ would
have two distinct fixed points  in $\oo$ $-$ 
which is disallowed by the no perfect fits condition.
This cannot happen. Therefore $\gamma_i(l)$ converges
to $u$ for any $l$ close enough to $v$
and this finishes the analysis of case 3.1.

\vskip .18in
\noindent
{\underline {Case 3.2}} $-$
Suppose that $v$ is an ideal point of a leaf $s$ of
$\oou$.

The proof of this subcase is very long.
In this case we do not necessarily obtain
that $\gamma_i | (\partial \oo - v)$ \ converges locally
uniformly to $u$.
Suppose $l$ is nonsingular, intersects $s$ and $W \cap s$ has no
singular points.
As in case 3.1 we only have to deal with the case that 
$\gamma_i(l)$ does not escape compact sets in $\oo$.

\vskip .05in
From now on in this case fix this leaf  $l$ of $\oos$.
\vskip .05in

Assume that $\gamma_i(l)$ converges to 
a line leaf $l^*$  of  a leaf
$l_0$ of $\oos$.
Let $l'$ be any stable leaf intersecting $s$ and
closer to $v$ than $l$ is.

The first situation is 
that up to subsequence $\gamma_i(l')$
converges 
to $l'_0$ different from $l_0$.
Then $l_0, l'_0$ do not share an ideal point $-$ because of
the no perfect fits hypothesis.
Since $\gamma_i(s)$ intersects $\gamma_i(l), \gamma_i(l')$
and $\gamma_i(l), \gamma_i(l')$ converge to $l_0, l'_0$
not sharing an ideal point then
$\gamma_i(s)$ cannot escape in $\oo$.
This follows directly from the escape lemma.

Hence assume $\gamma_i(s)$ converges to a leaf $s_1$ of $\oou$.
Notice that
$s_1$ intersects $l_0$ and $l'_0$ for
otherwise, by the escape lemma again,
$s_1$ will share ideal point 
with at least one of $l_0, l'_0$, again
disallowed by the no perfect
fits condition.
Therefore $\gamma_i(l \cap s)$ converges to $l_0 \cap s_1$
and $\gamma_i(l' \cap s)$ converges to $l'_0 \cap s_1$.
As seen before, if $n, m$ are big enough this produces
2 distinct fixed points of $\gamma_m \circ \gamma^{-1}_n$
$-$ one near $l_0 \cap s_1$ and one near $l'_0 \cap s_1$.
This is disallowed.

\vskip .1in
We conclude that for any $l'$ stable leaf intersecting $s$ and
separating $v$ from $l$, 
the sequence $\gamma_i(l')$ also converges to $l^*$.
Let 

$$w, w' \ \ \ {\rm
be \ the \ ideal \ points \ of} \ \ \  l^*.$$

\noindent
Let $z, z'$ be the ideal points of $l$.
Let $I, I'$ be the 
disjoint half open intervals of $\partial \oo$ with one ideal
point in $z, z'$ and the other in $v$, that is,
$z \in I$ but $v$ is not in $I$ (for some orientation 
of $\oo$ then $I = [z,v), I' = (v,z']$).
Assume wlog that $\gamma_i(z)$ converges to $w$. The arguments
above show that $\gamma_i(I)$ converges locally
uniformly  to $w$ and
$\gamma_i(I')$ converges locally uniformly to $w'$.

The strategy to prove case 3.2 is as follows:
Using the no perfect fits condition we will incrementally
upgrade the property above
to show that \ $\gamma_i | (\partial \oo - \partial s)$ \
converges locally uniformly to $\partial l^*$ $-$ this last
one is a set, not a single point.
This means that for any $C$ compact contained in 
$\partial \oo - \partial s$, then for $i$ big enough
$\gamma_i(C)$ is contained in a small neighborhood of
$\partial l^*$.
Notice that $s$ may be singular so the set \  $\partial \oo - \partial s$ \
may have more than 2 components.

Recall that in case 3.2 the leaf $l$ of $\oos$ is fixed.
Consider an arbitrary
unstable leaf $s'$ intersecting $l$,
with $s' \not = s$.
Then $s'$ has at least one
ideal point in either $I$ or $I'$. If $s'$ 
has an ideal point $t$ in $I$ then 
$\gamma_i(t)$ converges to $w$.
Since no unstable leaf has ideal
point $w$ it follows from the escape lemma
that $\gamma_i(s')$ converges to $w$ in $\cd$.
Let now $J$ ($J'$) be the component of $(\partial \oo - \partial s)$
containing $z$ ($z'$)
(hence $I \subset J, I' \subset J'$). 
The
above arguments imply that $\gamma_i(J)$ converges 
locally uniformly to $w$ and 
$\gamma_i(J')$ converges locally uniformly to $w'$. 
To prove this use the fact that for any 
$v \in J - \overline I$ there is $s'$ unstable leaf
with $s' \cap l \not = \emptyset$ and $\partial s'$ 
separating $\partial s$ from $v$ in $\partial \oo$.
This last statement follows from the escape lemma
and the fact that there are no leaves in $\oou$ non separated
from $s$.

If $s$ is nonsingular we are done.
This is because if $v, t$ are the ideal points of $s$,
then $\varphi(v) = \varphi(t) = y$ and
$\varphi^{-1}(y) = \{ v, t \}$.
For any compact set $C$ in \ $\rr - y$ \
there is a compact set $V$ in $\partial \oo - \{ v, t \}$ with
$C \subset \varphi(V)$,
since $\varphi^{-1}(y) = \{ v, t \}$.
Hence $V$ is contained in
the union of 2 compact intervals $V_1, V_2$
in $\partial \oo - \{ v, t \}$ so up to reordering $V_1 \subset J$
and $V_2 \subset J'$.
Hence

$$\gamma_i \ | \ V_1 \ \ {\rm converges \ to} \ \  w
\ \ \ {\rm and} \ \ \ 
\gamma_i \ | \ V_2 \ \ {\rm converges \   to} \ \ w'.$$

\noindent
Notice $\varphi(w) = \varphi(w')$ and let this be $x$.
This shows that in $\rr$, $\gamma_i | C$ converges uniformly
to $x$.
Hence $y, x$ is the source/sink pair for $\gamma_i$.
With more analysis one can show that 
$u$ is not an ideal point of $s$ and $u$ is
an ideal point of $l_0$.
We do not provide the
arguments as we will not use that.

\begin{figure}
\centeredepsfbox{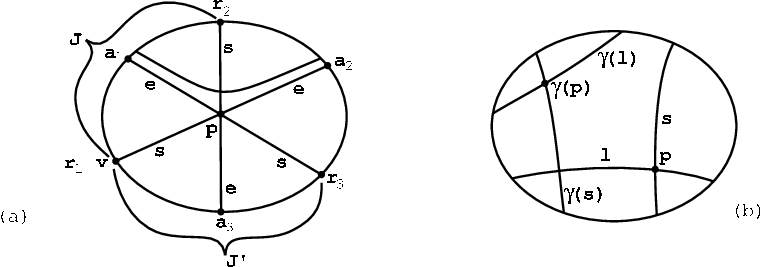}
\caption{
a. Trapping the orbits 
in the singular case ($n = 3$), b. An interesting
 counterexample.}
\label{singt}
\end{figure}

\vskip .1in
To finish the analysis of case 3.2, we
suppose from now on that $s$ is singular with $n$ prongs.
Let $r_1 = v$, \
let $r_2$ be the other endpoint of $J$  $-$
this is also an
ideal point of $s$  and let
$r_n$ be the similar endpoint of $J'$.
Complete the ideal points of $s$ circularly to $r_1, ..., r_n$.
Let $p$ be the singular point in $s$, see
fig. \ref{singt}, a.
Let $e$ be the stable leaf through $p$.
Then $e$ has a prong with ideal point
$a_1$ in $J$ and one with ideal point
$a_n$ in $J'$.
Order the other endpoints of $e$ as
$a_1, ..., a_n$.
Let $e^*$ be the line leaf of $e$ with
ideal points $a_1, a_n$.
We proved above that $\gamma_i(a_1)$ converges to $w$ and
$\gamma_i(a_n)$ converges to $w'$.
There are two options: 1) $\gamma_i(p)$ does not escape compact sets in $\oo$;
\ 2) $\gamma_i(p)$ escapes compact sets in $\oo$.

\vskip .1in
\noindent
{\underline {Option 1}} $-$
Suppose that $\gamma_i(p)$ does not
escape compact sets in $\oo$.

Up to subsequence  $\gamma_i(p)$ converges in $\oo$
so assume that
$\gamma_i(p) = p_0$ for $i \geq i_0$ (using the fact that 
$p$ is singular). 
Let
$f$ be the generator of the isotropy group of
$p_0$ fixing also 
$\gamma_{i_0}(a_1) = w, \ \gamma_{i_0}(a_n) = w'$
and $f$ associated to the forward direction in
the orbit $\{ p_0 \} \times \rrrr$.
Since $\gamma_i(a_1), \gamma_i(a_n)$ converge 
to $w, w'$, there is $i_0$ so that for $i > i_0$:

$$\gamma_i  \ = \ f^{m_i} \circ \gamma_{i_0}, \ \ \ 
m_i \in \ {\rm {\bf Z}}$$ 

\noindent
Here $w$ is an attracting point, so $m_i$ converges to
$+\infty$. 
Lemma \ref{soursi} now implies 
that for any compact  set $C$ in
\ $\partial \oo - \partial s$ \ then $\gamma_i(C)$ is in a small
neighborhood of $\partial \oos(p_0)$. All points
of $\partial s$ are identified under $\varphi$ 
and similarly for $\partial \oos(p_0)$. Let $y = \varphi(\partial s)$,
$x = \varphi(\partial \oos(p_0))$.
Then in $\rr$, \
$\gamma_i  | (\rr - y)$ \ converges 
locally uniformly to $x$.
This finishes the argument for option 1).


\vskip .1in
\noindent
{\underline {Option 2}} $-$ Suppose that 
$\gamma_i(p)$ escapes compact sets in $\oo$.

Since
$\gamma_i(e^*)$ converges to $l^*$ and 
$\gamma_i(p)$ is in $\gamma_i(e^*)$, the escape lemma implies that 
$\gamma_i(p)$ converges to either $w$ or $w'$.
Suppose without loss of generality that $\gamma_i(p)$ converges
to $w$. Then $\gamma_i(a_n)$ converges to $w'$ and
all the ideal points $a_1, ..., a_{n-1}$ converge
to $w$ under $\gamma_i$. 
Here is the justification of this statement:
If for some $j$ in $2,...,n-1$, \ 
$\gamma_i(a_j) \rightarrow w'$, then 
also $\gamma_i(a_{n-1}) \rightarrow w'$. Hence $\gamma_i(r_n)
\rightarrow w'$. But since $\gamma_i(p) \rightarrow w$, then
that unstable prong of $\oou(p)$ converges to an unstable leaf
with one ideal point in $w$ and another in $w'$. This is disallowed
under no perfect fits (in fact this cannot happen in general, but
we will not need that).
Therefore $\gamma_i(a_n) \rightarrow w'$ and 
$\gamma_i(a_j) \rightarrow w$ for $j = 1, ..., n-1$.

Let $J_2$ be the interval of $\partial \oo$ bounded by
$v$ ($= r_1$) 
and $r_n$ and so that $J_2$ is 
disjoint from $J'$ $-$ that is, 
$J_2 = \partial \oo - \overline {J'}$. If $A$ is the
region of $J_2$ between $a_{n-1}$ and $r_n$ we claim that
$\gamma_i(J_2 - A)$ converges to $w$.
Here $(\gamma_i(a_j))$ converges to $w$,
$1 \leq j \leq n-1$.
The nonsingular unstable leaves $s'$ intersecting
$\oos(p)$ in the prong with ideal point $a_{n-1}$
have one ideal point in $A$ and another
ideal point $y$ in $J_2 - A$. 
Since $\gamma_i(a_{n-1}) \rightarrow w$, then $\gamma_i(y) \rightarrow w$.
This implies that the ideal points of $s'$ 
both have to converge
to $w$ under $\gamma_i$.
Since $\gamma_i$ is a homeomorphism of $\partial \oo$ it now follows that
$\gamma_n(J_2)$ converges locally uniformly
to $w$. As $\gamma_i(J')$ converges locally uniformly to $w'$ then 

$$\gamma_i \ | \ (\partial \oo - \partial s)
\ \ \ {\rm converges \ locally \ uniformly \ to} \ \ \ 
\{ w, w' \}.$$

\noindent
Notice that $\partial \oo$ is the disjoint union
of $J_2, J', r_1, r_n$.
 If $y = \varphi(\partial s)$ and
$x = \varphi(w)$ then $y, x$ is the source/sink pair
for a subsequence $\gamma_i$ acting on $\rr$.
This finishes case 3.2.


This shows that $\pi_1(M)$ acts as a convergence
group on $\rr$ and finishes the proof of theorem 
\ref{cogr}.
\end{proof}

\noindent
{\bf {Remark}} $-$ We construct an example 
as in case 3.2 where the sequence 
$(\gamma_i)_{i \in {\bf  N}}$ does not have source/sink
pair the points $v, u$ for the action
on $\partial \oo$ as naively expected in case 3.2.
In fact the source is a collection of points and so is the sink.
We start with $\Phi$ a pseudo-Anosov flow without
perfect fits and not conjugate to suspension Anosov.
For simplicity assume that everything is orientable.
In addition assume that $\Phi$ is transitive.
The tricky thing
is to get $\gamma_i$ to act freely on $\oo$.
Let \ $s = \oou(p), l = \oos(p)$ \ where $p$ is periodic,
nonsingular.
Let $\gamma$ in $\pi_1(M)$ with
$\gamma(s)$ intersecting $l$ and so that $\gamma(l)$ does
not intersect $s$, see fig. \ref{singt}, b.
Since $\Phi$ is transitive it is always possible to find such
$\gamma$ unless there is a product region 
in that quarter of $p$ $-$ but then $\Phi$ would be conjugate
to a suspension Anosov flow, contrary to assumption.
Let $f$ be the generator of the isotropy group of $\gamma(p)$
leaving invariant all points in 
$\partial \oos(\gamma(p)), \ \partial \oou(\gamma(p))$ and
associated to the positive direction of the flow line.
Let $\gamma_i = f^i \circ \gamma$.
Then $\{ \gamma_i \}_{i \in {\bf N}}$ are all distinct. 

Suppose that  for some $j$, the 
 $\gamma_j$ has a fixed point. Fix $j$ and 
let $h = \gamma_j$. Notice that $l, h(l)$ ($= \gamma(l)$)  both intersect
a common unstable leaf $\gamma(s)$; \
also $s, h(s)$ ($= \gamma(s)$)
intersect the stable leaf $l$ and
$s, h(l)$ do not intersect.
If $h^m(l)$ converges to $r$ as $m \rightarrow \infty$
thee $h(r) = r$. This is because the leaf space $\hhs$ of $\wls$
is Hausdorff.
Hence $h$ has a fixed point $q$ in $r$.
Then 
$h(\oou(q))  =  \oou(q)$
%
and $\oou(q)$ intersects $h^m(l)$ for $m$ big enough
and hence for all $m$.
But since $h$ contracts $h^m(l)$ towards $\oos(q)$ then
it expands unstable leaves away.
In particular $s$ cannot intersect $r$.
However by construction $h$ moves $s$ and $l$ in the
same direction and hence $\gamma_j(s)$ is closer to $\oou(q)$
than $s$ is.
It follows that $h^m(s)$ converges to a leaf $t$ and $t$ does
not intersect $r$. Hence $h(r) = r$, $h(t) = t$
and $r \cap t = \emptyset$.
This produces two fixed points of $h$ in $\oo$.
Hence theorem \ref{chain} implies that there are perfect fits,
contrary to assumption.
This contradiction shows that $h$ does not have any fixed point in
the component of $\oo - l$ containing $h(l)$.
Now consider $h^{-1}$: \ $h^{-1}(l) = \gamma^{-1}(l)$ does not
intersect $l$ and $h^{-1}(s) = \gamma^{-1}(s)$ intersects
$l$. So the same argument as above shows that $h$ does not
have a fixed point in the component of $\oo - l$ containing
$h^{-1}(l)$. Hence $h$ does not have fixed points in $\oo$ and
acts freely.

It follows that each $\gamma_i$ acts freely on $\oo$ and
has 2 fixed points $v_i, u_i$ in $\partial \oo$.
In addition as $i \rightarrow \infty$, \   $u_i$ converges
to an ideal point of $\gamma(l)$ and $v_i$ converges
to an ideal point of $s$ $-$ the one separated 
from $\gamma(l)$ by $l$. So this is exactly
the situation in case 3.2 of theorem \ref{cogr}.
Notice also
that $\gamma_i(s) = \gamma(s)$ and $\gamma_i(l) = \gamma(l)$ so
the collection
$\{ \gamma_i \}_{i \in {\bf N}}$ does not act properly 
discontinuously in $\oo$.
Here $\gamma_i(\partial s) = \gamma(\partial s)$ and $\gamma_i(\partial l) 
= \partial \gamma(l)$, so there are not two points
in $\partial \oo$ forming a source/sink pair
for the action of $(\gamma_i)$ on $\partial \oo$.
Still $\gamma_i | (\partial \oo - \partial s)$
converges locally uniformly to $\partial \gamma(l)$.

\vskip .1in
The next goal is to show that 
every point in $\rr$ is a conical limit point.

\begin{theorem}{}{}  
Let $\Phi$ a pseudo-Anosov flow without perfect fits
and not conjugate to a suspension Anosov flow.
Let $\rr$ be the associated sphere quotient of 
$\partial \oo$.
Then 
every point in $\rr$ is a conical
limit point for the action of $\pi_1(M)$ on $\rr$.
Hence $\pi_1(M)$ acts as a uniform convergence group
on $\rr$.
\label{unifco}
\end{theorem}

\begin{proof}{}
The last statement follows from the first because
theorem \ref{coni} implies that the
action of $\pi_1(M)$ on the space of distinct 
triples of $\rr$
is cocompact. 

We show that any $x$ in $\rr$ is a conical limit point
for the action of $\pi_1(M)$.
There are 3 cases:

\vskip .1in
\noindent
{\underline {Case 1}} $-$ $x = \varphi(z)$ where $z$ is
the ideal point of $l$ of $\oos$ or $\oou$ and there is
$\gamma \not =$ id in $\pi_1(M)$ with $\gamma(l) = l$.

Since all ideal points of $l$ are taken to $x$ under
$\varphi$ and $\gamma$ permutes the ideal
points of $l$, it follows that $\gamma(x) = x$.
Assume that $x$ is the repelling fixed point
of $\gamma$ $-$ up to taking an inverse if necessary.
Let $\gamma_i = \gamma^i, i \geq 0$.
Then $\gamma_i(x) = x$ so $\gamma_i(x)$ converges 
to $x$.
Let $c$ be the other fixed point of $\gamma$ in $\rr$.
For any $y$ distinct from $x$ in $\rr$ it
follows from lemma \ref{soursi},
that $\gamma_i(y) = (\gamma^i)(y)$ converges
to $c$. Hence $x$ is a conical limit point.

\vskip .1in
\noindent
{\underline {Case 2}} $-$
$x = \varphi(z)$ where $z$ is an ideal point of $l$ of
$\oos$ or $\oou$ and $l$ is not invariant under
any $\gamma$ of $\pi_1(M)$.

Suppose without loss of generality
that $l$ is an unstable leaf.
Let $L = l \times \rrrr$ a leaf of $\wlu$.
Here $\pi(L)$ does not have a periodic orbit of $\Phi$.
Let $\alpha$ be an orbit of $\wwp$ in $L$.
We look at the asymptotic behavior of $\pi(\alpha)$ in
the negative direction (all orbits
in $L$ are backward asymptotic, so this argument
is independent of the orbit $\alpha$ in $L$).
If $\pi(\alpha)$ limits only in a singular orbit then $\pi(\alpha)$
must be in the unstable leaf of a singular orbit,
contrary to assumption.

For each $i$ choose $p_i$ in $\alpha$ with
$(p_i)$ escaping in the
negative direction and $(\pi(p_i))$ converging to 
a nonsingular point $\mu$  in $M$.
By discarding a number of initial terms, we can assume that
all $\pi(p_i)$ are in a neighborhood $V$ of $\mu$ to which
the shadow lemma can be applied.
There are $\gamma_i$ in $\pi_1(M)$ with $\gamma_i(p_i)$
in $V$.
By the shadow lemma the $\gamma_i$ correspond to closed
orbits $\beta_i$ of the flow $\Phi$.
In particular we assume that $V$ is sufficiently small, so
that there is still a small neighborhood $U$ of $\mu$
with $\overline V \subset U$ and
there are lifts $\widetilde \beta_i$ of $\beta_i$ with
points in $U$.
We assume that $\overline U$ does not intersect any
singular orbit.
It follows that no $\beta_i$ is a singular orbit.
Since $\widetilde \beta_i, \gamma_i(\alpha)$ have points near
$p_1$, we may also
assume up to subsequence that both sequences converge.
Let $\tau$ be the limit of 
$(\widetilde \beta_i)$.
Hence $\tau$ is also not a singular orbit.
Notice that a priori there is no relation between $\mu$ and
$\tau$ except that $\mu$ is near $\tau$.
Let also $\delta$ be the limit of 
 $(\gamma_i(\alpha))$.
Notice that $\pi(\delta)$ has a point in $\overline V$. 

Each $\gamma_i$ takes $p_i$ to a point very close to $p_1$
and the $p_i$ escape in $\mi$ with $i$, so 
up to subsequence we can assume that the $\gamma_i$ are all distinct.
Hence the length of $\beta_i$ goes to infinity
(the $\beta_i$ does not have to be an
indivisible closed orbit).
Let $q_i = \Theta(\widetilde \beta_i)$, 
so $(q_i)$ converges to $q_0 = \Theta(\tau)$.
Let

$$\partial \oou(q_0) = \{ s, s' \}, \ \ 
\partial \oos(q_0) = \{ t, t' \}, \ \ 
\partial \oou(q_i) = \{ s_i, s'_i \}, \ \ 
\partial \oos(q_i) = \{ t_i, t'_i \}$$

\noindent
Since the points $p_i$ are flow backwards of $p_1$ in $\alpha$
and $p_i$ is sent near $p_1$ by $\gamma_i$, then
$\gamma_i$ corresponds to the flow lines $\beta_i$ being
traversed in the forward direction.
By lemma \ref{soursi},
$\{ s_i, s'_i \}$ is the repelling
set for the action of $\gamma_i$ on $\partial \oo$ and
$\{ t_i, t'_i \}$ is the attracting set.

Here $\oos(q_i)$ intersects $l = \oou(\alpha)$ and $\oou(\gamma_i(\alpha))$
for every $i$.
As described above $\gamma_i(l)$ converges to
the unstable leaf $r := \oou(\Theta(\delta))$.
Since $\oou(\gamma_i(\alpha))$ converges and the length
of $\beta_i$ goes to infinity, then the arguments of 
case 2.c of the proof of theorem \ref{cogr} show that the only
possibility is that $\oou(q_i)$ converges to $l =  \oou(\alpha)$ $-$ otherwise
$l$ would be pushed farther and farther away from $\oou(q_i)$.
This shows that $\tau$ is in $L$ and
$s, s'$ are the ideal
points of $l$.
Up to renaming the ideal points of $l$, \ $z = s$.
Up to another subsequence assume
that 

$$s_i \rightarrow s, \ \ \ s'_i \rightarrow s', \ \ \ 
t_i \rightarrow t, \ \ \ t'_i \rightarrow t'$$

\noindent
Again the arguments in case 2.c of theorem \ref{cogr}
show that in $\partial \oo$, we have  \
$\gamma_i \ | \ (\partial \oo - \{ s, s' \})$ \ converges
locally uniformly to the set $\{ t, t' \}$. 
Also $\gamma_i(z)$ converges to a point $d$  in 
$\partial r$.
As $z = s$ then
$x = \varphi(s)$ and  we
have in $\rr$ that \
$\gamma_i \ | \ (\rr - \{ x \})$ \ converges locally uniformly
to $\varphi(t)$. 
Since $d$ is a unstable ideal point and $t$ is a stable 
ideal point, it follows that $\varphi(t) \not = \varphi(d)$.
Summing it all up:

$$\gamma_i(x) \ \rightarrow \ \varphi(d) \ \ \ {\rm and} \ \ \ 
\gamma_i(y) \ \rightarrow \ \varphi(t) \ \ \ {\rm for \ any}
\ \ y \ \in \ {\cal R} - \{ x \}$$

\noindent
This shows that 
$x$ is a conical limit
point.

\vskip .15in
\noindent
{\bf {Remark}} $-$ Obviously it is crucial in this proof that 
$z$ is an ideal point of a leaf of $\oos$ or $\oou$. Since
$z$ is an unstable ideal point and we want to push points
away from the unstable ideal point,
then in the proof above
we use $\gamma_i$ associated to positive flow direction
(recall lemma \ref{soursi}),
while keeping track of what $\gamma_i$ does to $z$.
The only difference is that here we were careful
to make sure $\gamma_i(z)$ did not converge to a certain stable
ideal point in the limit. This proof does not work at
all in the case $z$ is not ideal point of a leaf of $\oos$ or
$\oou$. 

\vskip .1in
\noindent
{\underline {Case 3}} $-$ $x = \varphi(z)$ where $z$ is not
ideal point of a leaf of $\oos$ or $\oou$.

This case is much more interesting.
Since $z$ is not ideal point of a leaf of $\oos$ or $\oou$,
by proposition \ref{options} there is a neighborhood system of
$z$ in
$\cd$ defined by a sequence of stable leaves, which can be
assumed to be all nonsingular.
Let $l_1$ be
one of these leaves. The construction here will be
inductive.
Let $W$ be the component of $\oo - l_1$ which has $z$ in its closure.
Let 
$\partial l_1 = \{ b_0, b_1 \}$.
Let $(b_0,z)$ be the interval of $\partial \oo$
contained in the closure of $W$ in $\cd$
and similarly define $(b_1,z)$.
Let $s$ be a leaf of $\oou$ intersecting $l_1$.
If
$s$ is near $b_0$ then all ideal
points of $s$ are near $b_0$ $-$ by the escape lemma (lemma \ref{duo}, part iii).
If $s$ is near $b_1$ then all ideal points are near $b_1$. The ideal points
of the prongs of $s$ entering $W$
vary monotonically in $\partial \oo$ as one moves $s$
across $l_1$. Since no unstable leaf has ideal point $z$ and
the leaf space of $\oou$ is Hausdorff, then
there is a single leaf $-$ call it $s_1$ intersecting $l_1$ and
having at least one prong contained in $W$ with an ideal point
in $(b_0,z)$ and another prong with ideal point in 
$(b_1,z)$, see fig. \ref{split}, a.
Let $p_1$ be the singular point in this leaf which has to be in $W$.
Let $v_1$ be the ideal point of $\oou(p_1)$ in $(b_0,z)$ closest
to $z$ and $u_1$ the one in $(b_1,z)$ closest to $z$.
Let $a_1$ be the ideal point of the (unique) prong of $\oou(p_1)$
intersecting $l_1$.

\begin{figure}
\centeredepsfbox{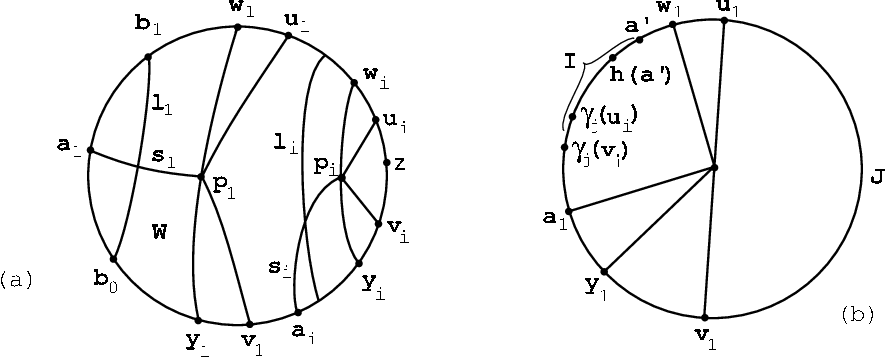}
\caption{
a. Spliting in the stable leaves,
b. Mapping back to a compact region.}
\label{split}
\end{figure}

We can now proceed inductively: assuming that $l_{i-1}$ has
been chosen and $s_{i-1}, p_{i-1}$ have been constructed, 
let $l_i$ be a stable leaf separating
$z$ from $\oou(p_{i-1})$. As before construct $s_i, p_i,
u_i, v_i$,
see fig. \ref{split}, a. 
Let $w_i$ be the ideal point of 
$\oos(p_i)$ in $(u_i, b_1)$ closest to $b_1$
and $y_i$ the 
ideal point of $\oos(p_i)$ in $(v_i,b_0)$ closest
to $b_0$ $-$ do this also for $i = 1$. 
There are such points because $\oou(p_i)$ intersects $l_i$ which
is a stable leaf.
Let $a_i$ be the ideal point of the prong of $\oou(p_i)$ which 
intersects $l_i$.

We will now take subsequences at will and rename points and
transformations, in order to simplify notation.
Every $p_i$ is singular, so 
up to subsequence assume the $p_i$ are all translates
of each other.
Hence there are $\gamma_i$ in $\pi_1(M)$ with $\gamma_i(p_i) = p_1$.
Up to another subsequence either every 
$\gamma_i$ preserves orientation in $\oo$,
or every $\gamma_i$ reverses orientation in $\oo$. In the second
case throw out $p_1$ (that is start with $p_2$ which will be renamed $p_1$
and also rename the
$\gamma_i$ to have $\gamma_i(p_i) = p_1$ for the new $p_1$, etc..).
So we can assume that every $\gamma_i$ preserves orientation in $\oo$.
Up to a further subsequence assume that $\gamma_i(a_i) = a_1$
(where as before throw out initial terms and rename if necessary).
Under these conditions, it now follows that 
$\gamma_i(u_i) = u_1$,
$\gamma_i(v_i) = v_1$, 
$\gamma_i(w_i) = w_1$,
$\gamma_i(y_i) = y_1$.
Let $(a_i,w_i)$ be the interval in $\partial \oo$
defined by $a_i, w_i$ and not containing $z$.
Assume also up to subsequence that for $j > i$ then
$y_i, v_i, u_i, w_i$ are in $(a_j, w_j)$, see
fig. \ref{split}, a.
This is because there are 2 possibilities for the placement of
$a_j$.

Since $p_1$ is singular, let $h$ 
be a generator of the isotropy group of 
$p_1$ which leaves all prongs of $\oos(p_1)$ (and hence of $\oou(p_1)$)
invariant.
Ideally we would like to obtain transformations which send
more and more of \ $\partial \oo - \{ z \}$ \ to
a compact set in $(a_1,w_1)$.
However in order to simplify the argument and the notation with indices
we will prove that this is true for a {\underline {fixed}} compact
set of $\oo - \{ z \}$ and then use that and the convergence group
property to show that $\varphi(z)$ is a conical limit point.
For each $i$ let $T_i$ be the closed interval of $\partial \oo$ defined
by $u_i, v_i$ and not containing $z$.

\vskip .15in
For the remainder of the proof we fix $i$ 
very big and let $C = T_i$ $-$ this is almost all of $\partial \oo
- \{ z \}$.
Let $a'$ be a point in
$(a_1,w_1)$.
By construction for any $j$ then
$\gamma_j(u_j) = u_1, \ \gamma_j(a_j) = a_1$. Let $j > i$.
Since $u_i$ is in $(a_j, w_j)$ then $\gamma_j(u_i)$ is
in $(a_1,w_1)$.
Now for each $j > i$ there is a single $n_j$ in ${\bf Z}$ so that

$$h^{n_j}(\gamma_j(u_i)) \ \ \ \ {\rm is \ in} \ \ \ \ 
[a', h(a'))$$

\noindent
where $[a', h(a')]$ is the subinterval of
$[a_1,w_1]$ bounded by these points.
Suppose that $w_1$ is a repelling fixed point of $h$
(that is, $h$ is associated to backwards flow direction).
Since $\gamma_j$ preserves orientation in $\cd$ 
then $t_j = h^{n_j}(\gamma_j(v_i))$ is closer to 
$a_1$ in $[a_1,w_1]$ than $h^{n_j}(\gamma_j(u_i))$ is.
We claim that $t_j$ is in a compact set $I$ of
$(a_1, w_1)$ as $j$ varies (in particular $\gamma_j(C) \subset I$).
Otherwise there are $j$ with $t_j$ arbitrarily close to $a_1$.
Here 

$$h^{n_j} \gamma_j(\oou(p_i))$$

\noindent
is an unstable leaf with a point in
$[a',h(a')]$ and another very close to $a_1$. 
Take a subsequence and find in the limit
an unstable leaf $\delta$ with an ideal
point in $[a',h(a')]$ and another in $a_1$ $-$ a consequence of the
escape lemma (lemma \ref{duo}, part iii).
Since $\delta$ is not
$\oou(p_1)$  this would force the existence
of perfect fits, contradiction.
Hence there is a compact subinterval $I$ in
$(a_1,w_1)$ with $t_j$ always in $I$.
We now define the transformations

$$g_j \ = \  h^{n_j} \gamma_j, \ \ \ \ j  \ >  \ i \ \ \ \  {\rm hence} \ \ \ \
g_j(C) \ \subset \ I.$$


\noindent
Let $J$ be the closed interval of $\partial \oo$
bounded by $u_1, v_1$ and not containing $a_1$. 
Then $g_j(z)$ is in $J$ for any $j \geq 2$ so up to a subsequence
we may assume that
$g_j(z)$ converges to a point $c$ in $J$.

\vskip .1in
We will show that there is a subsequence of
$(g_j)$ which proves that
$x$ is a conical limit point.

We first claim that $\varphi(I), \varphi(J)$ are disjoint.
Suppose that $\varphi(I)$ intersects $\varphi(J)$.
Then there has to be a leaf of $\oos$ or $\oou$ with ideal
points in both $I$ and $J$. Consider first the
unstable case. The endpoints of $J$ are 
ideal points of $\oou(p_1)$. The other ideal points
of $\oou(p_1)$ are not in $I$ $-$ by construction of
the interval $I$ in $(a_1,w_1)$. Any other leaf of $\oou$
either has all ideal points in $J$ or has no ideal point
in $J$. Hence no unstable leaf has ideal points in
$I$ and $J$.

Consider now stable leaves: 
$\oos(p_1)$ has one ideal point in $J$ and all
others in the interval of $\partial \oo$ defined
by $w_1, y_1$  and containing $I$.
Hence $\oos(p_1)$ it does not
have an ideal point in $I$. Let $r$ be any other leaf of $\oos$.
If $r$ has an ideal point in $J$ then $r$ is separated
from the interval $I$ by $\oos(p_1)$ $-$ hence $r$
cannot limit in $I$.
We conclude that $\varphi(I), \varphi(J)$ are 2 disjoint compact
subsets of $\rr$.

\vskip .1in
Recall that $(g_n(z))$ converges to $c$ and $x = \varphi(z)$.
Hence in $\rr$ the
sequence $(g_n(x))$ converges to $\varphi(c) \in \varphi(J)$.
In theorem \ref{cogr} we have already shown that $\pi_1(M)$ 
acts as a convergence group on $\rr$, so assume
up to subsequence that $(g_n)$ 
has a source/sink pair (for notational simplicity we still denote
this subsequence by $(g_n)$). 
That means there are $a, b \in {\cal R}$ so that 
$g_n(A)$ converges to $b$ for any compact set $A$ of ${\cal R} - \{ a \}$.
In particular if we find three distinct points $d_1, d_2, d_3$ of
${\cal R}$ so that $g_n(d_1), g_n(d_2)$ converge to $e_1$ and
$g_n(d_3)$ converges to $e_2$ with $e_1 \not = e_2$, then
$d_3$ is the source and $e_1$ is the sink.

The image $\varphi(C)$ contains infinitely many points, so take
3 distinct points $d_0, d_1, d_2$ in $\varphi(C)$. By the above,
for at least two of these points the sequence $(g_n(d_k))$ converges.
So assume wlog that $(g_n(d_1)), (g_n(d_2))$ converge $-$ the limit
is in $\varphi(I)$.
The sequence $(g_n(x))$  also converges and the limit
is in $\varphi(J)$.
As $\varphi(I), \varphi(J)$ are disjoint, it follows that
the limits of $(g_n(d_1)), (g_n(d_2))$ have to be the same point $t$.
By the previous paragraph $t$ is the sink and $x$ is the source
for the sequence $(g_n)$ acting on ${\cal R}$.
Since $t \in \varphi(I), x \in \varphi(J)$ it follows
that $t \not = x$. Hence the sequence $(g_n)$ 
of $\pi_1(M)$
shows that $x$ is a conical limit point.

This shows that all points of $\rr$ are conical
limit points for the action of $\pi_1(M)$.
Hence $\pi_1(M)$ acts
as a uniform convergence group in $\rr$. 
This finishes
the proof of theorem \ref{unifco}.
\end{proof}

We now analyse the space $\mi \cup \rr$. 
We first establish some
notation.
Let 

$$\eta: \  \cd \times [-1,1] \ \ \rightarrow \ \ \mi \cup \rr$$

\noindent
be the projection map. 
Recall also the sphere filling map $\varphi: \partial \oo \rightarrow {\cal R}$.
We consider the quotient topology in $\mi \cup \rr$.
Let $\ct$ be this topology.
Recall that 
$\partial (\cd \times [-1,1])$ is a sphere, let 
$\eta_1 = \eta \ | \ \partial (\cd \times [-1,1])$.
With the subspace topology from $\ct$, then $\rr$ is a sphere also.
We stress that in all arguments here we implicitly identify 
$\mi$ with $\oo \times (-1,1)$ and in particular also
think of $\mi$ as a subset of $\cd \times \mmp$.

Notice that $\pi_1(M)$ naturally acts on $\mi \cup \rr$ by
homeomorphisms as it preserves stable and unstable foliations.
Our main goal to finish this section is to show that this action
is a convergence group action.

One problem is that it is hard to verify directly whether a set 
in $\mi \cup \rr$ is open or not.
To make it more explicit we define another topology $\ct'$ in
$\mi \cup \rr$ and then show it is the same as the quotient topology.
The new topology will be defined using neighborhood systems.
Recall \cite{Ke} chapter 1 that a neighborhood
system $\uu_x$ of a point $x$ is a collection satisfying:

1) If $U$ is in $\uu_x$ then $x$ is in $U$.

2) If $U, V$ are in $\uu_x$ then $U \cap V$ is
in $\uu_x$.

3) If $U$ in $\uu_x$ and $U \subset V$ then $V$ is in $\uu_x$.

\noindent
Define $U$ to be open if $U$ is a neighborhood of 
any of its points. 
This defines a topology in the space.


\vskip .1in

\begin{define}{(neighborhood systems in $\mi \cup \rr$)}{}
Let $\Phi$ be a pseudo-Anosov flow without perfect fits and
not topologically conjugate to suspension Anosov.

i) If $x$ is in $\mi$ then $V$ is in $\uu_x$ if $V$ contains
an open set  of $\mi$ (with its usual topology) containing
$x$,

ii) Let $x$ in $\rr$ so that $\varphi^{-1}(x) = \{ b \}$,
a single point. 
The point $b$ of $\partial \oo$ is not identified to any 
other point of $\partial \oo$, hence $b$ is not an ideal point of a leaf
of $\oos$ or $\oou$.
In this case $b$ has a neighborhood system
in $\cd$ defined by sequences of nonsingular stable or
unstable leaves. Let $l$ be one such leaf
and $U_l$ the corresponding open set of $\cd$,
as in definition \ref{canon}, where $b$ is in $U_l$.
Let $V_l = U_l \times \mmp$ a subset of $\cd \times \mmp$.
We say that $V$ is in
$\uu_x$ if for some $l$ as above then $V_l \subset \eta^{-1}(V)$.
Notice $\eta^{-1}(V)$ is a subset of $\cd \times \mmp$.

iii) Let $x$ in $\rr$ with $\varphi^{-1}(x)
= \{ a_1, ..., a_n \}$.
For simplicity assume that 
$a_1, ..., a_n$ are the ideal points of a stable leaf
$l$. Let $g$ be the cellular decomposition element of ${\cal R}$ of
$\partial (\cd \times \mmp)$
associated
to $l$ (that is $g = l \times \{ 1 \} \cup  \ \cup_i (\{ a_i \} \times \mmp)$ \
or equivalently  $g$ is identified to the point $x$).
For each $i$,
choose $r_i$ unstable leaves defining small neighborhoods of $a_i$ in
$\cd$.
Let $V_{r_i} = U_{r_i} \times \mmp$ as in ii), where $a_i$ is in $U_{r_i}$.

Let $l_1, ..., l_n$ stable leaves ($\oos$) very near each line
leaf 
of $l$ and so that for each $i$, then
$l_i, l_{i+1}$ intersect
$r_i$ transversely ($i \ mod(i_0)$). 
Then $l_1, ..., l_n, r_1, ...., r_n$ bound a compact
region $B$ in $\oo$. Choose any section $\tau: B \rightarrow \mi$
of $\Theta$ restricted to $B$.
Let $H_{\tau}$ be the union of $B \times \{ 1 \}$
together with the
set of points $w$ in $\mi$ (or in $\oo \times (-1,1))$
with $w = \wwp_t(b)$ for some $b$ in $\tau(B)$ and $t \geq 0$.
Let $\delta$ denote the collection $(l_1, ...., l_n, r_1, ...
, r_n, \tau)$. We use the notation $A_{\delta}$ to denote
the following:

$$A \ = \ A_{\delta} \ = \ A(l_1, ..., l_n, r_1, ..., r_n, \tau)
\ = \ H_{\tau} \cup V_{r_1} \cup ... \cup V_{r_n} $$

\noindent
Let $\uu_x$ be the collection of the sets $Z$ so
that for some $\delta$ as above then 
$A_{\delta} \subset \eta^{-1}(Z)$.

In the case of ideal points of unstable leaves, one
switches stable and unstable objects and chooses points
flow {\underline {backwards}} from a section and backward
ideal points.
\label{neigh}
\end{define}

\begin{figure}
\centeredepsfbox{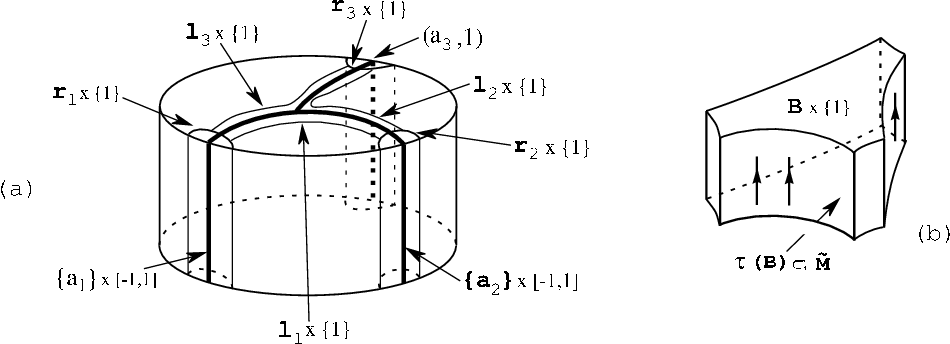}
\caption{
a. The neighborhoods of certain points, b. 
Flow forward of sections.}
\label{sets}
\end{figure}

\vskip .1in
\begin{lemma}{}{} The collection $\uu_x$ for $x$ in $\mi \cup \rr$
defines a neighborhood system and consequently
a topology $\ct'$ in $\mi \cup \rr$.
\end{lemma}

\begin{proof}{}
For $x$ in $\mi$ this is clear.
In the other 2 cases it is easy to see that properties 1) and 3) 
of neighborhood systems always hold:
3) is obvious by definition and 1) holds
because the cell decomposition elements (in $\partial (\cd
\times \mmp)$) are always contained in the sets $V_l$ or $A_{\delta}$.

We now check property 2). Suppose first that $x$ is of type ii).
Let $x = \varphi(b)$.
Let $V_1, V_2$ in $\uu_x$, with
$V_1$ defined by $l$ and $V_2$ defined by $r$ leaves of $\oos$ or 
$\oou$. Then there
is $l'$ in $\oos$ or $\oou$ so that $l' \cup \partial l'$ separates
$b$ from $r \cup l$ in $\cd$. Then $U_{l'}$ is contained
in $U_l \cap U_r$ and we are done. 

Let now $x$ of type iii).
Let $U_1, U_2$ be neighborhoods of $x$, where $U_i$ contains
$A_i$ of the form
$A_i = 
A(l^i_1, ..., l^i_n, r^i_1, ..., r^i_n, \tau_i)$
as in definition \ref{neigh}, so
that for each $i$, $l^1_i, l^2_i$ are close to the same line
leaf of $l$ and $r^1_i, r^2_i$ define small neighborhoods of $a_i$.
Choose $l^3_i$ closer to $l$ than both 
$l^1_i$ and $l^2_i$ and $r^3_i$ closer to $a_i$
than both $r^1_i$ and $r^2_i$.
Let $B_3$ be the compact region of $\oo$ defined by the $l^3_i, r^3_i$.
Choose a section $\tau_3$ in $B_3$ so that in
the intersection 
$B_3 \cap  (B_1 \cup B_2)$ then $\tau_3$ is greater than
$max(\tau_1,\tau_2)$.
Then $A_3 = A(l^3_1, ..., l^3_n, r^3_1, ..., r^3_n, \tau_3)$
is in $\uu_x$ and $A_3 \subset A_1 \cap A_2 \subset U_1 \cap U_2$.
Hence $\uu_x$ is a neighborhood system 
for $x$ in $\mi \cup \rr$.

Therefore the collection $\{ \uu_x, x \in \mi \cup \rr \}$ defines a topology
in $\mi \cup \rr$.
\end{proof}

\begin{lemma}{}{}
The quotient topology $\ct$ in $\mi \cup \rr$ and the 
neighborhood system topology $\ct'$ are the same topology.
This implies that the quotient topology in $\rr$ and the
subspace topology from $\ct'$ in $\rr$ are also
the same topology.
\end{lemma}

\begin{proof}{}
First let $U$ in $\ct'$ and let $x$ in $U$.
If $x$ is in $\mi$, then (i) of definition \ref{neigh} shows that
there is $V$ open in (usual topology) of $\mi$ with
$x \in V \subset U$. 
If $x$ is in $\rr$ let $g = \eta^{-1}(x)$.
By construction if $x$ is of type ii)
or iii) as in definition \ref{neigh}, then $\eta^{-1}(U)$ contains an open
set in $\cd \times \mmp$
which contains $g$.
This shows that $\eta^{-1}(U)$ is an open set in $\cd \times \mmp$ and
hence $U$ is in $\ct$.

Conversely let 
$U$ in $\ct$. Then $\eta^{-1}(U)$ is
open in $\cd \times \mmp$.
Let $x$ in $U$.
If $x$ is in $\mi$, then $x$ is in the open set
$\eta^{-1}(U) \cap \mi \subset \eta^{-1}(U)$ so
$\eta^{-1}(U)$ is in $\uu_x$.

Suppose then that $x$ is in $\rr$ and 
 let $g$ the cell element of ${\cal R}$ associated to
$x$.
For simplicity we assume that 
$x$ is of type iii) in definition \ref{neigh}, as type ii)
is analogous and easier to deal with.
Let 
$l$ (as in part iii) of def. \ref{neigh})
be the leaf of 
(say) $\oos$ with $l \times \{ 1 \}$ a subset of
$g$.
Then $\eta^{-1}(U)$ is 
an open set in $\cd \times \mmp$
containing $g$.
For any ideal point $b$ of $l$, then $\eta^{-1}(U)$ contains
an open neighborhood of $b \times \mmp$ in
$\cd \times \mmp$. Since $b$ is a stable
ideal point, there is an unstable leaf $z$ defining
a small neighborhood of $b$ in $\cd$ so that
$V_z 
\subset \eta^{-1}(U)$. 
We also consider for each line leaf of $l$ a regular
leaf $e$ of $\oos$ close to this line leaf.
Choose each $e$ sufficiently close to $l$ so that 
these $e$'s and the $z$'s as above define a compact polygon
$B$ in $\oo$. 
As $\eta^{-1}(U)$ is open and contains $l \times \{ 1 \}$, 
it follows that 
if the $e$'s are sufficiently close to $l$ and 
the $z$'s sufficiently close  to $\partial l$, then $B \times \{ 1 \}
\subset \eta^{-1}(U)$.
As $B$ is compact, there 
is a high enough section $\tau: B \rightarrow \mi$
so that $H_{\tau} \subset \eta^{-1}(U)$.
This shows that $\eta^{-1}(U)$ contains
one set of form $A_{\delta}$ as in iii) of definition
\ref{neigh} and so $U$ is in ${\cal U}_x$.
Since $U$ is in $\uu_x$ for any $x$ in $U$, it follows
that $U$ is open with respect to $\ct'$.
Hence $\ct$ is equal to $\ct'$.
\end{proof}

\begin{lemma}{}{}
The space $\mi \cup \rr$ is compact.
\end{lemma}

\begin{proof}{}
Let $\{ Z_{\alpha} \}_{\alpha \in {\cal I}}$ be an open cover
of $\mi \cup \rr$.
This provides an open cover of $\rr$ which is compact.
Hence there is a finite subcollection
$Z_{\alpha_1}, ..., Z_{\alpha_n}$ whose union contains 
$\rr$.
Then

$$C \ = \ \mi \cup \rr - ( \bigcup_{i=1}^n \ Z_{\alpha_i})
\ \ \subset \ \ \mi$$

\noindent
is closed. Since the topology
in $\mi$ is the same as the induced topology from
$\mi \cup \rr$, it follows that $C$ is closed in $\mi$
and hence compact and it has a finite subcover.
This finishes the proof.
\end{proof}

Here is another way to see that $\pi_1(M)$ acts
on $\mi \cup \rr$: Let $\gamma$ in $\pi_1(M)$. Then $\gamma$
takes sets of the form $V_l$ (of (ii) of definition \ref{neigh}) 
for $l$ in $\oos$ or
$\oou$ to $V_{\gamma(l)}$. 
Sections $\tau: B \rightarrow \mi$ 
over compact sets $B$ in $\oo$ are taken to sections
over compact sets $\gamma(B)$ by $\gamma$.
Hence $\pi_1(M)$ preserves the collection of sets described in 
\ ii), iii) of definition \ref{neigh}.
Therefore $\gamma$ takes neighborhoods  $\mi \cup \rr$
to neighborhoods and
consequently $\pi_1(M)$ acts  by homeomorphisms
on $\mi \cup \rr$.

We stress that it is hard to find open  sets
in $\mi \cup \rr$ explicitly: for example if $l$ is
a nonsingular leaf of $\oos$, with corresponding
open set $V_l$ in $\cd \times \mmp$, it is not
true that $\eta(V_l)$ is open in $\mi \cup \rr$, because $V_l$ is
not saturated by the equivalence relation defining
the quotient: Certainly $V_l \cap \mi$ is open
in $\mi$ and $V_l \cap (\cd \times \{ 1 \})$ 
is both open and saturated in $\cd \times \{ 1 \}$.
However $V_l \cap (\cd \times \{ -1 \})$ is
{\underline {not}} saturated. Take any leaf $s$ of 
$\oou$ intersecting $l$. Then $s \times \{ -1 \}$ 
interects $V_l$ but is not contained in $V_l$.
Those leaves $s \times \{ -1 \}$ would have to be
contained in a saturation of $V_l$. But their
ideal points propagate through $\partial \oo \times \mmp$
and then propagate in the top $\cd \times \{ 1 \}$ through
stable leaves.

\begin{lemma}{}{}
The space $\mi \cup \rr$ is first countable.
\label{first}
\end{lemma}

\begin{proof}{}
We only need to check this for 
$x$ in $\rr$
since $\mi$ is a manifold and is open in $\mi \cup \rr$.
Suppose $\varphi^{-1}(x) = \{ a_1, ..., a_{i_0} \}$, all
ideal points of a stable leaf $l$.
The other cases are either similar or simpler.
For each $1 \leq i \leq i_0$,  we will construct a nested sequence of
unstable leaves
$(s^n_i)_{n \in {\bf N}}$ 
forming a master sequence
defining $a_i$.
For each line leaf $l_i$ of $l$ we will construct a nested sequence of nonsingular
stable leaves
$(l^n_i)_{n \in {\bf N}}$ 
converging to $l_i$ in that
sector of $l$.
Suppose that $l^n_i, l^n_{i+1}$ ($i \ {\rm mod}(i_0)$) bound a small segment 
$T_i^n$ in $\partial \oo$ containing $a_i$ in its interior.
We do the construction so that for all $n$ and $i$,  the leaves 
$l^n_i, l^n_{i+1}$ intersect $s^n_i$ transversely.
Then for each $n$ 

$$l^n_1, ..., l^n_{i_0}, s^n_1, ..., s^n_{i_0}$$

\noindent
defines a compact set $B_n$ in $\oo$.
It is not true that $B_j \subset B_i$ if $j > i$.
Fix a section $\tau_1: B_1 \rightarrow \mi$.
We will choose sections $\tau_n: B_n \rightarrow \mi$
so that for each $n$, $\tau_n  (B_{n-1} \cap B_n)$ 
is flow forward of $\tau_{n-1}(B_{n-1} \cap B_n)$
and the flow length from $\tau_1(B_n \cap B_1)$ to $\tau_n(B_n \cap B_1)$
goes to infinity uniformly in $n$.

Let $A_n = A(l^n_1, ..., l^n_{i_0},
r^n_1, ..., r^n_{i_0}, \tau_n)$.
Notice that $\eta(A_n)$ is not open in $\mi \cup \rr$
because $A_n$ is not saturated. 
However we will choose $A_n$ inductively so that there is an
open set $U_n$ in $\mi \cup \rr$ satisfying

$$ \eta(A_{n-1}) \ \superset \ U_n \ \superset \ \eta(A_n)$$

Here is the construction.
Suppose that $l_1^{n-1}, ..., l_{i_0}^{n-1}, s^{n-1}_1, ..., s^{n-1}_{i_0}$
have been chosen. 
We choose one set $l^n_i, 1 \leq i \leq i_0$ closer to $l$ than $l^{n-1}_i$ and
$s^n_i$ closer to $a_i$ than $s^{n-1}_i$.
We will adjust these choices as needed.

Let $x$ in $\eta(A_n)$.
Certainly we can choose the section $\tau_n$ so that if $x$ is in $\eta(A_n)$ and
$x$ is in $\mi$ then $x$ is in the interior of $\eta(A_{n-1})$.
Therefore assume that $x$ is in $\rr$ and let $y$ in $\eta^{-1}(x)$.
There are 3 possibilities:

\vskip .05in
A) First suppose that $y$ is in $\oo \times \{ -1 \}$.

Then $y$ is in the region of $\cd \times \{ -1 \}$ bounded by some 
$s^n_i \times \{ -1 \}$, which is strictly smaller than the
region bounded by $s^{n-1}_i \times \{ -1 \}$.
Let $v$ be the leaf of $\oou$ with $y$ in $v \times \{ -1 \}$.
Then $v$ is contained in the region
$U_{s^{n-1}_i}$ and hence there is a set $A_{\delta}$ as in 
(iii) of definition \ref{neigh} associated to $v$ and
so that 

$$A_{\delta} \  \subset  \ U_{s^{n-1}_i} \ \subset \ A_{n-1}$$

\noindent
By definition $\eta(A_{n-1})$, is in $\uu_x$ because

$$\eta^{-1}(\eta(A_{n-1})) \ \superset \ A_{n-1} \ \superset \ A_{\delta}$$

\vskip .05in
B) The second case is that $y$ is in $\partial \oo \times \mmp$,
but $y$ is not equivalent to any point in
$\oo \times \{ 1 \}$ or $\oo \times \{ -1 \}$ $-$ that is, $y$ does not
come from an ideal point of a leaf of $\oos$ or $\oou$.
Then $y$ is in some $V_{s^n_i}$ and by part (ii) of definition \ref{neigh},
$\eta(A_{n-1})$ is a neighborhood of $x$ in $\mi \cup \rr$.

\vskip .05in
C) Finally suppose that $y$ is in $\oo \times \{ 1 \}$.

If $y$ is in the region of $\cd$ bounded by the 
$l^n_i, 1 \leq i \leq i_0$, then the proof as in part A) applies.
The last case to analyse is $y$ is $V_{s^n_i}$ for some
$i$.
Here $y$ is in $u \times \{ 1 \}$ with $u$ leaf of $\oos$.
In this case we adjust $s^n_i$ so that its endpoints are in the
open interval $T^{n-1}_i$.
Then all stable leaves near $u$ are in the region bounded 
by $l^{n-1}_i, l^{n-1}_{i+1}$. This shows that 
$\eta(A_{n-1})$ is a neighborhood of $x = \eta(y)$.

The modification in part C)  makes the $U_{s^n_i}$ smaller and hence one has to 
rechoose the $l^n_i$ closer to $l$ accordingly so that 
$s^n_i$ intersects both $l^n_i$ and $l^n_{i+1}$. With this modification 
it follows that $\eta(A_{n-1})$ is a neighborhood of any point $x$ in
$\eta(A_n)$ so there is an open set $U_n$ in $\mi \cup \rr$ with
$\eta(A_{n-1}) \superset U_n \superset \eta(A_n)$.

As the sequence $(l^n_i)$ converges to a line leaf of $l$ for each $i$,
$(s^n_i)$ converges to $a_i$ and $\tau_n(B_n)$ escapes in the positive 
direction, then
it is now clear that the collection $\{ U_n \}_{n \in {\bf N}}$ forms
a countable basis for the topology of $\mi \cup \rr$ at  $x$.
\end{proof}

This result will be used in section 5.

Finally we show 
that the action of $\pi_1(M)$ 
on $\mi \cup \rr$ is a convergence group action.
The description of the topology in $\mi \cup \rr$ using
neighborhood systems is extremely useful for this result.

\begin{theorem}{}{}
Let $\Phi$ be a pseudo-Anosov flow without perfect fits and
not conjugate to a suspension Anosov flow. Then
the induced action of $\pi_1(M)$ on $\mi \cup \rr$
is a convergence group action.
\label{bige}
\end{theorem}

\begin{proof}{}
Let $(\gamma_n)_{n \in {\bf N}}$ be a sequence of distinct
elements in $\pi_1(M)$.
Since the action of $\pi_1(M)$ on $\rr$ is a convergence group
action, then up to subsequence we can assume there
are $x, y$ in $\rr$ with $(\gamma_n)$
converging locally uniformly to $x$ in $\rr - \{ y \}$.
We want to show that $(\gamma_n)$ converges locally
uniformly to $x$ when acting on $(\mi \cup \rr) - \{  y \}$.
Let $C$ be a compact set in $\mi \cup \rr - \{ y \}$.
Recall the surjective map $\varphi: \partial \oo \rightarrow \rr$.

\vskip .18in
\noindent
{\underline {Case 1}} $-$ $\varphi^{-1}(y) = \{ e \}$ $-$ a single point.

Then $\eta^{-1}(y)$ is a vertical segment in
$\partial \oo \times I$. For any
neighborhood
$U$ of $y$ in $\mi \cup \rr$, \
there is $l$ an unstable (or stable) leaf defining a small
neighborhood of 
$e$ in $\cd$ so   
that $V_l \subset \eta^{-1}(U)$, $V_l$ as in definition
\ref{neigh}.
If $C$ is disjoint from $U$ then

$$\eta^{-1}(C) \ \subset \ \cd \times I - V_l$$

\noindent
Let $Z$ be the closure of the segment of $(\partial \oo -
\partial l)$ not containing $e$ (this is almost
all of $\partial \oo$). 
By the source/sink property of $y, x$ for the sequence
$(\gamma_n)$ acting on $\rr$, 
the set $\gamma_n(Z)$ 
is very near $\varphi^{-1}(x)$ for $n$ big. 
As $\gamma_n(Z)$ is a segment in $\partial \oo$, 
then there is a single 
point $b$ in
$\varphi^{-1}(x)$ with
$\gamma_n(Z)$ near $b$ for $n$ big.
It follows that 
$\gamma_n(\cd \times I - V_l)$ is very near $\{ b \} \times \mmp$
in $\cd \times \mmp$
and so $\gamma_n(\eta^{-1}(C))$ is very near $\{ b \} \times \mmp$
in $\cd \times \mmp$.
We conclude that 
$\gamma_n(C)$ is very near $x = \eta( \{ b \} \times \mmp)$ in $\mi \cup \rr$ as
desired. This finishes the analysis of case 1.


\vskip .2in
\noindent
{\underline {Case 2}} $-$ 
$\varphi^{-1}(y) = \{ a_1, ..., a_{i_0} \}$, with $i_0 \geq 2$.

Suppose for simplicity that 
$\{ a_1, ..., a_{i_0} \}$ 
are the ideal points of $\oos(p) = l$ for some $p$ in $\oo$.
Let $C$ be a compact set in $\mi \cup \rr - \{ y \}$.
As before there are $\{ l_i \}_{1 \leq i \leq i_0}$ regular leaves
of $\oos$
very near the line leaves of $l$ and there are
$\{ r_i \}_{1 \leq i \leq i_0}$,
regular leaves of $\oou$ defining small neighborhoods of $a_i$
so that the $l_i$'s together with the $r_i$'s define a compact
set $B$ in $\oo$ and there is a section $\tau: B \rightarrow \mi$
with 

$$A  \ = \ A(l_1, ..., l_{i_0}, r_1, ..., r_{i_0}, \tau) 
\ \ \ \ {\rm and \ so \ that} \ \ \ \  \eta^{-1}(C) \ \subset  \
\cd \times \mmp   - 
A$$

\noindent
Assume that $r_i$ intersects $l_i, l_{i+1}$ (mod $i_0$) and
has ideal points near $a_i$.
Since $r_i, l_i$ are regular we need to be careful.
Let $\widetilde r_i$ be the component of $\oo - r_i$ which
has $a_i$ in its closure (in $\cd = \oo \cup \partial \oo$).
Let also $\widetilde l_i$ be the component of $\oo - l_i$ not
containing the other $l_j$.
Then consider the sets
$U_{r_i}$  and
$U_{l_i}$ 
as in definition \ref{canon}.
The endpoints of $l_i$ bound a closed interval $I_i$ in $\partial \oo$
contained in the closure of $U_{l_i}$ (they do not contain any $a_j$).
Similarly the endpoints of $r_i$ bound a very small closed interval $J_i$ 
in $\partial \oo$ containing
$a_i$. 
As in definition \ref{neigh},
let $V_{r_i} = U_{r_i} \times \mmp$
and let $H_{\tau} = \{ \wwp_t (z) \ | \ z \in \tau(B) \ \ 
{\rm and} \ \ t \geq 0\} \cup  \ (B \times \{ 1 \})$.
The sets $C$ and $A$ will be fixed for the rest of the proof of case 2.

Let $\varphi^{-1}(x) = \{ b_1, ...,  b_{j_0} \}$.

\vskip .2in
\noindent
{\underline {Case 2.a}} $-$ The union \ $\cup _i \ \gamma_n(\partial l_i)$ \
is eventually (with $n$) always very near a single
point $b_1$ in $\varphi^{-1}(x)$.

Since $\gamma_n$ restricted to compact sets
of \ $(\partial \oo - \varphi^{-1}(y))$ \ 
has image very close to $\varphi^{-1}(x)$ for $n$ big, 
it follows that
for all $i$, $1 \leq i \leq i_0$ then $\gamma_n(I_i)$ is very close
to $b_1$ in $\cd$.
This implies that $\gamma_n(V_{l_i})$ is very close
to $\{ b_1 \} \times \mmp$ in 
$\cd \times \mmp$. In addition since the $\gamma_n$ are homeomorphisms
of $\partial \oo$, then there is a single $i$ (assume for simplicity
that $i = 1$)
so that $\gamma_n(J_1)$ is almost all of
$\partial \oo$ and hence $\gamma_n(\partial \oo - J_1)$ is very close to $b_1$.
Notice that

$$\cd \times \mmp - (H_{\tau} \cup V_{r_1} \cup ... \cup V_{r_n}) \ \ \subset \ \ 
\cd \times \mmp \ - \ V_{r_1}$$


\noindent
By the above \ 
$\gamma_n(\cd \times \mmp - V_{r_1})$ 
\ is 
very close to $\{ b_1 \} \times \mmp$ in $\cd \times \mmp$.
It follows that $\gamma_n(C)$ is very close to $x$ for $n$ big.
This finishes the analysis in this case.


%
%
%
%
%

\vskip .2in
\noindent
{\underline {Case 2.b}} $-$ The union \ $\cup_i \ \gamma_n(\partial l_i)$ \
gets closer to more 
than one point in $\varphi^{-1}(x)$.

We first explain why the $b_i$ are ideal points of an 
{\underline {unstable}} leaf in this case.
To start we claim that,  for a single $i$,the ideal points
of $\gamma_n(l_i)$ are close to a single
point in $\varphi^{-1}(x)$ for $n$ big.
Let $c_1, c_2$ be the endpoints of $l_i$. 
If the claim is not true, then up to subsequence
the sequences
\ $(\gamma_n(c_1)), \ (\gamma_n(c_2))$ \ converge to
two distinct points $d_1, d_2$ in $\varphi^{-1}(x)$.
If follows that $\gamma_n(U_{l_i} \cap \partial \oo)$
contains most of a segment with endpoints 
$d_1, d_2$. This contradicts the fact that
$\gamma_n(U_{l_i} \cap \partial \oo)$ converges to points
in $\varphi^{-1}(x)$. This proves the claim.

The hypothesis of case 2.b implies that 
there is some $i$ so that the ideal points of 
$\gamma_n(l_i), \gamma_n(l_{i+1})$  are not close.
But $\gamma_n(r_i)$ intersects both of these leaves,
hence the escape lemma implies that 
up to subsequence $(\gamma_n(r_i))$ 
converges to a leaf $s$ of $\oou$.
The source/sink property for $y, x$ implies that
the ideal points of $\gamma_n(r_i)$ have
to be getting close to points in $\varphi^{-1}(x)$. It follows
 that $\varphi^{-1}(x) = \partial s$ 
with $s$ an unstable leaf,
as we desired to show.

For any neighborhood $W$ of $x$ in $\mi \cup \rr$ there is 
a set $D$ in $\cd \times \mmp$ as in definition \ref{neigh}:
\ $D$ is defined by $s_1, ..., s_{j_0}$ regular leaves of
$\oou$ near line leaves of $s$; also $t_1, ..., t_{j_0}$ 
regular leaves of $\oos$, where $t_j$ defines a small neighborhood $U_{t_j}$ 
of $b_j$ in $\cd$. The $s_j, t_j, 1 \leq j \leq j_0$ 
jointly bound a compact
set $B'$ in $\oo$, consider a section $\nu: B' \rightarrow \mi$
and $E_{\nu}$ the set of points flow {\underline {backwards}}
from $\nu(B')$ union $B' \times \{ -1 \}$:

$$E_{\nu} \ \ = \ \ \wwp_{(-\infty,0]} (\nu(B')) \ \cup \ (B' \times \mmp)$$

\noindent
Let 

$$D \ \  = \ \ D(s_1, ..., s_{j_0}, t_1, ..., t_{j_0}, \nu) \ \ 
= \ \ \big( \bigcup_{1 \leq j \leq j_0} \ V_{t_i} \big) \ 
\cup \ E_{\nu}$$

\noindent
Then there is such a $D$ so that $D \subset \eta^{-1}(W)$.
Fix one such $D$.
We want to show that $\gamma_n(C)$ is eventually contained in $W$
in $\mi \cup \rr$.
It suffices to show that 
$\gamma_n(\cd \times \mmp - A) \subset D$ \ in $\cd \times \mmp$.
In case 2.b an argument in $\mi$ will be needed.
For the fixed $B$ as above with section $\tau: B \rightarrow \mi$, 
let $E_{\tau}$ be the set of points
flow backwards from the section $\tau(B)$ union $B \times \{ -1 \}$
(just as $E_{\nu}$ was defined).
Hence
$B \times \mmp$ is the union of $E_{\tau},H_{\tau}$ and the intersection
of $E_{\tau}, H_{\tau}$ is equal to $\tau(B)$.
Notice that

$$\cd \times \mmp \ - \ A \ \ \subset \ \
 \big( \bigcup_{1 \leq i \leq i_0} \ V_{l_i} \big) \ 
\cup \ E_{\tau}$$

Choose 
leaves $l_i$ close enough to $l$ so that
the length of any segment of 
$\wlu \cap \tau(B)$ from $\tau(l_i)$ to
$\tau(l_j)$ 
is very small.
This yields a smaller neighborhood of $y$
($V_{l_i}$ is bigger)
and we show that the complement of this neighborhood of $y$ 
goes near $x$ under $\gamma_n$. 
For any $i$, 
the endpoints of $\gamma_n(l_i)$ 
converge to
a single point in
$\varphi^{-1}(x)$ as $n \rightarrow \infty$.
Hence $\gamma_n(l_i)$ also does
and so $\gamma_n(V_{l_i})$  gets very near 
$\varphi^{-1}(x) \times \mmp$
in $\cd \times \mmp$. 

\vskip .08in
In order to finish the proof in case 2.b
we need to analyse $\gamma_n(E_{\tau})$.
We modify the leaves $s_j$ to be close enough to $s$, 
and the leaves $t_j$ to be close enough to $b_j$ and extend the section
$\nu(B')$ so that any unstable segment in
$\nu(B')$ connecting $t_j \times \rrrr$ to $t_k \times \rrrr$ has
very large length.
This decreases the set $D$, so we still have $D \subset \eta^{-1}(U)$.

Up to subsequence suppose
there are $z'_n$ in $E_{\tau}$ so that 
$\gamma_n(z'_n)$ are not in $D$. 
If $c$ is an ideal point of $r_i$ in $\partial \oo$, then
$\gamma_n(c)$ converges to a point in
$\varphi^{-1}(x)$
in $\partial \oo$.
Let $C_n$ be the closed, connected region in $\cd$ bounded by
the union of the 
$\gamma_n(r_i), 1 \leq i \leq i_0$ union its ideal points.
Then $\gamma_n(z'_n)$ is in 
$C_n \times [-1,1]$.
The bottom of this set is $C_n \times \{ -1 \}$ which
is contained in $D$ for $n$ big.
Hence if $\gamma_n(z'_n)$ is not in $D$ the following happens:
First 
$\gamma_n(z'_n)$ is in $B' \times (-1,1)$,
in particular $\gamma_n(z'_n)$ is in $\mi = \oo \times (-1,1)$.
Second, as $\gamma_n(z'_n)$ is not in $E_{\nu}$ then 
$\gamma_n(z'_n)$ is flow forward
from a point in $\nu(B')$.
As $z'_n \in E_{\tau}$, flow $z'_n$ forward to a point $z_n$ in the section $\tau(B)$.
Hence $\gamma_n(z_n)$ is still flow forward of a point in $\nu(B')$.

Now consider the segment $v_n$ which is the intersection of 
$\wu(z_n)$ with the section $\tau(B)$.
By construction this
segment has arbitrarily small length and hence so does the 
segment $\gamma_n(v_n)$ in $\wu(\gamma_n(z_n))$ $-$ because
$\gamma_n$ acts as an isometry on $\mi$.
This segment $\gamma_n(v_n)$ is entirely flow
forward of $\nu(B')$.
Flow $\gamma_n(v_n)$ backwards until it hits the section
$\nu(B')$. The unstable length gets decreased
when flowing backwards or at least it does not
increase too much, so it is a small length.

The segment $v_n$ has endpoints in $l_i \times \rrrr$ and
$l_j \times \rrrr$ for some $i, j$.
The endpoints of $\gamma_n(v_n)$ are in $G_i = \gamma_n(l_i \times \rrrr)$
and $G_j = \gamma_n(l_j \times \rrrr)$ which,
for $n$ sufficiently big, 
are contained in the union of
$V_{t_k}, 1 \leq k \leq j_0$. Notice that the boundary 
of $V_{t_k}$ 
is the stable leaf $t_k \times \rrrr$.
If both $G_i$ and $G_j$ are contained in the same $V_{t_k}$ this
forces $\gamma_n(l_i)$ to be contained in $V_{t_k}$ because its endpoints
are in this set and an unstable leaf cannot intersect 
the stable leaf boundary more than once.
But then $\gamma_n(z'_n)$ is in $D$ and we finish the analysis.
The remaining possibility is that
$G_i$ is in some $V_{t_k}$ and $G_j$ is in
some $V_{t_m}$ with $j \not = m$. 
Therefore $\gamma_n(v_n)$ flows back to a segment which has a subsegment
from $t_k \times \rrrr$ to $t_m \times \rrrr$ in $\nu(B')$. This 
subsegment
has fairly small length and this contradicts the choice
of leaves $\{ s_j, t_j, 1 \leq j \leq j_0 \}$ and the section $\nu$.
This shows that $\gamma_n(z'_n)$ are in $D$ contradiction to
assumption.

This shows that $\gamma_n(E_{\tau})$ is contained in 
$D \subset \eta^{-1}(W)$.
Hence in $\mi \cup \rr$, the sets
$\gamma_n(\mi \cup \rr - \{ y \})$ converge
locally uniformly to $x$.
This finishes the analysis of case 2.b and
hence finishes the  proof that $\pi_1(M)$ acts as
a convergence group on $\mi \cup \rr$.
\end{proof}

\section{Connections with Gromov hyperbolicity}

In this section we relate the flow ideal boundary and compactification
with the large scale geometry of $\mi$ and Gromov hyperbolic spaces.
Bowditch \cite{Bo1}, following ideas of Gromov,
gave a topological characterization of 
the action of a hyperbolic group on its ideal
boundary.

\begin{theorem}{(Bowditch \cite{Bo1})}{}
Suppose that $X$ is a perfect, metrisable compactum.
Suppose that a group $\Gamma$ acts on $X$,
such that the induced action on the space
of distinct triples is properly discontinuous
and cocompact. Then $\Gamma$ is a hyperbolic
group. Moreover there is a natural $\Gamma$-equivariant
homeomorphism of $X$ into $\partial \Gamma$, where 
$\partial \Gamma$ 
is the Gromov ideal boundary of $\Gamma$.
\end{theorem}

The $\Gamma$-equivariant homeomorphism $\alpha: X \rightarrow
\partial \Gamma$ satisfies: if $f$ is an element of $\Gamma$
and $a$ is the attracting fixed point of the action
of
$f$ in $X$, then $\alpha(a)$ is the attracting
fixed point of the action of $f$ in $\partial \Gamma$.
In our situation $X = \rr$ and $\Gamma = \pi_1(M)$,
which acts on $X$.

If $\pi_1(M^3)$ is Gromov hyperbolic, Gromov also showed
that $\mi$ has a compactification with an ideal boundary
\cite{Gr,Gh-Ha,CDP}. It is equivariantly homeomorphic
to the Gromov boundary of $\pi_1(M)$, which is
denoted by $\si$. 
The following is now an immediate consequence of theorem \ref{unifco}.

\begin{theorem}{}{}
Let $\Phi$ be a pseudo-Anosov flow without perfect fits and
not conjugate to a suspension Anosov flow.
Let $\rr$ be the flow ideal sphere.
Theorem \ref{unifco} shows  that $\pi_1(M^3)$ acts as a uniform
convergence group on $\rr$. Bowditch's theorem 
implies that $\pi_1(M)$ is Gromov hyperbolic
and the action of $\pi_1(M)$ on $\rr$ is topologically conjugate
to the action of $\pi_1(M)$ on the Gromov ideal
boundary $\si$ of $\mi$.
\label{rigid}
\end{theorem}

Let $\zeta: \rr \rightarrow \si$ be the conjugacy given
by theorem \ref{rigid}. It is uniquely defined.

In addition to theorem \ref{rigid}
we also prove that the group equivariant compactification $\mi \cup \rr$
is equivariantly homeomorphic to the 
Gromov compactification of $\mi$.
First we define a bijection

$$\xi: \mi \cup \rr \ \rightarrow \ \mi \cup \si \ \ \ - \ \ 
{\rm if} \ x \in \mi \ \ {\rm let} \ \xi(x) = x, \ \ \  \
{\rm if} \  x \in \rr \  \ {\rm let} \ \xi(x) = \zeta(x)$$

\noindent
Clearly this map $\xi$ is group equivariant:
if $\gamma$ is in $\pi_1(M)$ then $\xi(\gamma(x)) = \gamma(\xi(x))$.

\begin{theorem}{}{}
The map $\xi: \mi \cup \rr \rightarrow \mi \cup \si$
is a group equivariant homeomorphism.
The map $\varphi_1 = \xi \circ \varphi: \partial \oo \rightarrow \si$ is
a group invariant Peano curve.
\label{equi}
\end{theorem}

\begin{proof}{}
We only need to show that $\xi$ is a homeomorphism.
We know that $\mi$ is open in both $\mi \cup \rr$ and
in $\mi \cup \si$ and the induced topology 
from both of these is the original topology of
$\mi$.
Hence $\xi$ is continuous in $\mi$.
Let $x$ in $\rr$.
Lemma \ref{first} showed that
$\mi \cup \rr$ is first countable. Hence to check continuity of
$\xi$ at $x$  we only need to verify
what happens for sequences.
Let then $p_n$ in $\mi \cup \rr$ converging 
to $x$ as $n$ converges to infinity.
Theorem \ref{rigid} shows that $\xi$ restricted to $\rr$ 
is continuous. Hence we may assume that $p_n$ is
in $\mi$.
Then there are $q_n$ in a fixed compact set in 
$\mi$ and $\gamma_n$ in $\pi_1(M)$ with
$\gamma_n(q_n) = p_n$.
We may assume that the $\gamma_n$ are distinct otherwise
up to subsequence all $\gamma_n = \gamma$ and $\gamma_n$ sends
$q_n$ into a fixed compact set, contradiction.

By the convergence group action of 
$\pi_1(M)$ on $\mi \cup \rr$ (theorem \ref{bige}),
there is a source/sink pair $y,z$ for some subsequence
of $(\gamma_n)$ (still denoted $(\gamma_n)$). 
Since $\pi_1(M)$ also acts as a convergence
group on $\mi \cup \si$ \cite{Fr,Ge-Ma}, then
for this subsequence there is another
subsequence (denoted 
$(\gamma_{n_i})$) 
with a source/sink pair $b, a$ for the action in $\mi \cup \si$.
As the action of $\pi_1(M)$ on $\rr$ is equivariantly conjugate
to the action on $\si$, it follows that 
$\xi(y) = b$ and $\xi(z) = a$.
Now 

$$p_{n_i} \ = \ \gamma_{n_i}(q_{n_i}) \ \ {\rm converges \ to}
\ \  x  \ \ {\rm in}   \ \ \mi \cup \rr$$

\noindent
with $q_{n_i}$ in a fixed compact
set of $\mi$. It follows that $x$ is the 
sink of the sequence $(\gamma_{n_i})$  acting on $\mi \cup \rr$,
so
$x = z$.

Consider now the situation in $\mi \cup \si$.
Here $\xi(p_{n_i}) = \gamma_{n_i}(\xi(q_{n_i}))$ with $q_{n_i}$ 
in a compact set of $\mi$. 
Then $\xi(q_{n_i})$ is in a compact set of $\mi$.
By the convergence group property
of $\pi_1(M)$ acting on $\mi \cup \si$, then up to subsequence
we may assume that $\gamma_{n_i}(\xi(q_{n_i}))$ converges to 
the sink $a = \xi(z) = \xi(x)$.
This shows that for any sequence $(p_n)$ converging to
$x$ in $\mi \cup \rr$, there is a subsequence 
$(p_{n_i})_{i 
\in {\bf N}}$ with $\xi(p_{n_i})$ converging to $\xi(x)$
in $\mi \cup \si$.
It follows that $\xi$ is continuous at $x$ and so 
$\xi$ is continuous.
Since $\mi \cup \rr$ is compact and Hausdorff then
$\xi$ is a homeomorphism.

Using this fact the second statement follows from the fact that
the map $\varphi: \partial \oo \rightarrow \rr$ is group equivariant.
This finishes the proof of the theorem.
\end{proof}




\section{Quasigeodesic flows and quasi-isometric singular foliations}

In the last two sections of the article we obtain geometric
consequences for flows and foliations.
A flow $\Phi$ in a manifold $N$ is quasigeodesic if
in $\widetilde N$, distance along flow lines of $\wwp$ is
a bounded multiplicative distortion of ambient distance.
Quasigeodesic flows are extremely useful \cite{Th1,Gr,Ca-Th,Fe-Mo}.
In this section we show that if $\Phi$ is a 
pseudo-Anosov flow without perfect fits,
then $\Phi$ is
quasigeodesic. This will produce new examples
of quasigeodesic pseudo-Anosov flows.
A foliation ${\cal E}$ (singular or not) is quasi-isometric
if distance along leaves of $\widetilde {\cal E}$ is a bounded multiplicative
distortion of ambient distance in $\widetilde N$.
This property is very important \cite{Th1,Th2,Mor,Gr,Ca-Th,Fe5,Fe8}.
We show that the stable/unstable foliations of 
pseudo-Anosov flows without perfect fits are quasi-isometric.
These results are consequences of theorems 
\ref{rigid}, \ref{equi} and previous results.
Notice that both properties are invariant under quasi-isometries:
if $\Phi$ is a quasigeodesic flow and $\Phi'$ is topologically
conjugate to $\Phi$, then $\Phi'$ is also quasigeodesic.
The same holds for the quasi-isometric property for 
foliations. A quasi-isometry is a map so that when lifted
to the universal cover it is bi-lipschitz in the large.

\begin{theorem}{}{}
Let $\Phi$ be a pseudo-Anosov flow without perfect fits.
Then $\Phi$ is a quasigeodesic flow.
In addition the foliations $\Lambda^s, \Lambda^u$
are quasi-isometric foliations.
\label{qg}
\end{theorem}

\begin{proof}{}
Suppose first that $\Phi$ is topologically conjugate to a suspension
Anosov flow. If $\Phi'$ is a suspension Anosov flow and $M$ has
the solv metric, then $\wwp'$ is a flow by minimal geodesics
and the stable and unstable foliations $\wls', \wlu'$ are
foliations by totally geodesic surfaces.
Therefore $\Phi$ is quasigeodesic and $\wls, \wlu$ are quasi-isometric
foliations.

For the remainder of the proof assume that $\Phi$ is not
conjugate to a suspension Anosov flow.
Since $\Phi$ has no perfect fits,
theorem \ref{rigid} shows that $\pi_1(M)$ is Gromov hyperbolic.
We now show that $\Phi$ is quasigeodesic.
We will prove 3 topological properties of the flow lines
in $\mi \cup \rr$ (and then transfer them to $\mi \cup \si$):

\vskip .1in
\noindent
{\underline {Property 1}} $-$ For each flow line
$\alpha$ of $\wwp$ then it limits in 
a single point of $\rr$ denoted by $\alpha_+$
and similarly for the backwards direction.

$\alpha$ can be seen as a vertical segment
$\{ y \} \times (-1,1)$ in $\cd \times \mmp$
where $y$ is in $\oo$.
Let $q$ in $\alpha$.
Let $z = (y,1)$ and let
$x = \varphi(z)$ a point in $\rr$.
We claim that $x$ is the limit 
of $\alpha$ in $\mi \cup \rr$.
Let $g$ be the decomposition element of $\partial (\cd \times \mmp)$
associated to $z$.
For any neighborhood $U$ of $x$ in $\mi \cup \rr$
there is a set
type $A = A(l_1, ..., l_n, r_1, ..., r_n, \tau)$ (definition \ref{neigh})
with $A \subset \eta^{-1}(U)$.
The description of type (iii) in definition \ref{neigh}, shows that
since $z$ is in $g$ any such set $A$ as above contains
$\wwp_t(q)$ for all $t$ bigger than some $t_0$.
This shows that in $\mi \cup \rr$ the flow line
$\alpha$ forward converges to $x$.

Similarly let $\alpha_-$ be the negative ideal point
of $\alpha$. In fact for any $q$ in $\mi$ let 
$\alpha = \wwr(q)$ and define
$\mu_+(q) = \alpha_+$ and
$\mu_-(q)  = \alpha_-$.
This defines functions $\mu_+, \mu_-: \mi \rightarrow \rr$.

\vskip .1in
\noindent
{\underline {Property 2}} $-$ For each flow line
$\alpha$ of $\wwp$, then the ideal points $\alpha_+, \alpha_-$
are distinct.

Let $\alpha$ be an orbit of $\wwp$ which is $\{ y \} \times (-1,1)$
for some $y$ in $\partial \oo$.
Suppose that \ $(y,1), \ (y,-1)$ \ project to the same point in $\rr$.
By the construction of theorem \ref{flobo},
if a point in $\cd \times 1$ is identified
to a point in $\cd \times \{ -1 \}$ then at least
one of them has to be in $\partial \oo \times \mmp$.
Since $y$ is in $\oo$,
this is not the case here. Therefore  $\alpha_+, \alpha_-$ are
distinct in $\rr$.

\vskip .1in
\noindent
{\underline {Property 3}} $-$ The endpoint functions
$\mu_+, \mu_-: \mi \rightarrow \rr$ are continuous.

Given $p$ in $\mi$, $p$ is in $\{ y \} \times (-1,1)$ for some $y$
in $\oo$.
For any neighborhood $U$ of $\mu_+(p)$ then $\eta^{-1}(U)$
contains
a set of type $A(l_1, ..., l_n, r_1, ..., r_n, \tau)$.
By the description of neighborhoods in definition \ref{neigh},
then
for any $q$ sufficiently near $p$ then the forward orbit of $q$ is
eventually in $A(l_1, ..., l_n, r_1, ..., r_n, \tau)$
and so $\mu_+(q)$ is in $U$.
This shows continuity of the map $\eta_+$ at $x$.

\vskip .1in
Since the map $\xi: \mi \cup \rr \rightarrow \mi \cup \si$
is a homeomorphism then as seen in
$\mi \cup \si$ properties 1) through 3) also hold for
orbits of $\wwp$.
This is the key fact here: properties in $\mi \cup \rr$ get
transferred to $\mi \cup \si$.
We now use a result of Fenley-Mosher \cite{Fe-Mo}
which  states that if $\pi_1(M)$ is Gromov hyperbolic
and properties 1) through 3) hold for orbits 
of a flow $\wwp$ then $\Phi$ is a uniform quasigeodesic flow.
Hence $\Phi$ is a quasigeodesic flow.

We now prove that $\Lambda^s, \Lambda^u$ are quasi-isometric
singular foliations.
Given that $\Phi$ is a quasigeodesic pseudo-Anosov flow, then
it was proved in \cite{Fe5}, theorem 3.8, that $\Lambda^s$ is 
quasi-isometric if and only if $\wls$ has
Hausdorff leaf space and similarly for $\Lambda^u$.
Suppose that $\wls$ does not have Hausdorff leaf space
and let $F, L$ not separated in $\wls$.
Theorem \ref{theb} shows that $F, L$ are connected by
a chain of lozenges. A lozenge has 2 perfect
fits $-$ which are disallowed by hypothesis.
Hence $\Lambda^s, \Lambda^u$ are quasi-isometric foliations.
This finishes the proof of theorem \ref{qg}.
\end{proof}


\section{Asymptotic properties of foliations}

Here we show that $\rrrr$-covered foliations and
foliations with one sided branching in atoroidal manifolds 
are transverse to pseudo-Anosov flows without perfect
fits and therefore satisfy the continuous extension
property.
This parametrizes and characterizes their limit sets.
In addition this shows that pseudo-Anosov flows without perfect fits are 
very common.

\begin{theorem}{}{}
Let $\fol$ be a Reebless $\rrrr$-covered foliation in
$M^3$ closed, atoroidal and not finitely covered
by ${\bf S}^2 \times {\bf S}^1$.
Then $\pi_1(M)$ is Gromov hyperbolic and $\fol$
satisfies the continuous extension property.
This produces new examples of group invariant Peano curves.
\label{conti}
\end{theorem}

\begin{proof}{}
Up to a double cover, we may assume that $\fol$ is transversely
orientable.
Recall that $\rrrr$-covered means that the leaf space
of $\fn$ is homemorphic to the reals $\rrrr$.
If $\fol$ is $\rrrr$-covered and
$M$ is not finitely covered by ${\bf S}^2 \times {\bf S}^1$,
then it was proved in \cite{Fe6,Cal2} that
either there is a $\zz$ subgroup of $\pi_1(M)$ or 
there is a pseudo-Anosov flow $\Phi$ transverse to $\fol$
and regulating for $\fol$.
Since $M$ is (homotopically) atoroidal the second option occurs.
Regulating means that every orbit of $\wwp$ intersects an
arbitrary leaf of $\fn$ and vice versa.
Therefore the orbit space of $\wwp$ can be identified
to the set of points in a leaf $F$ of $\fn$.
Using Candel's theorem \cite{Ca}
we can assume that all leaves of $\fol$ are hyperbolic.
In this situation the set $\cd = \oo \cup \partial \oo$ is
naturally identified to the compactification of
$F$ with a circle at infinity $\pin F$.
Here is why:
The construction of $\Lambda^s, \Lambda^u$ in
\cite{Fe5,Cal2} is obtained by blowing down
2 transverse laminations which intersect the
leaves of $\fol$ in geodesics.
Therefore there are 2 geodesic laminations
(stable and unstable) in $F$, whose complementary regions are finite
sided ideal polygons \cite{Fe6,Cal2}.
It follows that the ideal points of $F$ are either ideal points of leaves of 
$\wls \cap F, \wlu \cap F$ or have neighborhood systems defined
by leaves of these. Hence $\pin F$ is naturally homeomorphic
to $\partial \oo$ and $F \cup \pin F$ is homeomorphic
to $\oo \cup \partial \oo$.
This works for any $F$ in $\fn$.

Suppose there is a perfect fit between a leaf $L$ of $\wls$ and
a leaf $H$ of $\wlu$.
Then in $\oo$ there are rays of $\Theta(L), \Theta(H)$ 
defining
the same ideal point in $\partial \oo$.
By the above description there is a pair of geodesics
in $F$, one stable and one unstable with the
same ideal point in $\pin F$.
By hyperbolic geometry considerations these 2 geodesics
are asymptotic in $F$, so projecting to $M$ and taking 
limits we obtain a leaf of $\fol$ so that there is a 
geodesic which is a leaf of both the stable and unstable
laminations. This contradicts the fact that the
stable and unstable laminations are transverse.

It follows that $\Phi$ has no perfect fits. By theorem
\ref{rigid} it follows that $\pi_1(M)$ is Gromov hyperbolic
(this particular
fact was already known, by the Gabai-Kazez theorem \cite{Ga-Ka}
and results in \cite{Fe6,Cal2,Fe7}).
By theorem  \ref{qg}
it folows that $\Phi$ is a quasigeodesic
pseudo-Anosov flow and in addition the map
$\varphi_1: \partial \oo \rightarrow \si$ is a group
equivariant Peano curve.
The previously known  examples of such group invariant Peano curves
occurred for fibrations \cite{Ca-Th} and slitherings by work
of Thurston \cite{Th5}. 
The results here are useful 
because 
Calegari \cite{Cal1}, 
showed that there are many examples of $\rrrr$-covered
foliations in hyperbolic $3$-manifolds which are
not slitherings or uniform foliations. 
The results here imply the previous results for fibrations
and slitherings.

Now we analyse the continuous extension property for the leaves
of $\fn$. Since $\Phi$ is quasigeodesic and transverse
to $\fol$, then 
the main theorem in \cite{Fe8} implies that leaves
of $\fn$ extend continuously to $\si$. Hence $\fol$
has the continuous extension property.
This finishes the proof of theorem \ref{conti}.
We remark that there is a direct proof of the continuous
extension property in this case since
$\partial \oo$ is naturally identified to $\pin F$.
For simplicity we just quote the result of \cite{Fe8}.
Notice that the leaves of $\fn$ have limit set the whole
sphere, so each leaf $F$ of $\fn$ produces a sphere
filling curve. 
\end{proof}

We now turn to foliations with one sided branching.

\begin{theorem}{}{}
Let $\fol$ be a Reebless foliation with 
one sided branching in $M^3$ closed, atoroidal
and not finitely covered by ${\bf S}^2 \times {\bf S}^1$.
Then $\pi_1(M)$ is Gromov hyperbolic.
There is a pseudo-Anosov flow $\Phi$ transverse
to $\fol$ which has no perfect fits and hence 
is a quasigeodesic flow and its
stable/unstable foliations are quasi-isometric.
It follows that $\fol$ has the continuous extension
property.
\end{theorem}

\begin{proof}{}
Recall that $\fol$ has one sided branching if the leaf space of
$\fn$ is not Hausdorff, but the non Hausdorff behavior 
occurs only in (say) the negative direction.
Since $\fol$ has one sided branching it is transversely
oriented.
Suppose that $\fn$ has branching only in the negative
direction.
When $M$ is atoroidal and not finitely covered
by ${\bf S}^2 \times {\bf S}^1$, Calegari \cite{Cal3} produced
a pseudo-Anosov flow $\Phi$ which is transverse
to $\fol$ and forward regulating for $\fol$.
Forward regulating means that if $x$ is in
a leaf $F$ of $\fn$ and $L$ is a leaf of $\fn$,
for which there is a positive transversal
from $F$ to $L$, then the forward orbit of $x$ 
intersects $L$.

As in the $\rrrr$-covered case this is obtained from
2 laminations transverse to $\fol$ which intersect
the leaves of $\fn$ in a collection of geodesics.
Suppose there is $G$ in $\wls$ and $H$ in $\wlu$ 
forming a perfect fit.
Then $G$ intersects $F_0$ leaf of $\fn$
 and $H$ intersects $F_1$.
Since $\fol$ has one sided branching there is a leaf
$F$ of $\fn$ with positive transversals from
$F_0$ to $F$ and from $F_1$ to $F$.
By the above property $G$ and $H$ intersect $F$.
There are rays in $\Theta(G)$ and $\Theta(H)$ with
same ideal point $p$ in $\partial \oo$.

The ideal circle of $\oo$ is the same as the universal
circle for the foliation $\fol$ in this case \cite{Cal3}.
The universal circle is obtained as the inverse limit
of circles at infinity escaping in the positive direction.
Given $A, B$ leaves of $\fn$ we write $A < B$ if there
is a positive transversal from $A$ to $B$.
Given $A < B$ in $\fn$ then there is a dense set of 
directions in $A$ which are asymptotic to $B$ \cite{Cal3}.
This is not symmetric $-$ there is not a dense
set of directions from $B$ which is asymptotic to $A$.
In our situation with $F_0 < F$ and $F_1 < F$ then the 
asymptotic directions  from $F$ to $F_0$ form an
unlinked set with the asymptotic directions
from $F$ to $F_1$ \cite{Cal1,Cal3}.
This implies there are natural surjective, continuous,
weakly circularly monotone maps 
$\pin F \rightarrow \pin F_i$. The universal circle
${\cal V}$ is obtained as an inverse limit of these
maps.

The stable/unstable laminations are obtained by analysing
the action of $\pi_1(M)$ in ${\cal V}$ and producing
laminations $-$ that is, a collection
of pairs of points in ${\cal V}$ 
which are unlinked.
They produce a collection of geodesics in leaves
of $\fn$, without transverse intersections,
 which vary continuously in the transversal
direction.
Therefore if rays of $\Theta(G), \Theta(H)$ define the
same ideal point in $\partial \oo$, then in the leaf $F$ which they
jointly intersect the following happens:
the associated stable/unstable geodesics
are asymptotic.
As in the $\rrrr$-covered case this leads to a contradiction
to $\wls, \wlu$ being transverse.

We conclude that $\Phi$ has no perfect fits.
From this point on the proof follows the same arguments as 
in the $\rrrr$-covered
case.
\end{proof}

\begin{corollary}{}{}
Let $\fol$ be a Reebless foliation with one sided branching
in $M^3$ atoroidal and not finitely covered by ${\bf S}^2 
\times {\bf S}^1$.
For any leaf $F$ of $\fn$, then the limit set
of $F$ is not the whole sphere $\si$.
\end{corollary}

\begin{proof}{}
The limit set of a set $B$ in $\mi$ is the set of accumulation
points of $B$ in $\si$.
Suppose that there is branching of $\fn$ only in the negative direction.
Choose $E, L$  non separated from each other and
so that $F < E$.
Branching in the negative direction means that there is
a sequence of leaves $(G_n)$ on the {\underline {positive}}
side of $E, L$ which converges to both $E, L$.
Since $E, L$ are non separated from each other, then they
do not intersect the same orbit of $\wwp$. Recall the projection
map $\Theta: \mi \rightarrow \oo$.
The sets $\Theta(E), \Theta(L)$ are disjoint.
Since $E, L$ are non separated from each other
in their positive sides,
then the analysis in section 4 of \cite{Fe8}, 
shows that there is a slice leaf $S$ of a unstable leaf of $\wls$, so that
$s = \Theta(S)$ is a boundary component of $\Theta(L)$ and $s$
separates $\Theta(L)$ from $\Theta(E)$.

Then the limit set of $S$, $\Lambda_S$ is a Jordan curve $C$ $-$
this is shown in   \cite{Fe1,Fe5}.
This uses the fact that $\ls$ is a quasi-isometric foliation.
The construction implies that the leaf $E$ separates $F$ from $S$ $-$
here we use that $E, L$ are non separated from 
each other on their positive sides and $F$ is in the back of $E$.
Since $S$ is disjoint from $F$ then the limit set of $F$ is
contained in the closure of one complementary component
of $\Lambda_S$ in $\si$. Therefore $\Lambda_F$ is not $\si$.
\end{proof}

\noindent
{\bf {Remarks}} $-$  

1) The remaining open situation for the continuous extension property
is that of $\fol$ with two sided branching. This means that
the leaf space of $\fn$ has non Hausdorff behavior in both
the positive and negative directions. The particular case of finite
depth foliations was recently solved in \cite{Fe8} using completely
different methods than this article. In particular in \cite{Fe8} one starts
with strong geometric properties, namely that $M$ is
hyperbolic and there is a leaf which is quasi-isometrically
embedded (the compact leaf) and this has enormous
geometric consequences. The tools here are purely from
dynamical systems.

2) Many $\rrrr$-covered examples in hyperbolic $3$-manifolds
which are not slitherings were constructed
by Calegari in \cite{Cal1}. Many explicit examples of 
foliations with one sided braching were constructed
by Meigniez in \cite{Me}.


3) Suppose that $\fol$ is Reebless in $M^3$ with
$\pi_1(M)$ negatively curved.
It is asked in \cite{Fe3,Fe8}: is $\fol$ \ $\rrrr$-covered
if and only if for some $F$ in $\fn$ then the limit
set $\Lambda_F = \si$?
If $\fol$ is $\rrrr$-covered then $\Lambda_F = \si$ for
every $F$ in $\fn$ \cite{Fe3}. The converse is
true if there is a compact leaf  in $\fol$ \cite{Go-Sh,Fe3}.
The previous theorem shows that if $\fol$ has one
sided branching then $\Lambda_F$ is not $\si$ for
any $F$ in $\fn$.
Therefore the remaining  open case for this question
 is also when $\fol$ has 2 sided branching.

4) The results of this article show that foliations
in manifolds with Gromov hyperbolic fundamental
group are very similar to surface Kleinian groups:
the $\rrrr$-covered case corresponds to doubly degenerate
surface Kleinian groups, where the limit sets are the whole sphere.
The foliations with one sided branching correspond to
singly degenerate Kleinian groups where there is a single
component of the domain of discontinuity.
It remains to be seen whether foliations with 2 sided
branching behave like non degenerate surface
Kleinian groups.

{\footnotesize
{
\setlength{\baselineskip}{0.01cm}

\noindent
Florida State University

\noindent
Tallahassee, FL 32306-4510, USA

and 

\noindent
Princeton University

\noindent
Princeton, NJ 08544-1000, USA

}
}

\end{document}